\documentclass[twoside, 11pt]{amsart}
\usepackage[active]{srcltx} 
\usepackage[all]{xy}
\CompileMatrices
\usepackage{latexsym}
\usepackage{ifthen}
\usepackage{amsfonts}
\usepackage{amssymb}
\usepackage{amsthm}
\usepackage{url}
 \newtheorem{thm}{Theorem}[section]
 \newtheorem{prop}[thm]{Proposition}
 \newtheorem{cor}[thm]{Corollary}
 \newtheorem{lem}[thm]{Lemma}
  \newtheorem{ques}[thm]{Question}
 
 \newtheorem*{bthm}{Theorem}
\newtheorem*{bconj}{Conjecture}
\newenvironment{ithm}{\bigskip\noindent{\bf Theorem}\it}{\bigskip}

\theoremstyle{definition}
\newtheorem{defn}[thm]{Definition}

\theoremstyle{remark}
\newtheorem{rem}[thm]{Remark}

\font\russ=wncyr10 \font\Russ=wncyr10   scaled\magstep 1
\def\sha{\hbox{\russ\char88}}
\def\Sha{\hbox{\Russ\char88}}

\newcommand{\fp}{\ifmmode {\mathbb{F}_p}\else$\mathbb{F}_p$\ \fi}
\newcommand{\zp}{\ifmmode {\mathbb{Z}_p}\else$\mathbb{Z}_p$\ \fi}
\newcommand{\z}{\mathbb{Z}}
\newcommand{\zpMod}{\ifmmode\mbox{$\zp$-Mod}\else$\zp$-Mod \fi}
\newcommand{\Mod}{\ifmmode\mbox{$\Lambda$-Mod}\else$\Lambda$-Mod \fi}
\renewcommand{\mod}{\ifmmode\mbox{$\Lambda$-mod}\else$\Lambda$-mod
\fi}
\newcommand{\La}{\ifmmode\Lambda\else$\Lambda$\fi}

\newcommand{\Hom}{{\mathrm{Hom}}}
\newcommand{\Ext}{{\mathrm{Ext}}}
\newcommand{\Tor}{{\mathrm{Tor}}}
\newcommand{\tor}{{\mathrm{tor}}}
\newcommand{\E}{{\mathrm{E}}}

\newcommand{\pd}{{\mathrm{pd}}}
\newcommand{\cd}{{\mathrm{cd}}}
\newcommand{\rk}{{\mathrm{rk}}}

\newcommand{\D}{\mathbf{D}}

\renewcommand{\H}{\mathrm{H}}
\newcommand{\du}{\mathrm{D}}
\renewcommand{\t}[1]{\mbox{$\mbox{\rm T}_{#1}$}}

\newcommand{\M}{\ifmmode {\frak M}\else${\frak M}$ \fi}
\newcommand{\m}{\ifmmode {\frak m}\else$\frak m$ \fi}
\newcommand{\p}{\ifmmode {\frak p}\else${\frak p}$\ \fi}
\renewcommand{\P}{\ifmmode {\frak P}\else${\frak P}$\ \fi}
\newcommand{\e}{\ifmmode {\mathcal{E}}\else$\mathcal{E}$ \fi}

\newcommand{\h}{{\mathcal{ H}}}
\newcommand{\C}{\mbox{$\mathcal{ C}$}}
\renewcommand{\O}{\mbox{$\mathcal{ O}$}}
\newcommand{\G}{\ifmmode {\mathcal{G}}\else${\mathcal{G}}$\ \fi}
\renewcommand{\d}{\ifmmode {\mathcal{ D}}\else${\mathcal{D}}$\ \fi}
\newcommand{\A}{\ifmmode {\mathcal{A}}\else${\mathcal{ A}}$\ \fi}
\newcommand{\la}{\mathbb{A}} 
\newcommand{\lu}{\mathbb{U}} 
\newcommand{\gas}{\mathbb{A}_S}  
\newcommand{\gla}{\mathbb{A}_S}  
\newcommand{\glu}{\mathbb{U}_S}  
\newcommand{\gus}{\mathbb{E}_S} 
\newcommand{\gu}{\mathbb{E}} 
\newcommand{\tab}{\mathbb{T}}
\newcommand{\F}{\mbox{$\mathcal{F}$}}
\renewcommand{\projlim}[1] {{\lim\limits_{\stackrel{\displaystyle
\longleftarrow}{#1}}}}
\newcommand{\projlimsc}[1]
{{\lim\limits_{\stackrel{\scriptstyle \longleftarrow}{#1}}}}
\newcommand{\projlimssc}[1]
{{\lim\limits_{\scriptstyle \longleftarrow}}_{ #1}}
\newcommand{\dirlim}[1]
{{\lim\limits_{\stackrel{\displaystyle \longrightarrow}{#1}}}}
\newcommand{\dirlimssc}[1]
{{\lim\limits_{\scriptstyle \longrightarrow}}_{ #1}}

\renewcommand{\in}{\ \epsilon\ }

\newcommand{\kl}{[\![}
\newcommand{\kr}{]\!]}
\newcommand{\Qp}{\ifmmode {{\Bbb Q}_p}\else${\Bbb Q}_p$\ \fi}
\newcommand{\qp}{\ifmmode {{\Bbb Q}_p}\else${\Bbb Q}_p$\ \fi}
\newcommand{\Q}{\ifmmode {\Bbb Q}\else${\Bbb Q}$\ \fi}

\newcommand{\Coind}{\mbox{\mbox{\rm Coind}}}
\newcommand{\Ind}{\mathrm{Ind}}
\newcommand{\ho}{\mathrm{Ho}}
\newcommand{\Sel}{\mathrm{Sel}}
\newcommand{\csel}{\widehat{Sel}}
\newcommand{\ts}{\mbox{$\Sha$}}

\begin{document}

\title{On the Iwasawa theory of $p$-adic Lie Extensions}%
\author{Otmar Venjakob}%
\address{Universit\"{a}t Heidelberg\\ Mathematisches Institut\\
 Im Neuenheimer Feld 288\\
 69120 Heidelberg, Germany.}
\email{otmar@mathi.uni-heidelberg.de}
\urladdr{http://www.mathi.uni-heidelberg.de/\textasciitilde
otmar/}



\begin{abstract}
In this paper the new techniques and results concerning the
structure theory of modules over non-commutative Iwasawa algebras
 (\cite{css},\cite{ven1})
are applied to arithmetic:  we study Iwasawa modules over $p$-adic
Lie extensions $k_\infty$ of number fields $k$ ``up to
pseudo-isomorphism".  In particular, a close relationship is
revealed  between  the Selmer group of abelian varieties,  the
Galois group of the maximal abelian unramified $p$-extension of
$k_\infty$ as well as the Galois group
 of the maximal abelian   $p$-extension unramified outside $S$
where $S$ is a certain finite set of  places of $k.$ Moreover, we
determine the Galois module structure of  local units and other
modules arising from Galois cohomology.
\end{abstract}

\maketitle

\PUSH{intro.tex}%
\section{Introduction}

The starting point of the {\em Iwasawa theory of (non-commutative)
$p$-adic Lie groups} was M. Harris' thesis \cite{harris} in 1979.
For an elliptic curve $E$ over a number field $k$ without complex
multiplication he studied the Selmer group $\Sel(E,k_\infty)$ over
the extension $k_\infty=k(E(p))$ which arises by adjoining the
$p$-division points of $E$ to $k.$ Then, the Galois group
$G=G(k_\infty/k)$ is an open subgroup of $Gl_2(\zp)$ - due to a
celebrated theorem of Serre \cite{serre72} - and so a (compact)
$p$-adic Lie group. Following Iwasawa's general idea, he studied
the Pontryagin dual $\Sel(E,k_\infty)^\vee$ of the Selmer group as
module over  the Iwasawa algebra
 $$\La(G)=\zp\kl G\kr, $$
i.e.\ the completed group algebra of $G$ with coefficients in $\zp.$\\
In the late $90s$ J. Coates and S. Howson (\cite{coates-howsonI},
\cite{coates-howsonII}, \cite{coates}, \cite{howson}) as well as
Y. Ochi \cite{ochi} revived this Iwasawa-theoretic approach. Among
other things they proved a remarkable Euler characteristic formula
for the Selmer group, studied ranks, torsion-properties and
projective dimensions of
standard local and global Iwasawa modules.\\
The contributions of this work to the Iwasawa theory of  $p$-adic
Lie groups are obtained by applying some new techniques we have
developed in \cite{ven1}.  There we introduced the concept of
pseudo-null modules over $\La=\La(G),$ which is based on a general
dimension theory for Auslander regular rings (for the definition
see subsection \ref{dimensiontheory} and note that \La\ is a
non-commutative ring in general). Therefore it is fundamental for
our applications that $\La(G)$ is an Auslander regular ring if $G$
is a compact $p$-adic Lie group without $p$-torsion (cf.\
\cite[Thm. 3.26]{ven1}). As a first example in which this approach
proves  effective we consider the following generalization of a
theorem of R. Greenberg \cite{green} and T. Nguyen-Quang-Do
\cite{Ng} (who considered the case  $G\cong \mathbb{Z}_p^d$): For
a finite set $S$ of places of a number field $k$ let $k_\infty|k$
be a Galois extension unramified outside $S$ such that the Galois
group $G=G(k_\infty/k)$ is a torsion-free $p$-adic Lie-group and
let $k_S$ be the maximal outside $S$ unramified extension of $k.$
Then there is a basic result on the structure of the  Galois group
 $$X_S=G(k_S/k_\infty)^{ab}(p)$$
of the maximal abelian $p$-extension of $k_\infty$ unramified
outside $S$ considered as $\La(G)$-module, which is by definition
the maximal abelian pro-$p$ quotient of the Galois group
$G_S(k_\infty):=G(k_S/k_\infty).$

\begin{ithm} (Theorem \ref{nospeudonullglobalX})
If\/ $\H^2(G_S(k_\infty),\qp/\zp)=0,$ then the $\La(G)$-module
$X_S$ does not contain any non-trivial pseudo-null submodule.
\end{ithm}



Once having available the concept of pseudo-null modules one is
tempted to study Iwasawa modules ``up to pseudo-isomorphism". We
will write $M\sim N$ if there exists a \La-homomorphism $M\to N$
whose  kernel and cokernel is pseudo-null. Since in general $\sim$
is not a symmetric relation we consider also the quotient category
$\mod/{\mathcal{PN}}$ with respect to subcategory $\mathcal{PN}$
of pseudo-null \La-modules, which is a Serre subcategory, i.e.\
closed under subobjects, quotients and extensions.\\
Now it turns out that - as in the classical $\zp$-extension case -
the  $\La(G)$-module $X_S$ is closely related to the modules
$X_{nr}$ and  $X_{cs}^S$ which denote the Galois groups of the
maximal abelian unramified pro-$p$-extension of $k_\infty$ and the
maximal abelian unramified pro-$p$-extension of $k_\infty$  in
which every prime above $S$ is completely decomposed,
respectively. In the next theorem, $G_v$ denotes the decomposition
group of $G$ at a place $\nu,$ $S_f$ is the set of finite places
in $S,$ $\E^1(M)$  the Iwasawa adjoint $\Ext_\La^1(M,\La)$ of a
\La-module $M $ and $M(-1)$ means the twist of $M$ with the Galois
module $\zp(-1):=\Hom(\mu,\qp/\zp).$

\begin{ithm}  (Theorem \ref{XcspseudotorX})
If  $\mu_{p^\infty}\subseteq k_\infty,$ and $\dim(G_\nu)\geq 2$
for all $\nu\in S_f,$ then
 \[X_{nr}(-1)\sim
 X_{cs}^S(-1)\sim\E^1(\tor_\La
 X_S).\]
If, in addition, $G\cong \mathbb{Z}^r_p,$ $r\geq2,$ then even the
following holds:
 $$X_{nr}(-1)\sim X_{cs}^S(-1)\sim (\tor_\La X_S)^\circ,$$
 where $^\circ$ means that  $G$ operates via the involution
 $g\mapsto g^{-1}.$
 \end{ithm}

In this context we should mention that it is still an open
question - even for $G\cong \mathbb{Z}^r_p,$ $r\geq2$ - whether in
general $X_{nr}$ is pseudo-isomorphic to the dual
$(Cl(k_\infty)(p)^\vee)^\circ$ of the direct limit of the
$p$-primary ideal class groups in a $p$-adic tower of number
fields with involution $-^\circ$ (which can also be defined for
non-commutative $p$-adic Lie groups under additional assumptions,
see prop.\ \ref{circ})
\[X_{nr}\sim
(Cl(k_\infty)(p)^\vee)^\circ \;\;?\] In the $\zp$-case this is a
well-known theorem due to Iwasawa, see \ref{quescl-xnr} for
further discussion.

Drawing our attention to cohomology groups associated with an
abelian varieties $\A$  defined over $k,$ we set
$k_\infty=k(\A(p))$ and mention that $\H^1(G_S(k_\infty),
\A(p))^\vee$ has no non-zero pseudo-null submodule (Theorem
\ref{nospeudonullglobalXAbvariety}). With respect to the
($p$-)Selmer group $\Sel(\A,k_\infty)$ of $\A$ over
$k_\infty=k(\A(p))$  we generalize a result of P.\ Billot in the
case of good, {\em supersingular} reduction, i.e.\
$\widetilde{\A_{k_\nu}}(p)=0,$ at any place dividing $p.$ Over a
$\zp$-extension an analogous statement was proved by K. Wingberg
\cite[cor.\ 2.5]{wing}. We shall write $\A^d$ for the dual abelian
variety of $\A.$ Assume that $G(k_\infty/k)$ is a
 pro-$p$-group without any $p$-torsion. Then the
following holds (corollary \ref{supersingular-ab-var}):
 $$X_{cs}\otimes_{\zp}(\A^d(p))^\vee\sim\E^1(\tor_\La(\Sel(\A,k_\infty)^\vee)).$$
We also refer the reader to our joint work with Y. Ochi
\cite{ochi-ven} where we prove under certain conditions that the
Pontryagin dual of the Selmer group of an elliptic curve without
CM and good {\em ordinary} reduction at any place dividing $p$
does not contain any non-zero pseudo-null \La-submodule.

Furthermore,  we proved a structure theory for the $\zp$-torsion
part of a \La-module $M$ in \cite{ven1}. Up to pseudo-isomorphism
any $\zp$-torsion \La-module is of the form \[\bigoplus_i
\La/p^{n_i}.\] In particular, we obtained a natural  definition of
the $\mu$-invariant \[\mu(M):=\sum_i n_i\] of $M.$ Defining
$\mu(M):=\mu(\tor_\zp M)$ for an arbitrary \La-module $M,$ this
invariant   is additive on short exact sequences of \La-torsion
modules. Hence, we can formulate and prove a generalization of
theorem 11.3.7 of \cite{nsw}:

\begin{ithm} (Theorem \ref{free-prop})
Let $k_\infty|k$ be a $p$-adic pro-$p$ Lie extension such that $G$
is without $p$-torsion.  Then $\G=G(k_S(p)/k_\infty)$ is a free
pro-$p$-group if and only if $\mu(X_S)=0$ and the weak Leopoldt
conjecture holds, i.e.\ $\H^2(G_S(k_\infty),\qp/\zp)=0.$
\end{ithm}

In theorem \ref{base-change}  we describe how the weak Leopoldt
conjecture and the vanishing of $\mu(X_S)$ - if considered
simultaneously -  behave under change of the base field.
Furthermore,  we get a formula for the $\mu$-invariants for
different $S.$

We briefly outline further results. In  section \ref{local} we
generalize Wintenberger's result on the Galois module structure of
local units. Let $k$ be a finite extension of $\qp$ and assume
that $k_\infty|k$ is a Galois extension  with Galois group
$G\cong\Gamma\rtimes_\rho\Delta,$ where $\Gamma$
 is a pro-$p$ Lie group  of dimension $2$ (e.g.\ $\Gamma=\mathbb{Z}_p\rtimes\zp$)  and $\Delta$ is a profinite group of
{\em possibly infinite} order  prime to $p,$ which acts on
$\Gamma$ via $\rho:\Delta\to Aut(\Gamma).$ Then we characterize
the $\La(G)$-module structure of the Galois group
$G_{k_\infty}^{ab}(p)=G(k_\infty(p)/k_\infty),$ where
$k_\infty(p)$ is  the maximal abelian $p$-extension of $k_\infty,$
see theorem \ref{localAd2}.\\
Then we apply these results to the local study of elliptic curves
$E$ {\em with CM,} i.e.\  we determine the structure of local
cohomology groups with certain division points
of $E$ as coefficients.\\
Section \ref{global} is devoted to the study of ``global" Iwasawa
modules. Besides the themes already mentioned above  we study the
norm-coherent $S$-units of $k_\infty$
$$\gus:=\projlim{k'} ({\mathcal{ O}}_{k',S}^\times\otimes\zp)$$ by
means of Jannsen's spectral sequence for Iwasawa adjoints. Using
Kummer theory, we compare $\gus$ to
 $${\mathcal{E}}_S(k_\infty):=(E_S(k_
 \infty)\otimes_\z\qp/\zp)^\vee,$$
where $E_S(k_\infty)=\dirlimssc{k'} E_S(k')$ denotes the (discrete
module of) $S$-units of $k_\infty.$ In particular, we show that
$\E^0(\gus)\cong\E^0\E^0(\e_S(k_\infty)),$  where $\E^0(M)$
denotes $\Hom_\La(M,\La)$ for any \La-module $M,$ and thus
$$\rk_\La\gus=\rk_\La\e_S=r_2(k)$$ under some assumptions, see
corollary \ref{E0}. If $\E^0(\gus)$ is  projective, its structure
can be described more precisely. A criterion which tells us when
this is the case is given in proposition \ref{E0projectiv}.\\

\textsc{General Notation and Conventions}

We follow the notation in the paper \cite{ven1}, which is similar
to that used in \cite{nsw}. In particular, this means:

\begin{enumerate}

\item   For a  discrete  (resp.\ compact) $\zp$-module $N$
with continuous action by some profinite group $G$,
\[N^\vee=\Hom_{\zp,cont}(N,\Q_p/\zp)\] is the compact (resp.\
discrete) Pontryagin dual of $N$ with its natural $G$-action. If
$N$ is  $p$-divisible, \[T_p(N)=\Hom(\Q_p/\zp,N)= \projlim{i}\;
{_{p^i}}N\] denotes the Tate module of
 $N,$ where ${_{p^i}}N$ denotes the kernel of the multiplication by
 $p^i.$ For $G=G_k$ the  absolute Galois group of number or local field
 $k,$ we define the $r$th Tate twist   of $N$ by
 $N(r):=N\otimes_\zp T_p(\mu)^{\otimes
 r}$ for $r\in\mathbb{N}$ and $N(r):=N\otimes_\zp \Hom(T_p(\mu)^{\otimes
 r},\zp)$ for $-r\in \mathbb{N},$ where $\mu$ denotes the $G_k$-module of all
 roots of unity and by convention
 $T_p(\mu)^{\otimes 0}=\zp$ with trivial $G$-action. Finally,
we set \[N^\ast:=\varinjlim_i\Hom({_{p^i}}N,
 \mu_{p^\infty})(=T_p(N)^\vee(1)).\]

\item For a finitely generated abelian $p$-primary group $A$ we
denote by $A_{div}$ the quotient of $A$ by its maximal
$p$-divisible subgroup.

\item Let $G$ be a profinite group and $H$ a closed subgroup of
$G$. For a  $\Lambda(H)$-module $M$, we define $\Ind_{H}^G M :=
M\widehat{\otimes}_{\Lambda(H)}\Lambda(G)$
 (compact or completed induction), where $\widehat{\otimes}$ denotes completed
 tensor product, and  $\Coind_{H}^G M := \Hom_{\Lambda(H)}(
 \Lambda(G),M)$ (co-induction).

\item If $G$ is any profinite group, by $G(p)$ and $G^{ab}$
we denote the maximal pro-$p$ quotient and the maximal abelian
quotient $G/[G, G]$ of $G,$ respectively. For an abelian group $A$
we also denote by $A(p)$ its $p$-primary component.

\item Let $k$ be a field. For a $G_k$-module $A$, we write
$A(k):=H^0(G_k, A)$.

\item By a Noetherian ring, we mean a left and right Noetherian ring
(with  a multiplicative unit). By $pd_\Lambda(M)$ we denote the
projective dimension of $M$ while $pd(\Lambda)$ denotes the global
dimension of $\Lambda.$

\item The dual of an abelian variety $\A$ is denoted by
$\A^d.$

\item We refer the reader to subsection \ref{diagram} for
the definitions of $R^{ab}(p),$ $ N^{ab}(p),$ $ X,$ $Y,$ $J$ and $
Z.$
\end{enumerate}

\textsc{Acknowledgements.}\\

This article is based on  a part of my Dissertation, Heidelberg
2000, and I would like to thank my supervisor Kay Wingberg most
warmly for leading me to the nice field of ``higher dimensional"
Iwasawa theory. John Coates receives my special gratitude for his
great interest and valuable suggestions. I wish to express my warm
thanks to Yoshi Ochi for our fruitful interchange of ideas when he
stayed in Heidelberg in 1999/2000. Susan Howson is  heartily
acknowledged for providing me with a copy of her Ph.D. thesis and
for the discussions in Heidelberg as well as in Cambridge.
Finally, I thank the referee for pointing out some inaccuracies
and for his or her comments which helped to improve the
exposition.

%
\POP
\PUSH{basics.tex}%

\section{Algebraic properties of \La-modules}

\subsection{Notation and Preliminaries}\label{not} We  recall some basic facts on $p$-adic
Lie groups and their Iwasawa algebras which are thoroughly
discussed in \cite{ven1}; the reader who is not familiar  with
them is recommended to read first or parallel sections 1-3 of
(loc.\ cit.).
 For any compact $p$-adic Lie group $G$
the completed group algebra $\La=\La(G)$ is Noetherian (see
\cite{la}V 2.2.4). If, in addition, $G$ is pro-$p$ and has no
element of finite order, e.g.\ if $G$  is uniform, then $\La(G)$
is an integral domain, i.e.\ the only zero-divisor in $\La(G)$ is
$0$
 (see \cite{Ne}); literally the same proof shows that the corresponding statements
  hold also for the completed group algebra $\fp\kl
G\kr $ with coefficients in the finite field $\fp$ with $p$
elements. For instance, for $p\geq n+2,$ the group $Gl_n(\zp)$ has
no elements of order $p,$ in particular, $GL_2(\zp)$ contains no
elements of finite $p$-power order if $p\geq 5$ (see \cite{howson}
4.7). In any case, the normal subgroup
$\Gamma_i:=\ker(Gl_n(\zp)\to Gl_n(\z/p^i))$ of $Gl_n(\zp)$ is a
uniform pro-$p$ group for $i\geq 1$ if $p\neq 2$ or $i\geq 2$ if
$p=2$ by \cite[thm.\ 5.2]{dsms}. We should also mention that $G$
has finite cohomological dimension $\cd_p G=m$ if and only if $G$
has no element of finite $p$-power order and its dimension as
$p$-adic analytic manifold equals $m.$

 If \La\ is  Noetherian and without zero-divisors we can form a skew
field $Q(G)$ of fractions of \La\ (see \cite{GoWa}). This allows
us to define the rank of a \La-module:

\begin{defn}
The rank $\rk_\Lambda M$ is defined to be the dimension of
\linebreak $M\otimes_\La Q(G)$ as a left vector space over $Q(G)$
 $$\rk_\Lambda M=\dim_{Q(G)}(M\otimes_\La Q(G)).$$

\end{defn}

Obviously, the rank is finite for any $M$ in the category $\mod$
of finitely generated $\La$-modules. For the rest of this section,
we assume that {\em all \La-modules considered are finitely
generated.}

 By $\ho(\La)$ we denote the category of ``\La-modules up to
 homotopy"and we write $M\simeq
N,$ if $M$ and $N$ are homotopy equivalent,  i.e.\ isomorphic in
$\ho(\La),$ which holds if and only if $M\oplus P\cong N\oplus Q$
with projective \La-modules $P$ and $Q.$ In particular,
$M\simeq 0$ if and only if $M$ is projective.\\

For $M\in\mod$ we define the Iwasawa adjoints of $M$ to be
 $$\E^i(M):=\Ext^i_\La(M,\La),\quad i\geq 0,$$
which are a priori right \La-modules by functoriality and the
right \La-structure of the bi-module \La\ but will be considered
as left modules via the involution of \La. By convention we set
$\E^i(M)=0$ for $i<0.$ The \La-dual $\E^0(M)$ will also be denoted
by $M^+.$

It can be shown that for $i\geq 1$ the functor $\E^i$ factors
through $\ho(\La)$ defining  a functor
 $$\E^i:\ho(\La)\to\mod.$$

By $\du$ we denote the  transpose
 $$\du :\ho(\La)\to\ho(\La),$$
 which is a contravariant duality functor, i.e.\ it satisfies
  $\du^2=\mathrm{Id}.$ Furthermore, if $\pd_\La M\leq 1$
then $\du M\simeq \E^1(M).$  The next property will be of
particular importance:

\begin{prop}(cf.\ \cite[prop.\ 5.4.9]{nsw})\label{canonicalsequ}
For $M\in\mod$ there is a canonical exact sequence
 $$\xymatrix@1{
    {\ 0 \ar[r] } &  {\ \E^1\du M \ar[r] } &  {\ M \ar[r]^{\phi_M }} &  {\ M^{++} \ar[r] } &  {\
    \E^2\du M\ar[r] }& 0, \
 }$$
where $\phi_M $ is the canonical map from $M$ to its bi-dual. In
the following we will refer to the sequence as ``the" canonical
sequence (of homotopy theory).
\end{prop}

A \La-module $M$ is called {\em reflexive} if $\phi_M $ is an
isomorphism from $M$ to its bi-dual $M\cong M^{++}.$

As Auslander and Bridger \cite{aus} suggest the module $\E^1\du M$
should be considered as torsion submodule of $M.$ Indeed, if \La\
is a Noetherian integral domain this submodule  is a torsion
module while $M^{++}$ is torsion-free and thus $\E^1\du M$
coincides exactly with the set \footnote{A priori it is not clear
whether this sets forms a submodule if \La\ is not commutative.}
of torsion elements $\tor_\La M.$ Hence, a \La-module $M$ is
called {\em \La-torsion module} if $\phi_M\equiv0,$ i.e.\ if
$\tor_\La M:=\E^1\du M=M.$ We say that $M$ is {\em
\La-torsion-free} if $\E^1\du M=0.$ It turns out that
 a finitely generated \La-module $M$ is  a \La-torsion module if
and only if $M$ is  a $\La(G')$-torsion module (in the strict
sense) for some open pro-$p$ subgroup $G'\subseteq G$ such that
$\La(G')$ is integral. Since $M^{++}$ embeds into a free
\La-module  the torsion-free \La-modules are exactly the
submodules of free modules (see \cite[before prop.\ 2.7]{ven1} for
details).

 Sometimes it is also
convenient to have the notation of the $1^{st}$ {\em syzygy} or
{\em loop space} functor $\Omega:\mod\to\ho(\La)$ which is defined
as follows (see \cite[1.5]{ja-is}): Choose a surjection $P\to M$
with $P$ projective. Then $\Omega M$ is defined by the exact
sequence
 $$\xymatrix@1{
    {\ 0 \ar[r] } &  {\ \Omega M \ar[r] } &  {\ P \ar[r] } &  {\ M \ar[r] } &  {\
    0 }. \
 }$$

\subsection{Dimension theory for the Auslander regular ring
\La(G)} \label{dimensiontheory} Let $G$  be any compact $p$-adic
group without $p$-torsion. In \cite{ven1} we proved that
$\La=\La(G)$ is an {\em Auslander regular ring}, i.e.\ \La\ has
finite projective dimension $d:=\pd \La=\cd_pG+1$ (by a result of
Brumer) and satisfies the Auslander condition: For any \La-module
$M$, any integer $m$ and any submodule $N$ of $\E^m(M),$ the grade
of $N$ satisfies $j(N)\geq m.$ Recall that the {\em grade} $j(N)$
is the smallest number $i$ such that $\E^i(N)\neq
0$ holds.\\
 Therefore there is a nice dimension theory for \La-modules which
we will recall briefly (for proofs and further references see
\cite{ven1}). A priori, any $M\in\mod$ comes equipped with a
finite filtration
$$ \t0 (M) \subseteq  \t1 (M)\subseteq\cdots\subseteq\t{d-1}
(M)\subseteq\t{d} (M)=M.$$ If we call the number
$\delta:=min\{i\mid \t{i} (M)=M\}$ the {\em dimension} $\delta(M)$
then $\t{i}(M)$ is just the maximal submodule of $M$ with
$\delta$-dimension less or equal to $i.$ We should mention that
for abelian $G$ the dimension $\delta(M)$ coincides with the Krull
dimension of  $supp_\La(M).$

 The filtration is related to the Iwasawa adjoints via a spectral
 sequence, in particular we have
 \[\t{i} (M)/\t{i-1} (M)  \subseteq \E^{d-i}\E^{d-i}(M) \]
and either of these two terms is zero if and only the other is.
Furthermore, the equality $\delta(M)+j(M)=d$ holds for any $M\neq
0.$

Note that $M$ is a $\La$-torsion module if and only if its
codimension $codim(M):= d-\delta(M)$ is greater or  equal to $1.$

A \La-module $M$ is called {\em pseudo-null} if its codimension
$codim(M)$ is greater or equal to $2.$ As in the commutative case
we say that a homomorphism $\varphi:M\to N$ of \La-modules is a
{\em pseudo-isomorphism} if its kernel and cokernel are
pseudo-null. A module $M$ is by definition pseudo-isomorphic to a
module $N,$ denoted
 $$M\sim N,$$
if and only if there exists a pseudo-isomorphism from $M$ to $N.$
In general, $\sim$ is not symmetric even in the $\zp$-case. While
in the commutative case $\sim$ is symmetric at least for torsion
modules, we do not know whether this property still holds in the
general case.\\
If we want to reverse pseudo-isomorphisms, we have to consider the
quotient category \linebreak $\mod/{\mathcal{PN}}$ with respect to
subcategory $\mathcal{PN}$ of pseudo-null \La-modules, which is a
Serre subcategory, i.e.\ closed under subobjects, quotients and
extensions. By definition, this quotient category is the
localization $({\mathcal{ PI}})^{-1}\mod$ of \mod\ with respect to
the multiplicative system ${\mathcal{ PI}}$ consisting of all
pseudo-isomorphisms. Since \mod\ is well-powered, i.e.\ the family
of submodules of any module $M\in\mod$ forms a
set, 
these localization exists, is an abelian category and the
universal functor $q:\mod\to\mod/{\mathcal{PN}}$ is exact.
Furthermore, $q(M)=0$ in $\mod/{\mathcal{PN}}$ if and only if
$M\in\mathcal{PN}.$  Recall that a morphism $h:q(M)\to q(N)$ in
the quotient category can be represented for instance  by two
\La-homomorphisms $f:M'\to M$ and $g:M'\to N$ where $f$ is a
pseudo-isomorphism and such that $h\circ q(f) =q(g);$ it is an
isomorphism if and only if $g$ is a pseudo-isomorphism. If there
exists an isomorphism between $q(M)$ and $q(N)$ in the quotient
category we also write $M\equiv N \;\mathrm{mod} \;\mathcal{PN}.$

Note that for any pseudo-isomorphism $f:M\to N$ the induced
homomorphism $\E^1(f)$ is a pseudo-isomorphism, too. If $M, N$ are
\La-torsion modules, also the converse statement holds. By the
universal property of the localisation, we obtain a functor
\[ \E^1:\mod/{\mathcal{PN}}\to \mod/{\mathcal{PN}},\]
which is exact if it is restricted to the full subcategory
$\mbox{\La-mod}^{\geq 1}/{\mathcal{PN}}$ of\linebreak
$\mod/{\mathcal{PN}}$ consisting of all \La-modules of codimension
greater or equal to $1,$ i.e.\ \La-torsion modules. More
precisely, there is a natural isomorphism of functors:
\[\E^1\circ\E^1\cong \mathrm{Id}:\mbox{\La-mod}^{\geq 1}/{\mathcal{PN}}\to \mbox{\La-mod}^{\geq 1}/{\mathcal{PN}}.\]

It is known that any torsion-free module $M$ embeds into a
reflexive module with pseudo-null cokernel while any torsion
module $M$ is pseudo-isomorphic to $\E^1\E^1(M)$ (cf.\
\cite[Prop.3.13]{ven1}). Moreover, there is a canonical
pseudo-isomorphism $\E^1(M)\sim\E^1(\tor_\La M)$ for any
\La-module $M.$

By $\mod(p)$ we shall write the plain subcategory of $\mod$
consisting of $\zp$-torsion modules while by
${\mathcal{PN}}(p)``={\mathcal{PN}}\cap\mod(p)"$ we denote the
Serre subcategory of $\mod(p)$ the objects of which are
pseudo-null \La-modules. In other words $M$ belongs to ${\mathcal{
PN}}(p)$ if and only if it is a $\La/p^n$-module for an
appropriate $n$ such that $\E^0_{\La/p^n}(M)=0.$ Recall that there
is a canonical exact functor
$q:\mod(p)\to\mod(p)/{\mathcal{PN}}(p).$ Then, there is the
following structure theorem on the $\zp$-torsion part of a
finitely generated \La-module:

\begin{thm}(cf.\ \cite[Thm.\ 3.40]{ven1})\label{p-structure}
Assume that $G$ is  a $p$-adic analytic pro-$p$ group without
$p$-torsion. Let $M$ be in $\mod(p).$ Then there exist uniquely
determined natural numbers $n_1,\ldots,n_r$ and an isomorphism in
$\mod(p)/{\mathcal{PN}}(p)$
 $$M\equiv \bigoplus_{1\leq i\leq r} \La/p^{n_i}\ \mbox{\rm  mod
 }{\mathcal{ PN}}(p).$$
\end{thm}

We define the $\mu$-invariant of a \La-module $M$ as
\[\mu(M)=\sum_i n_i(\tor_\zp M),\]
where the $n_i=n_i(\tor_\zp M)$ are determined uniquely by the
structure theorem applied to $\tor_\zp M.$ This invariant is
additive on short exact sequences of \La-torsion modules and
stable under pseudo-isomorphisms. Alternatively, it can be
described as
\[\mu(M)=\rk_{\fp\kl G\kr }\bigoplus_{i\geq 0}
{_{p^{i+1}}M}/{_{p^{i}}M}=rk_{\fp\kl G\kr }\bigoplus_{i\geq 0}
{{p^{i}}\tor_\zp M}/{{p^{i+1}}\tor_\zp M}.\]

Very recently, J. Coates, R.Sujatha and P. Schneider \cite{css}
found a general structure theorem for \La-torsion modules. They
proved that any finitely generated $\La(G)$-torsion module
decomposes into the direct sum of cyclic modules up to
pseudo-isomorphism, i.e.\ in the quotient category
$\La\mathrm{-mod}^{\geq 1}/{\mathcal{PN}}.$

\begin{bthm}[Coates-Schneider-Sujatha]
Let $G$ be a  $p$-valued compact $p$-adic analytic group. Then,
for any finitely generated $\La(G)$-torsion module $M$ there exist
finitely many reflexive left  ideals $J_1,\ldots,J_r$ and an
injective $\La(G)$-homomorphism

\[\bigoplus_{1\leq i\leq r}\La/J_i\hookrightarrow M/M_{ps}\]
with pseudo-null cokernel, where $M_{ps}=T_{\dim(G)-2}(M)$ denotes
the maximal pseudo-null submodule of $M.$ In particular, it holds
 \[M\equiv \bigoplus_{1\leq i\leq r}\La/J_i \ \mbox{\rm  mod
 }{\mathcal{ PN}.} \]
\end{bthm}

For the precise definition of a $p$-valued compact Lie group see
(loc.\ cit.) or directly in Lazard's article \cite{la}; we just
want to mention that any uniform pro-$p$-group belongs to this
class of pro-$p$ Lie groups, which is stable under taking closed
subgroups.

 It is still not known whether the ideals $J_i$ can be chosen as principal ideals as in the commutative case.
 Anyway, if we restrict to this kind of modules, we can define a second involution
 \[^\circ:\La\mathrm{-mod}^{\geq 1}_{pr}/{\mathcal{PN}}\to \La\mathrm{-mod}^{\geq 1}_{pr}/{\mathcal{PN}}\]
on the full subcategory $\La\mathrm{-mod}^{\geq
1}_{pr}/{\mathcal{PN}}$ of $\La\mathrm{-mod}^{\geq
1}/{\mathcal{PN}}$ consisting of those objects which are
isomorphic (in the quotient category) to a direct sum of cyclic
modules of the form $\La/\La f,$ $f\in\La.$  For any such $f$ we
set $(\La/\La f)^\circ:=\La/\La f^\circ,$  where
$^\circ:\La\to\La$ also denotes the involution of the group
algebra (induced by $g\mapsto g^{-1}$). The following proposition
implies among other things that this definition is invariant under
pseudo-isomorphism and therefore it extends to the whole category
$\La\mathrm{-mod}^{\geq 1}_{pr}/{\mathcal{PN}}.$

\begin{prop}\label{circ} Let $G$ be a  profinite group such that $\La=\La(G)$ is a Noetherian integral
ring.Then the following holds:
\begin{enumerate}
\item For any $f\in \La$ there is an isomorphism $\E^1(\La/\La f)\cong \La/\La
f^\circ.$
\item Assuming that $G$ is a $p$-adic analytic group without
$p$-torsion the above two involutions coincide: \[-^\circ\cong
\E^1(-):\La\mathrm{-mod}^{\geq 1}_{pr}/{\mathcal{PN}}\to
\La\mathrm{-mod}^{\geq 1}_{pr}/{\mathcal{PN}}\]
\end{enumerate}
\end{prop}

The proof is standard, see for example the proof of proposition
\ref{contragredient}, where we denote the involution on \La\ by
$\iota.$

We conclude this section with a technical result which will be
needed in the arithmetic applications.
 \begin{prop}\label{reflexivepd3}
Let $\Lambda$ be an Auslander regular ring. For any
$\Lambda$-module $M$ such that $\pd_\Lambda \E^0(M)\leq 1$ (e.g.
if $\pd\ \Lambda=3$ or if $\pd\ \Lambda=4$ and $\E^4\E^1(M)=0$)
its double dual $\E^0\E^0(M)$  is a $2$-syzygy of $\E^1\E^0(M),$
i.e.\ there is an exact sequence
 $$\xymatrix@1{
    {\ 0 \ar[r] } &  {\E^0\E^0(M) \ar[r] } &  {P_0 \ar[r] } &  {P_1 \ar[r] } &  {\E^1\E^0(M)\ar[r] } & 0\
 }$$
with projective modules $P_0$ and $P_1.$ Furthermore, in the case
of\/ $\pd\ \Lambda=3$ or $4,$ it holds that
$\E^1\E^0(M)\cong\E^3\E^1(M).$ If, in addition, $M$ itself is
reflexive and $\pd\ \Lambda =3,$ then
$\E^3\E^1M\cong\E^1(M)^\vee.$
\end{prop}

\begin{proof}
First observe that $\E^0(M)$ is a $2$-syzygy of $\mathrm{D}(M)$
due to the definition of the latter module, i.e.\ $\pd_\Lambda
\E^0(M)\leq\pd\ \Lambda-2=1,$ if $\pd\ \Lambda=3.$ In the case of
$\pd\ \Lambda=4$ it holds $\E^3\E^0(M)=\E^4\E^0(M)=0$ and
$\E^2\E^0(M)\cong\E^4\E^1(M)$ due to Bj\"{o}rk's spectral sequence
(see \cite[3.1]{ven1}). Hence, if $\E^4\E^0(M)$ vanishes, it
follows that $\pd_\Lambda \E^0(M)\leq 1$.  Now, choosing a
projective resolution of $E^0(M)$
 $$\xymatrix@1{
    {\ 0 \ar[r] } &  {\ \E^0(P_1) \ar[r] } &  {\ \E^0(P_0)  \ar[r] } &  {\ \E^0(M) \ar[r] } &  {\
    0 ,} \
 }$$
we derive the exact sequence
 $$\xymatrix@1{
    {\ 0 \ar[r] } &  {\ \E^0\E^0(M) \ar[r] } &  {\ P_0 \ar[r] } &  {\ P_1 \ar[r] } &  {\
    \E^1\E^0(M)\ar[r] } & 0. \
 }$$
But $\E^1\E^0(M)\cong\E^3\E^1(M)$ due  to Bj\"{o}rk's spectral
sequence for $\pd\ \Lambda\leq 4.$ If $M$ itself is reflexive and
$\pd\ \Lambda=3,$ then $\E^1\E^1(M)=\E^2\E^1(M)=0,$ i.e.\
$\E^1(M)$ is finite, respectively $\E^3\E^1(M)\cong\E^1(M)^\vee.$
\end{proof}

\subsection{Some representation theory}

 In the following lemma we shall write
$I(\Gamma)$ for the kernel of the canonical map $\zp\kl G\kr
\to\zp\kl G/\Gamma\kr ,$ where $\Gamma$ is any closed normal
subgroup of the profinite group $G.$ By $\mathrm{Rad}_G$ we denote
the radical of $\zp\kl G\kr ,$ i.e.\ the intersection of all open
maximal left (right) ideals of $\zp\kl G\kr .$ Finally, we write
$$M_G=M/I_GM$$ for the module of coinvariants of $M$ and  $\H_\bullet(G,M)$  for the
$G$-homology of a compact \La-module $M,$ which can be defined as
left derived functor of $-_G$ or alternatively as
$\Tor^\La_\bullet(\zp,M),$ where $\Tor$ denotes the left derived
functor of the complete tensor product $-\widehat{\otimes}_\La-.$

\begin{lem}\label{free-asympt}
Let $G=\Gamma\rtimes\Delta$ be the semi-direct product of a
uniform pro-$p$-group $\Gamma$ of dimension $t$ and a finite group
$\Delta$ of order $k$ prime to $p$.  If we write
$U_n=\Gamma^{p^n}\trianglelefteq G,$ then for any compact
$\Lambda=\Lambda(G)$-module $M$, the following statements are
equivalent:

\begin{enumerate}
\item $M\cong\Lambda^d,$
\item $M_\Gamma\cong\zp[\Delta]^d$ as $\zp[\Delta]$-module and for all $n$ $$\rk_\zp M_{U_n}=\rk_\zp \zp[G/U_n]^d=dkp^{tn},$$
\item $M_\Gamma/p\cong\mathbb{F}_p[\Delta]^d$ as $\mathbb{F}_p[\Delta]$-module
and for all $n$ $$\log_p \# (M_{U_n}/p^n)=\log_p \#
(\z/p^n[G/U_n]^d)={ndkp^{tn}}.$$
\end{enumerate}
\end{lem}

\begin{proof}
Obviously, (i) implies (ii) and (iii). For the converse let us
first assume that (ii) holds and let $m_1,\ldots , m_d\in M$ be
lifts of a $\zp[\Delta]$-basis of $M_\Gamma.$ Then the map
$\phi:\bigoplus_{i=1}^d \Lambda e_i\to M,$ which sends $e_i$ to
$m_i,$ is surjective, because $I(\Gamma)\subseteq \mathrm{Rad}_G$
(compare to the proof of \cite{nsw}. 5.2.14 (i), $d\Rightarrow b$)
and therefore we can apply  Nakayama's lemma \cite{nsw}, 5.2.16
(ii), (with $\mathrm{Rad}_G$ instead of $\M$). Hence, the induced
maps $\phi_{U_n}:\bigoplus_{i=1}^d \Lambda(G/U_n) e_i\to M_{U_n},$
are surjective, too. But since both modules have the same
$\zp$-rank by assumption, these maps are  isomorphisms and (i)
follows. The implication $\mathrm{(iii)}\Rightarrow\mathrm{ (i)}$
is proved analogously noting that $p\Lambda+I(\Gamma)\subseteq
\mathrm{Rad}_G.$
\end{proof}

For a finite group $G$ we denote by $K_0(\qp[G])=K'_0(\qp[G])$ the
Grothen\-dieck group of finitely generated $\qp[G]$-modules (which
are  projective by Maschke's theorem). If $G$ is a profinite group
and $U\trianglelefteq G$ an open normal subgroup we define the
Euler characteristic $h_U(M)$ of a finitely generated
$\Lambda=\Lambda(G)$-module $M$ to be the class
 $$h_U(M):=\sum (-1)^i [\H_i(U,M)\otimes_\zp\qp]\in K_0(\qp[G/U]).$$
Before stating the next result we recall some facts about the
representation theory of finite groups. So let $\Delta$ be a
finite group  of order  $n$ prime to $p.$ Then,  there is a
decomposition  $$\zp[\Delta]\cong\prod\zp[\Delta]e_i,\
e_i=\frac{n_i}{n}\sum_{g\in \Delta} \chi_i(g^{-1})g$$ of
$\zp[\Delta]$ in ``simple" components (in the sense that  they are
simple algebras after tensoring with $\qp$). If
$G=\Gamma\times\Delta,$ this induces a decomposition of
$\Lambda=\prod\Lambda^{e_i},$ $\Lambda^{e_i}=\zp\kl \Gamma\kr
[\Delta]e_i$ into a product of rings. Here $\{\chi_i\}$ is the set
of irreducible $\qp$ characters ($\widehat{=}$
$\mathbb{F}_p$-characters because $n$ is prime to $p$)  of
$\Delta$ and $n_i$ are certain natural numbers associated with
$\chi_i$ (see below). The simple algebras $\qp[\Delta]e_i$
decompose into the direct sum of their simple left ideals which
all belong to the same isomorphism class, say $N_i,$ i.e.\ there
is a isomorphism of $\qp[\Delta]$-modules
 $$\qp[\Delta]e_i\cong N_i^{n_i}.$$
In particular, $n_i$ is the length of $\qp[\Delta]e_i$ and can be
expressed as \[n_i=\chi(e_i)(\dim_\qp N_i)^{-1},\] where $\chi$ is
the character of the left regular representation of $\qp[\Delta].$

Now let $G$ be again a $p$-adic Lie group and set $\La:=\La(G).$
Recall that a finitely generated \La-module $M$ is a \La-torsion
module if and only if $M$ is  a $\La(G')$-torsion module for some
open pro-$p$ subgroup $G'\subseteq G$ such that $\La(G')$ is
integral.

 \begin{prop}\label{euler}
Let $G=\Gamma\times\Delta$ be the product of a pro-$p$ Lie group
$\Gamma$ of finite cohomological dimension $\cd_p(\Gamma)=m$ and a
finite group $\Delta$ of order $n$ prime to $p$ and let
$U\trianglelefteq\Gamma$ be an open normal subgroup. Then, for any
finitely generated
$\Lambda$-torsion module $M,$ 
 it holds
 $$h_U(M)=0.$$
 \end{prop}

 \begin{rem}
 For semi-direct products this statement is false in general. For
 example, it is easy to see that for $G=\zp\rtimes_\omega\Delta$
 with non-trivial $\omega$ the Euler characteristic of $\zp$ is not zero:
 $h_U(\zp)=[\qp]-[\qp(\omega)]\neq 0.$
 \end{rem}

\begin{proof}(of prop.\ \ref{euler})
We claim that under the assumptions of the theorem $M$ possesses a
finite free resolution. Indeed, since the Noetherian ring
$\Lambda$ has finite global dimension $\pd\Lambda=m+1$, there is
always a resolution of the form
 $$\xymatrix@1{
    {\ 0 \ar[r] } &  {\ P \ar[r] } &  {\ \Lambda^{d_m} \ar[r] } &  {\ \cdots \ar[r] } &  {\
    \Lambda^{d_0}\ar[r] } &0, \
 }$$
 with a projective  module $P.$
Since $M^{e_i}$ is a $\Lambda(\Gamma)$-torsion module (it is even
$\Lambda(\Gamma')$-torsion!) and  since $P^{e_i}$ is a free
$\Lambda(\Gamma)$-module, it must hold that
$P^{e_i}\cong(\Lambda(\Gamma))^{k_i d_{m+1}}$ as
$\Lambda(\Gamma)$-modules, where $k_i=\chi(e_i)$ denotes the
$\zp$-rank of $\zp[\Delta]e_i$  and $d_{m+1}=\sum_{i=0}^m (-1)^i
d_{m-i}.$ Consequently, $P^{e_i}_\Gamma\cong\z_p^{k_id_{m+1}}$ as
$\zp$-modules, respectively
$P^{e_i}_\Gamma\otimes\qp\cong\mathbb{Q}_p^{k_id_{m+1}}$ as
$\qp$-modules holds. But $P^{e_i}_\Gamma\otimes\qp$ must be
isomorphic to the direct sum of $m$ copies of $N_i$ for some $m$
due to the semi-simplicity of $\qp[\Delta].$ Counting
$\qp$-dimensions, we obtain $m=n_id_{m+1}$ and hence
$P^{e_i}_\Gamma\otimes\qp\cong\qp[\Delta]e_i^{d_{m+1}}.$ Since
$P^{e_i}_\Gamma$ is a projective $\zp[\Delta]$-module, this
implies $P^{e_i}_\Gamma\cong\zp[\Delta]e_i^{d_{m+1}},$
respectively $P^{e_i}\cong\Lambda(G)e_i^{d_{m+1}}$ (compare to the
proof of lemma \ref{free-asympt}) and
$P\cong\Lambda(G)^{d_{m+1}}.$
  This  proves the claim.\\
Furthermore, we observe that $\sum (-1)^id_i=0$
 and denote the resolution by \linebreak $F^\bullet\to M.$ Using the fact
 that the Euler characteristic of a bounded complex equals the
 Euler characteristic of its homology, 
 we calculate
  \begin{eqnarray*}
 \sum (-1)^i [\H_i(U,M)\otimes_\zp\qp]&=&\sum (-1)^i[\H_i(F^\bullet\otimes_\Lambda\qp[G/U])]\\
                  &=&\sum (-1)^i[F^\bullet\otimes_\Lambda\qp[G/U\kr \\
                  &=&\sum (-1)^i[\qp[G/U]^{d_i}]\\
                  &=&(\sum (-1)^id_i)[\qp[G/U\kr =0.
  \end{eqnarray*}
\end{proof}

\begin{lem}\label{kuz}
Let $G$ be a  profinite group, $H\subseteq G$ a closed subgroup
and $U\unlhd G$ an open normal subgroup. Then for any compact
$\zp\kl H\kr $-module $M$ the following is true:
\begin{enumerate}
\item $(\Ind^H_G(M))_U\cong\Ind^{HU/U}_{G/U}(M_{U\cap H})$ and
\item $\H_i(U,(\Ind^H_G(M)))\cong\Ind^{HU/U}_{G/U}\H_i(U\cap
H,M)$ for all $i\geq 0.$
\end{enumerate}
\end{lem}

\begin{proof}
The dual statement of (i) is proved in \cite{kuz} while (ii)
follows from (i) by homological algebra.
\end{proof}

\begin{lem}
Let $G=\Gamma\times\Delta$ be the product of a pro-$p$-group
$\Gamma$ and a finite group $\Delta$ of order prime to $p$. Then,
for any $\Lambda=\zp\kl \Gamma\kr [\Delta]$-module $M$ and for any
irreducible character $\chi$ of $\Delta$ with values in $\qp,$ the
following is true:
\begin{enumerate}
\item
$\Hom_\Lambda(M^{e_\chi},\Lambda)\cong\Hom_\Lambda(M,\Lambda)^{e_{\chi^{-1}}},$
\item $\mathrm{E}_\La^i(M^{e_\chi})\cong\mathrm{E}_\La^i(M)^{e_{\chi^{-1}}}$ for any $i\geq
0.$
\end{enumerate}
\end{lem}

\begin{proof}
While (ii) is a consequence of (i) by homological algebra the
first statement can be verified at once using the decompositions
$M\cong\bigoplus M^{e_\chi}$ and
$\Lambda\cong\bigoplus\Lambda^{e_\chi}:$
\begin{eqnarray*}
\Hom_\Lambda(M,\Lambda)^{e_{\chi^{-1}}}&\cong&\Hom_\Lambda(M,\Lambda^{e_\chi})\\
                                       &\cong&\Hom_\Lambda(M^{e_\chi},\Lambda^{e_\chi})\\
                                       &\cong&\Hom_\Lambda(M^{e_\chi},\Lambda).
\end{eqnarray*}
\end{proof}
\POP
\PUSH{diagram.tex}%
\subsection{Modules associated with group presentations}
\label{diagram}

Let $\C$ be a  class of finite groups closed under taking
subgroups, homomorphic  images and group extensions. Given an
exact sequence of  pro-$\C$-groups
 $$1\rightarrow \h\rightarrow \G\rightarrow G\rightarrow 1,$$
where $\G$ is assumed to be finitely generated, we choose a
presentation $\F\twoheadrightarrow\G$ of  $\G$ by a free
pro-$\C$-group $\F_d$  of rank $d$ and we associate the following
commutative diagram to it:
\begin{equation}\label{defR}
\xymatrix{ &1 &1 & & \\
              1\ar[r] &{\h}\ar[r]\ar[u] &{\G}\ar[r]\ar[u] &G\ar[r] &1 \\
              1\ar[r] &R\ar[r]\ar[u] &{\F_d}\ar[r]\ar[u] &G\ar[r]\ar@{=}[u] &1 \\
               &N\ar[u]\ar@{=}[r] &N\ar[u] & & \\
               &1\ar[u] &1\ar[u] & & \\  }
\end{equation}
Here, $R$ and $N$ are defined by the exactness of the
corresponding sequences. In general, the {\em $p$-relation module
$N^{ab}(p)$ of $\G$ with respect to the chosen free presentation}
 (and similarly $R^{ab}(p)$
   with respect to $G$ instead of $\G$)  fits into the following
exact sequence, which
 is called Fox-Lyndon resolution associated with the above free
 representation of $\G:$

\begin{equation}\label{lyndon}
\xymatrix@1{
   {\ 0 \ar[r] } &  {N^{ab}(p) \ar[r] } &  {\Lambda(\G)^d \ar[r] } &  {\Lambda(\G) \ar[r] } &  {\
   \zp\ar[r] }&0.
}
\end{equation}
\noindent Hence, if $\cd_p(\G)\leq 2$, then $N^{ab}(p)$ is a
projective $\Lambda(\G)$-module.

Furthermore, the augmentation ideal $I_{\F_d}$, i.e.\ the kernel
of $\Lambda(\F_d)\to\zp$, is a free $\Lambda(\F_d)$-modules of
rank $d$:
 $$
 I_{\F_d}\cong \Lambda(\F_d)^d $$
(for a proof of these facts, see \cite{nsw} Chap V.6).

 Let $A$ be a $p$-divisible $p$-torsion abelian group of finite $\zp$-corank $r$ with a continuous
 action of $\G.$

 \begin{defn}
 For a finitely generated $\La=\La(\G)$-module $M$ we define the finitely generated
 \La-module $$M[A]:=M\otimes_{\zp}A^{\vee}=\Hom_{cont.,\zp}(M,A)^{\vee}$$ with diagonal
 $\G$-action.We shall also write $M(\rho)$ for this
 $r$-dimensional twist where $\rho:G\rightarrow Gl_r(\zp)$ denotes
 the operation of $G$ on $A^\vee$.
 \end{defn}

 Note that the functor $-[A]$ is exact and that $\La[A]$ is a free \La-module of rank $r$ (cf. \cite[lem 4.2]{ochi-ven}).





\begin{prop}\label{contragredient} For every $i\geq0$
$$\E^i(M(\rho))\cong\E^i(M)(\rho^d), $$ where $\rho^d$ is the
contragredient representation, i.e.\ $\rho^d(g)=\rho(g^{-1})^t$ is
the transpose matrix of $\rho(g^{-1}).$
\end{prop}

\begin{proof}
By homological algebra (and using a free presentation of $M$) it
suffices to prove the case $i=0$ for free modules. Finally, we
only have to check the commutativity of the following diagram
which is associated to an arbitrary homomorphism $\phi:\La\to\La$
 $$\xymatrix{
   {\Hom_\La(\La(\rho),\La)}\ar[r]^{\phi(\rho)^\ast}\ar[d]&    {\Hom_\La(\La(\rho),\La)}\ar[d] \\
   {\La^r}\ar[d]     &  {\La^r}\ar[d] \\
   {\Hom_\La(\La,\La)(\rho)}\ar[r]^{\phi^\ast (\rho^d)}  &
   {\Hom_\La(\La,\La)(\rho)}.}$$
First note that via the  identification
 $\xymatrix@1{{\La^r}\ar@{=}[r]^{\psi_\rho} & {\La(\rho)} }$ the
 matrix representing $\phi(\rho)$ is $A:=\sum a_g g\rho(g^{-1}),$
 where we assume for simplicity that $\phi(1)=:a=\sum a_g
 g\in\zp[G].$ We denote by $\iota$ both, the involution
 $\La\to\La,\ g\mapsto\ g^{-1}$ (also extended to matrices with coefficients in \La) and the isomorphism of
 left \La-modules $\La\to\Hom_\La(\La,\La), g\mapsto (1\mapsto
 g^{-1}).$ Then its easy to see that the following two  diagrams
 commute
  $$\xymatrix{
  {\Hom_\La(\La(\rho),\La)}\ar[r]^{\phi(\rho)^\ast}\ar[d]^{(\psi_\rho)^\ast}&    {\Hom_\La(\La(\rho),\La)}\ar[d]^{(\psi_\rho)^\ast} \\
  {\Hom_\La(\La^r,\La)}\ar[d]^{i^r}&      {\Hom_\La(\La^r,\La)}\ar[d]^{i^r} \\
  {\La^r\ar[r]^{\iota(A^t)}} & {\La^r} , }$$

  $$\xymatrix{
     {\La^r\ar[r]^{B}}\ar[d]^{\psi_{\rho^d}} & {\La^r}  \ar[d]^{\psi_{\rho^d}}  \\
      {\La(\rho)\ar[r]^{\iota(a)(\rho^d)}}\ar[d]^{i(\rho^d)} &       {\La(\rho)}\ar[d]^{i(\rho^d)} \\
      {\Hom_\La(\La,\La)(\rho)}\ar[r]^{\phi^\ast (\rho^d)}  &
   {\Hom_\La(\La,\La)(\rho)},}$$
where $B=\sum a_g g^{-1}\rho^d(g),$ because $\iota(a)=\sum
a_gg^{-1}.$ We are done if we can verify $B=\iota(A^t).$ But
\begin{eqnarray*}
\iota(A^t)&=&\sum a_g g^{-1}\rho(g^{-1})^t\\
          &=&\sum a_g g^{-1}\rho^d(g)=B.
\end{eqnarray*}
\end{proof}

With the notation
 \begin{eqnarray} X&:=&X_{\h,A}:=\H^1(\h ,A)^{\vee}\label{X}\\
                   Y&:=&Y_{\h,A}:=(I_{\G}[A])_{\h}\label{Y}\\
                   J&:=&J_{\h,A}:={\rm ker}( \La (\G)[A]_{\h}\rightarrow
                   (A^{\vee})_{\h}),\label{J}
\end{eqnarray}
we get the following proposition:

\begin{prop}(\cite[lem.\ 4.5]{ochi-ven})\label{diagram-prop}
We have a commutative and exact diagram

$$\xymatrix{
    &   &   & 0 & 0 &   \\
    &   &   & J\ar[u] & J\ar[u] &   \\
  0\ar[r] & {\H^2(\h, A)^\vee} \ar[r] & (N^{ab}(p)[A])_{\h}\ar[r] & {\Lambda(G)^{dr}}\ar[u]\ar[r] & Y\ar[u]\ar[r] & 0 \\
  0\ar[r] & {\H^2(\h, A)^\vee}\ar[r]\ar@{=}[u] & (\H^1(N, A)^{\h})^\vee\ar[r]\ar@{=}[u] & {\H^1(R, A)^\vee}\ar[u]\ar[r]  & X\ar[u]\ar[r] & 0. \\
   & & & 0\ar[u] & 0\ar[u] &
}$$
Furthermore, if $\cd_p(\G)\leq 2$, then $N^{ab}(p)[A]$ is a
projective $\La(\G)$-module and $(N^{ab}(p)[A])_{\h}$ a projective
$\La(G)$-module.

\end{prop}



\begin{rem}
Assume $A$ is trivial as $\h$-module. Then the above diagram can
be easily  obtained by  twisting Jannsen's original diagram (i.e.\
with coefficients $\Qp/\zp$):
$\mathrm{diagram}(A)=\mathrm{diagram}(\qp/\zp)[A].$ Also the
higher Iwasawa adjoints of the occurring modules can be calculated
via proposition \ref{contragredient}: \begin{eqnarray*}
\E^i(X_{\h,A})&\cong&\E^i(X_{\h,\qp/\zp})(\rho^d)\\
\E^i(Y_{\h,A})&\cong&\E^i(Y_{\h,\qp/\zp})(\rho^d)\\
              &\cdots&
\end{eqnarray*}
\end{rem}

 The following theorem is a consequence of the diagram.
The restriction to $p$-adic Lie groups without $p$-torsion is
necessary in order to apply the dimension theory developed in
\cite{ven1}.

\begin{thm}\label{nospeudonullXA} Let $\cd_p(\G)\leq 2$ and $G$ a $p$-adic Lie group of dimension $h$ without
$p$-torsion. If the ``weak Leopoldt conjecture holds for $A$ and
$\h$", i.e.\ if $\H^2(\h,A)=0$, then neither $Y$ nor $X$ have
non-zero pseudo-null submodules: $\t{h-1}(X)=\t{h-1}(Y)=0$.
\end{thm}

\begin{proof}
Apply proposition 3.10 of \cite{ven1} to $Y$, which has
$\pd(Y)\leq 1$ according to the above diagram, and note that
$\t{h-1}(X)\subseteq\t{h-1}(Y)$ by proposition 3.2 of (loc.\
cit.).
\end{proof}

Let \begin{equation}\label{Z}
Z=Z_{\h,A}:=(D_2^{(p)}(\G,A)^\h)^\vee
\end{equation} where
$$D_2^{(p)}(\G,A)=\dirlim {U \subseteq_o
\G,n}(H^2(U,{_{p^n}A}))^\vee$$ and the direct limit is taken with
respect to the $p$-power map and the dual of the corestriction.
Then there is a description of the $\Lambda(G)$-module $Y$ as
follows:

\begin{prop}\label{DZ}
Assume that $\cd_p(\G)=2$ and that $N^{ab}(p)$ is a finitely
generated $\Lambda(\G)$-module. Then
 $$Y\simeq DZ\; \mbox{ and }\;  \E^0(Z)\cong\H^2(\h,A)^\vee,$$
thus $Y$ is determined by $Z$ up to projective summands. Suppose,
in addition, that $\H^2(\h,A)=0.$ Then
 $$\E^1(Y)\cong Z.$$
\end{prop}

For a proof of the proposition see \cite{nsw} 5.6.8 and
\cite{ochi} thm 3.13. 



%
\POP
\PUSH{local.tex}%

\section{Local Iwasawa modules}
\label{local}
\subsection{The general case} In this section we study the
structure of Iwasawa modules arising from ``$p$-adic
representations" $\G\to\mathrm{Aut}(A),$ where $\G=G_k$ is the
absolute Galois group of a finite extension $k$ of
$\mathbb{Q}_\ell$  and $A$ is a $p$-divisible $p$-torsion abelian
group of finite $\zp$-corank $r.$ Having fixed a $p$-adic Lie
extension $k_\infty$ of $k$ with Galois group $G$, we write
$\h=G(\overline{k}/k_\infty)\subseteq\G$ where $\overline{k}$
denotes the algebraic closure of $k.$ We are going to apply the
general results of  section \ref{diagram} to the module
$$X_{A}:=X_{\h,A}=\H^1(\h,A)^\vee=\H^1(k_\infty,A)^\vee,$$ i.e.\
we will determine the $\Lambda(G)$-modules occurring in the
canonical exact sequence
 $$\xymatrix@1{
    {\ 0 \ar[r] } &  {\ \E^1\mathrm{D}(X_A) \ar[r] } &  {\ X_A \ar[r] } &  {\ \E^0\E^0(X_A) \ar[r] } &  {\
   \E^2\mathrm{D}(X_A) \ar[r]}& 0. \
 }$$

The statements in this section often say that  the module $X_A$
(or another one) fits into an exact sequence of $\La(G)$-modules.
In general, this will not determine its Galois-module structure
uniquely. But if it happens that  such a sequence describes $X_A$
as 1st or 2nd syzygy  of some $\La(G)$-module with well-known
structure, then the Galois-module structure of $X_A$ is uniquely
determined  {\em up to homotopy,} i.e.\ up to projective summands
(see subsection \ref{not}).

For the sake of completeness and for the convenience of the reader
we restate some general results from \cite{ochi-ven2}, but see
also \cite{ochi}. Since we have fixed $\h,$ we shall omit it in
the notation and write $Y_A,$ $Z_A,$ etc. Recall that $G$ has
finite cohomological dimension $\cd_p G=m$ if and only if $G$ has
no element of finite $p$-power order and its dimension as $p$-adic
analytic manifold equals $m.$

\begin{lem}(cf.\ \cite{ochi})\label{localtor}
\begin{enumerate}
\item If $k$ is a finite extension of $\Q_\ell$ and $k_\infty$ is a Galois extension of $k$, then
 $Z=A^\ast(k_\infty)^\vee,$ where $A^\ast=(T_pA)^\vee(1)$ by
 definition,
\item $\E^1\mathrm{D}(X_A)\cong\E^1(A^\ast(k_\infty)^\vee),$
\item $\E^2\mathrm{D}(X_A)\subseteq\E^2\mathrm{D}(Y_A)\cong\E^2(A^\ast(k_\infty)^\vee),$
\item If $\cd_p(G)
\leq 2$ or $\cd_p(G)=3$ and $A(k_\infty)^\vee$ is
$\zp$-torsion-free, then $\mathrm{D}X_A\simeq\E^1(X_A).$
 \end{enumerate}
 \end{lem}
\begin{proof}
(i) is just local Tate duality while (ii) is a consequence of (i):
 $$\E^1\mathrm{D}(X_A)\cong\E^1\mathrm{D}(Y_A)\cong\E^1(Z_A)\cong\E^1(A^\ast(k_\infty)^\vee)$$
(Note that the first isomorphism holds because $J_A$ is
torsion-free as $\Lambda(U)$-module for a suitable open
pro-$p$-subgroup $U\subseteq G$, such that $\Lambda(U)$ is
integral). By the same reason and using the snake lemma, one sees
that $\E^2\mathrm{D}(X_A)\subseteq\E^2\mathrm{D}(Y_A).$ To prove
(iv) just note that in these cases $\pd X_A\leq 1$ by the diagram
\ref{diagram-prop}, the defining sequence (\ref{J}) of $J_A,$
corollary \cite[Cor. 6.3]{ven1} and \cite[Cor. 4.8]{ven1}.
\end{proof}

Recall that for a finitely generated abelian $p$-primary group $A$
we denote by $A_{div}$ the quotient of $A$ by its maximal
$p$-divisible subgroup. The next result generalizes  a result of
Greenberg \cite{green-adic}:

\begin{prop}(cf.\ \cite[{\S} 2]{ochi-ven2})\label{local-d1}
Let $n=[k:\mathbb{Q}_\ell],$ $\ell=p,$ be the finite degree of $k$
over $\mathbb{Q}_p$ and $k_\infty$  a Galois extension of $k$ with
Galois group $G\cong\Gamma\rtimes_\omega\Delta$, where
$\Gamma\cong\zp$ and  $\Delta$ is a finite group of order $t$
prime to $p,$ which acts on $\Gamma$ via the character
$\omega:\Delta\to\mathbb{Z}_p^\ast$. If $\chi=\omega^{-1}$ denotes
the inverse of the character which determines the action on the
$p$-dualizing module of $G$, the canonical sequence becomes
 $$\xymatrix@1{
    {\ 0 \ar[r] } &  {\ T_pA^\ast(k_\infty)(\chi) \ar[r] } &  {\ X_A \ar[r] } &  {\ P \ar[r] } &  {\
    M\ar[r] }&0, \
 }$$
where $P$ is a projective $\Lambda(G)$-module of
$\rk_{\Lambda(\Gamma)}P=rnt$ and $M$ fits into the exact sequence
 $$\xymatrix@1{
    {\ 0 \ar[r] } & M\ar[r]& {\ A^\ast(k_\infty)_{div}(\chi) \ar[r] } &  {\ \tor_\zp(A(k_\infty)^\vee)  }  \
 }.$$
Furthermore,
\begin{enumerate}
 \item if $A^\ast(k_\infty)$ is
finite, then ${T_p}A^\ast(k_\infty)(\chi)=0.$ If, in addition,
$A(k_\infty)^\vee$ is $\zp$-free, then $M\cong A^\ast(k_\infty).$
\item if $A^\ast(k_\infty)^\vee$ is $\zp$-free, then $X_A\cong P\oplus
T_pA^\ast(k_\infty)(\chi).$ In particular, $X_A$ is projective, if
$A^\ast(k_\infty)=0.$
\end{enumerate}
\end{prop}

\begin{proof}
First note that according to lemma \ref{localtor} and \cite[Cor.
4.8]{ven1}
 \begin{eqnarray*}
\E^1\mathrm{D}(X_A)& \cong & \E^1(A^\ast(k_\infty)^\vee)\\
                    & \cong &
                    \E^1(A^\ast(k_\infty)^\vee/\mathrm{tor}_\zp)\\
                    &\cong  &
                    (A^\ast(k_\infty)^\vee\otimes\qp/\zp(\chi^{-1}))^\vee\\
                    &\cong &{T_p}A^\ast(k_\infty)(\chi).
 \end{eqnarray*}
To determine $E^2\mathrm{D}(X_A)\cong\E^2\E^1(X_A)$ we use the
short exact sequences ((\ref{J}) and prop.\ \ref{diagram-prop})
 \begin{eqnarray*} &&\xymatrix@1{
   {\ 0 \ar[r] } &  {\ X_A \ar[r] } &  {\ Y_A \ar[r] } &  {\ J_A \ar[r] } &  {\
   0, } }\\ &&\xymatrix@1{
   {\ 0 \ar[r] } &  {\ J_A \ar[r] } &  {\ \Lambda(G)^d \ar[r] } &  {\ A(k_\infty)^\vee \ar[r] } &  {\
   0, }
}
 \end{eqnarray*} i.e.\ $\E^1(J_A)\cong\E^2(A(k_\infty)^\vee)\cong
 A(k_\infty)_{div}(\chi)$ by \cite[Cor. 4.8]{ven1} and
 $$\xymatrix@1{
     {\ A(k_\infty)_{div}(\chi) \ar[r] } &  {\E^1(Y_A) \ar[r] } &  {\ \E^1(X_A) \ar[r] } &  {\
    0 } \
 }$$ is exact.
 Forming  the long exact $\Ext$-sequence and applying lemmas
 \ref{localtor} and  \cite[Cor. 4.8]{ven1} again, gives the desired result.
\end{proof}

Let us now consider the case $\ell\neq p:$

\begin{prop}(cf.\ \cite[{\S} 2]{ochi-ven2})
In the situation of the last theorem but with  $\ell\neq p$ there
is an isomorphism
 $$X_A\cong T_p A^\ast(k_\infty)(\chi).$$
\end{prop}
\begin{proof}
In \cite{ochi}, prop.\ 3.12, it was calculated that the
$\Lambda(\Gamma)$-rank of $X_A$ equals the
$\Lambda(\Gamma)$-corank of $\H^2(k_\infty,A)$, but the latter
module vanishes because the order of $G$ is divisible by
$p^\infty$ (cf.\ \cite{nsw} 7.1.8).
\end{proof}

\begin{prop}(cf.\ \cite[{\S} 2]{ochi-ven2})\label{local-d2}
Let $n=[k:\mathbb{Q}_p]$  be the finite degree of $k$ over
$\mathbb{Q}_p$ and $k_\infty$  a $p$-adic Lie extension of $k$
such that its Galois group $G$ has cohomological dimension
$\cd_p(G)=2.$ Let $\Gamma\subseteq G$ be an arbitrary  open
uniform pro-$p$-subgroup, i.e.\   $\Lambda(\Gamma)$ is integral,
and let $t$ be the index $(G:\Gamma)$. If $\chi$ denotes the
inverse of the character which determines the action of $G$ on the
$p$-dualizing module, then the canonical sequence becomes
 $$\xymatrix@1{
    {\ 0 \ar[r] }  &  {\ X_A \ar[r] } &  {\ R \ar[r] } &  {\
   \E^2\mathrm{D}(X_A)\ar[r] }&0, \
 }$$
where $R$ is a reflexive $\Lambda(G)$-module with
$\rk_{\Lambda(\Gamma)}R=rnt.$ If, in addition, $A(k_\infty)^\vee$
is $\zp$-free, then $\E^2\mathrm{D}(X_A)$ fits into the exact
sequence
 $$\xymatrix@1{
    {\ 0 \ar[r] } &{\E^2\mathrm{D}(X_A)\ar[r]}&   {\ T_pA^\ast(k_\infty)(\chi) \ar[r] } &  {\ \Hom(T_p A(k_\infty),\zp)  }  \
 }.$$
\end{prop}

\begin{proof}
Using again the lemmas \ref{localtor} and  \cite[Cor. 4.8]{ven1},
the proof is completely analogous to that in the one-dimensional
case \ref{local-d1}.
\end{proof}

Note that in the case $p\neq l$ and $\cd_p(G)\geq 2$ we have
$\h=0$, i.e.\ $X_A=0,$ because the Galois group $G_k(p)$
 of the maximal $p$-extension of any local
field $k$ over $\mathbb{Q}_\ell$ is isomorphic to
$\zp(1)\rtimes\zp$ (resp.\ \zp) if $\mu_p\subseteq k$ (otherwise).
Thus it  does not have any non-trivial quotient $G$ which
satisfies these conditions.

\begin{prop}(cf.\ \cite[{\S} 2]{ochi-ven2})\label{localcd>2}
Let $n=[k:\mathbb{Q}_p]$  be the finite degree of $k$ over
$\mathbb{Q}_p$ and $k_\infty$  a $p$-adic Lie extension of $k$
such that its Galois group $G$ has cohomological dimension
$\cd_p(G)\geq 3.$ Let $\Gamma\subseteq G$ be an arbitrary  open
uniform pro-$p$-subgroup, i.e.\   $\Lambda(\Gamma)$ is integral,
and let $t$ be the index $(G:U)$. Then
 $$X_A\cong\E^0\E^0X_A$$
is a reflexive $\Lambda(G)$-module 
with $\rk_{\Lambda(\Gamma)}X_A=rnt.$
\end{prop}

\begin{proof} This follows from the lemmas \ref{localtor} and  \cite[Cor. 4.8]{ven1} as above.
\end{proof}

At the end of this part we want to restate the results concerning
the ranks of the considered modules . The  result was obtained
independently by S. Howson \cite[6.1]{howson} and Y. Ochi
\cite[thm 3.3]{ochi}, see also \cite[thm.\ 2]{ochi-ven2}.

\begin{prop}(Howson, Ochi) \label{localranks}
Let $k$ be a finite extension of $\Q_\ell$ and $k_\infty$  be a
pro-$p$ Lie extension of $k$ with Galois group $G=G(k_\infty/k)$.
As before $r$ denotes the $\zp$-rank of $\mathrm{rank}(A^\vee)$.
Assume that $\Lambda=\Lambda(G)$ is integral. Then
\[\rk_\Lambda \H^1(k_\infty,A)^\vee=\left\{\begin{array}{cl}
   r[k:\Qp] & \mbox{if }\ell=p \\
   0 & \mbox{otherwise}
 \end{array}\right..\]
\end{prop}

\begin{proof} (cf. Ochi)
Noting the vanishing of $\H^2(k_\infty,A)$ and that
$N^{ab}(p)\cong\Lambda(\G)$ for $d=[k:\Q_p]+2$ (conferring
\cite{ja-is} thm.\ 5.1 c)), the result follows from the diagram
\ref{diagram-prop} and the above remarks with respect to the case
$\ell\neq p.$
\end{proof}

 \subsection{The
case  $A=\mathbb{Q}_p/\mathbb{Z}_p$}  \subsubsection{Local units}
If we specialize to the important case $A=\qp/\zp$ with trivial
Galois action,  we are able to determine the module structure more
exactly using local class field theory: $X:=X_{\qp/\zp}\cong
\h^{ab}(p)\footnote{This notation refers to the diagram \ref{defR}
of section \ref{diagram} where we represent the absolut local
Galois group $\G$ of $k$ by a free profinite group of rank
$d=[k:\Q_\ell]+2$ according to \cite{nsw} theorem 7.4.1.} $ is the
Galois group of the maximal abelian $p$-extension of $k_\infty,$
which is canonically isomorphic to the inverse limit $$X\cong\la
(k_\infty):=\projlim{k'} \la (k')$$ of the $p$-completions $\la
(k')$ of the multiplicative groups of finite subextensions $k'$ of
$k$ in $k_\infty:$
 $$\la (k')=\projlim{m} (k')^\ast/(k')^{\ast p^m},$$
where the limit is taken via the norm maps. Since the Galois
module structure of $\la (k')$ is well known if tensored with
$\qp$, we get

\begin{thm}\label{localAd1}
Let $n=[k:\mathbb{Q}_\ell],$ $\ell=p,$ be the finite degree of $k$
over $\mathbb{Q}_p$ and $k_\infty$  a Galois extension of $k$ with
Galois group $G\cong\Gamma\rtimes_\omega\Delta$, where
$\Gamma\cong\zp$ and  $\Delta$ is a finitely generated profinite
group of order prime to $p,$ which acts on $\Gamma$ via the
character $\omega:\Delta\to\mathbb{Z}_p^\ast$. We write $k_0$ for
the fixed field of $\Gamma$ and denote by $\chi=\omega^{-1}$ the
inverse of the character which determines the action on the
$p$-dualizing module of $G$.

\begin{enumerate}
\item If $\mu_{p^\infty}\subseteq k_\infty,$ i.e.\ $k_\infty$ is the
cyclotomic $\zp$-extension of $k_0$ and $G=\Gamma\times\Delta,$
then it holds
 $$\la (k_\infty)\cong\Lambda^n\oplus\zp(1).$$
\item Let $\mu(k_\infty)(p)$ be finite. Then there is an exact
sequence of $\Lambda$-modules
 $$\xymatrix@1{
    {\ 0 \ar[r] } &  { \la (k_\infty)\oplus I_G \ar[r] } &  {\Lambda^{n+1} \ar[r] } &  { \mu(k_\infty)(p)(\chi) \ar[r] } &  {\
    0. } \
 }$$

For any presentation
 $$\xymatrix@1{
    {\ 1 \ar[r] } &  {\ K \ar[r] } &  {\ \F_{d'} \ar[r] } &  {\ G \ar[r] } &  {\
    1 } \
 }$$ by a free profinite group $\F_{d'}$ on $d'\leq n+1$ generators,
 there exists an exact sequence
 $$\xymatrix@1{
    {\ 0   \ar[r] } &  {\ \la (k_\infty)\ar[r] } &  {\Lambda^{n-d'+1}\oplus K^{ab}(p) \ar[r] } &  {\ \mu(k_\infty)(p)(\chi) \ar[r] } &  {\
    0. } \
 }$$
\end{enumerate}
\end{thm}

\begin{rem}
(i) The existence of a  presentations in (ii) is always guaranteed
by  \cite{ja-galois} theorem 4.3. Indeed, one can choose $d'=2$.\\
(ii) Using the Krull-Schmidt theorem and Maschke's theorem, it is
 easily proved (see the proof below) that
\begin{eqnarray*}
\E^0(I_G)(\omega)\oplus I_G& \cong& \mathbb{Z}_p \kl G\kr ^2 \\
\bigoplus_{i=1}^{m-1} I_G (\omega^i)\oplus I_G & \cong&
\mathbb{Z}_p\kl G\kr ^m,
\end{eqnarray*}
where $m$ denotes the order of $\omega.$  Hence, from the
isomorphism $K^{ab}(p)\oplus I_G\cong\zp\kl G\kr ^d$ according to
the Lyndon sequence (\ref{lyndon}), we get isomorphisms (for
$m\leq d$)
 \begin{eqnarray*}
K^{ab}(p)&\cong&\mathbb{Z}_p\kl G\kr
^{d-2}\oplus\E^0(I_G)(\omega)\\
          &\cong&\mathbb{Z}_p\kl G\kr ^{d-m}\oplus\bigoplus_{i=1}^{m-1} I_G
          (\omega^i).
 \end{eqnarray*}
 In particular, if $\omega$ is an involution and $d=2$, then
 $K^{ab}(p)\cong\E^0(I_G)(\omega)\cong I_G(\omega)$ holds.
\end{rem}

\begin{proof}
Let us first consider the case that $\Delta$ is a finite group,
which grants that $\La(G)$ is Noetherian. Then the statements are
consequences of theorem \ref{local-d1} once having determined the
structure of $P=\E^0\E^0X.$ We will apply the Krull-Schmidt
theorem  and we first observe that for any open normal subgroup
$U\trianglelefteq\Gamma$ and $\bar{G}:=G/U$ it holds:
$X_U\otimes\qp\cong P_U\otimes\qp$ and, if $k'$ denotes the fixed
field of $U,$ there are exact sequences of $\bar{G}$-modules
\begin{eqnarray*}
&&\xymatrix@1{
    {\ 0 \ar[r] } &  {\ U^{ab}(p) \ar[r] } &  {\ (I_G)_U \ar[r] } &  {\ \zp[\bar{G}] \ar[r] } &  {\
    \zp\ar[r] }& 0, \
 }\\
&& \xymatrix@1{
   {\ 0 \ar[r] } &  {\ X_U \ar[r] } &  {\ \bar{G}_{k'}^{ab}(p) \ar[r] } &  {\ U^{ab}(p) \ar[r] } &  {\
   0. }
}
\end{eqnarray*}
Hence, by Maschke's theorem  and using
$\bar{G}_{k'}^{ab}(p)\otimes\qp\cong\qp[\bar{G}]^n\oplus\qp$ (cf.\
\cite{nsw} 7.4.3), we get $$P_U\otimes\qp\oplus
(I_G)_U\otimes\qp\cong\qp[\bar{G}]^{n+1},$$ i.e.\ $$P\oplus
I_G\cong \Lambda^{n+1}.$$
 Now, taking $U$-coinvariants of the augmentation sequence
 $$\xymatrix@1{
    {\ 0 \ar[r] } &  {\ I_G \ar[r] } &  {\ \mathbb{Z}_p\kl G\kr  \ar[r] } &  {\ \zp \ar[r] } &  {\ 0 } \
 }$$  and tensoring with $\qp(\omega^i)$ gives
 $$\qp[\bar{G}]\oplus\qp (\omega^{i+1})\cong
 (I_G(\omega^i))_U\otimes\qp\oplus\qp(\omega^i).$$
  For (i) just note that $I_G$ is projective and $\omega$ trivial
 because $\Delta$ acts trivially on $\Gamma,$ hence:
 $I_G\cong\zp\kl G\kr.$ The first sequence in (ii) is immediate while the second one
results from the isomorphism $K^{ab}(p)\oplus I_G\cong\zp\kl G\kr
^d$ according to the Lyndon sequence (\ref{lyndon}).\\ Now let us
assume that $\Delta$ is infinite. If $\Delta'\subseteq\Delta$ is
an open subgroup then the functor obtained by  taking
$\Delta'$-coinvariants is exact because $\H_1(\Delta',M)=0$ for
any   $\La$-module $M.$ Since the automorphism group is virtually
pro-$p,$ there is an open normal subgroup $\Delta_0$ of $\Delta$
which acts trivially on $\Gamma,$ in particular any open normal
subgroup $\Delta'$ of $\Delta$ which is contained in $\Delta_0$ is
normal  in $G.$ Now a free presentation of $G$
 $$\xymatrix@1{
    {\ 1 \ar[r] } &  {\ K \ar[r] } &  {\ \F_{d'} \ar[r] } &  {\ G \ar[r] } &  {\
    1 } \
 }$$
induces a free presentation of $G':=G/\Delta'$
 $$\xymatrix@1{
    {\ 1 \ar[r] } &  {\ K_{\Delta'} \ar[r] } &  {\ \F_{d'} \ar[r] } &  {\ G/\Delta' \ar[r] } &  {\
    1. } \
 }$$
Using the Lyndon sequence, it is easy to verify that
$(I_G)_{\Delta'}\cong I_{G/\Delta'}$ and $K^{ab}(p)_{\Delta'}\cong
K_{\Delta'}^{ab}(p).$ Now the strategy is as follows. Take a
 $\La(G)$-module $M$ and show that for any $\Delta'$ as above its
$\Delta'$-coinvariants are isomorphic to certain finitely
generated $\La(G')$-modules of the same type, e.g.\ $\la
(k')\oplus I_{G'},$ where $k'$ is the fixed field of $k_\infty$ by
$\Delta'.$ Then it follows easily (using a compactness argument to
grant the existence of a compatible system of isomorphisms) that
$M\cong\la (k_\infty)\oplus I_G.$ As an example we prove the first
statement in (ii):  choose a surjection
$\La(G)^{n+1}\twoheadrightarrow \mu(k_\infty)(p)(\chi)$ and define
$M$ to be the kernel of it. Taking $\Delta'$-coinvariants and
comparing it with the result for $k',$ i.e.\ for (finite)
$\Delta/\Delta',$ we obtain an isomorphism $M_{\Delta'}\cong\la
(k')\oplus I_{G'}$ by Schanuel's lemma (see \cite[1.3]{ja-is} for
a generalized version). The other statements follow by similar
arguments.\\

The second isomorphism of the remark can be
 deduced by summing up $(I_G(\omega^i))_U\otimes\qp$ for $0\leq
 i\leq m.$ For the first one, use that due to the projectivity of
 $I_G$
  \begin{eqnarray*}
  \E^0(I_G)_U\otimes\qp&\cong&\Hom_{\mathbb{Z}_p\kl G\kr }(I_G,\mathbb{Z}_p\kl G\kr )_U\otimes\qp\\
  &\cong&\Hom_{\mathbb{Z}_p[\bar{G}]}((I_G)_U,\mathbb{Z}_p[\bar{G}])\otimes\qp\\
  &\cong&\Hom_\qp((I_G)_U,\qp)\otimes\qp
  \end{eqnarray*} holds.
\end{proof}

\begin{thm}
In the situation of the last theorem but with  $\ell\neq p$ there
is an isomorphism
 $$X\cong \left\{\begin{array}{cl}
   \mathbb{Z}_p(1)(\chi) & \mbox{ if } \mu_p\subseteq k_0 \\
  0 & \mbox{otherwise} \
 \end{array}\right..$$
\end{thm}

The next theorem generalizes results of Wintenberger \cite{winten}
who restricts himself to the case in which $G$ is abelian. It
applies for example to $\Gamma\cong\mathbb{Z}_p\rtimes\zp.$ Recall
that $R$ respectively $R^{ab}(p)$ were defined via diagram
\ref{defR}.

\begin{thm}\label{localAd2}
Let $n=[k:\mathbb{Q}_p]$  be the finite degree of $k$ over
$\mathbb{Q}_p$ and $k_\infty$  a Galois extension of $k$ with
Galois group $G\cong\Gamma\rtimes_\rho\Delta,$ where $\Gamma$
 is a pro-$p$ Lie group  of dimension $2$  and $\Delta$ is a profinite group of
order  prime to $p,$ which acts on $\Gamma$ via $\rho:\Delta\to
Aut(\Gamma).$ Let $k_0$ be the fixed field of $\Gamma$ and let
$\chi=det\rho^{-1}$ denote the inverse of the character which
determines the action on the $p$-dualizing module of $G$.
\begin{enumerate}
\item If $\mu(k_0)(p)=1$, then $X\oplus\Lambda\cong R^{ab}(p).$ If
$\rho$ is trivial, then $X\cong\Lambda^n.$
\item If $\mu_{p^\infty}\subseteq k_\infty$ and $G$ is without $p$-torsion and such that
its dualizing module is not isomorphic to $\mu_{p^\infty},$ then
there is an exact sequence of $\Lambda$-modules

 \[ \xymatrix@1{
    {\ 0 \ar[r] } &  {\ X\oplus\Lambda \ar[r] } &  {\ R^{ab}(p) \ar[r] } &  {\ \mathbb{Z}_p(1)(\chi) \ar[r] } &
    {\  0 .}    } \]

If $\rho$ is trivial, then
 $$\xymatrix@1{
    {\ 0 \ar[r] } &  {\ X \ar[r] } &  {\ \Lambda^n \ar[r] } &  {\ \mathbb{Z}_p(1) \ar[r] } &  {\
    0 } \
 } $$
is exact.
\item If $\mu(k_\infty)(p)$ and $\Delta$ are finite, then $X\cong\E^0\E^0(X)$ is
reflexive, i.e.\ there is an exact sequence
 $$\xymatrix@1{
    {\ 0 \ar[r] } &  {\ X \ar[r] } &  {\ R^{ab}(p)\ar[r] } &  {\ \Lambda \ar[r] } &  {\
    \mu(k_\infty)(p). } \
 }$$
If, in addition, $\mu(k)(p)=1,$ but $\mu(k_\infty)(p)\neq1$
 and $\chi^{-1}\neq\chi_{cycl},$ then the right map
is also surjective (in particular $X$ is not free in this case).
\end{enumerate}
\end{thm}

\begin{rem}
For extensions $k_\infty|k$ of the type $G\cong\Gamma\times\Delta$
with $\Gamma\cong\z_p^s,$ $s\geq 3$  and finite $\Delta,$ we can
consider the relative situation
 $$\xymatrix@1{
    {\ 0 \ar[r] } &  {\ X(k_\infty)_{\Gamma'} \ar[r] } &  {\ X(K_\infty) \ar[r] } &  {\ \zp \ar[r] } &  {\
    0, } \
 }$$
where $\Gamma'$ is direct factor of $\Gamma$ isomorphic to $\zp$,
i.e.\ $\Gamma\cong\Gamma'\times\z_p^{s-1},$ and $K_\infty$ is the
fixed field of $k_\infty$ with respect to $\Gamma'.$ By induction
and applying Diekert's theorem (\cite{nsw}) one reobtains at once
Wintenberger's results (but now more generally with not
necessarily abelian $\Delta$): For any irreducible character
$\chi\neq1,\chi_{cycl}$
 the component $X(k_\infty)^{e_\chi}$ is a free $\Lambda(G)^{e_\chi}$-module  of
rank $n$
 $$X(k_\infty)^{e_\chi}\cong(\Lambda(G)^{e_\chi})^n.$$
But since we already know that $\pd_\Lambda X=s-2$ for $s\geq 3,$
$X$ can not be projective in this case, i.e.\ $X(k_\infty)^{e_1}$
or    $X(k_\infty)^{e_{\chi_{cycl}}}$ is definitely not of this
type.
\end{rem}

We will prove the theorem  only for finite $\Delta$ because the
general case follows similarly as in theorem \ref{localAd1}. Just
note   that also in this case the automorphism group of $\Gamma$
is virtually pro-$p$ (see \cite[5.6]{dsms}).
\\ But before
giving the proof we need some preparation:

\begin{lem}\label{Rab}
Let $G=\Gamma\times\Delta$ be the product of a pro-$p$ Lie group
$\Gamma$ with  $\cd_p(\Gamma)=2$ and a finite group $\Delta$ of
order prime to $p.$ Then $$R^{ab}(p)\cong\Lambda^{n+1}.$$
\end{lem}

\begin{proof}
Let $U_n:=p^n\Gamma\trianglelefteq G.$ By the Lyndon sequence
(\ref{lyndon}) and using proposition \ref{euler}, we calculate the
Euler characteristic
$h_{U_n}(R^{ab}(p))=h_{U_n}(\zp)+h_{U_n}(\Lambda^{n+1})=h_{U_n}(\Lambda^{n+1}).$
The result follows.

\end{proof}

\begin{lem}\label{projective}
If in the situation of the theorem $\mu(k_\infty)(p)$ is infinite,
then  both $\E^0(X)$ and $\E^0\E^0(X)$ are projective.
\end{lem}

\begin{proof}


Since $\E^0(-)$ preserves projectives and
$\E^0\E^0\E^0(X)\cong\E^0(X)$ by \cite[Prop.\ 3.11]{ven1}, it is
sufficient to prove the statement for $\E^0\E^0(X).$ But according
to proposition \ref{reflexivepd3} the latter module is the
$2$-syzygy of $\E^3\E^1(X).$ We claim that
$$Y\simeq X\oplus\Lambda,$$ i.e.\ that
$\E^3\E^1(X)\cong\E^3\E^1(Y)\cong\E^3(\mu(k_\infty)(p)^\vee)=0,$
which implies the lemma. Indeed, due to Poincar\'{e}-duality
$$\H^2(G,\mu(k_\infty)(p))^\vee\cong\Hom_G(\mu(k_\infty)(p),D_2^{(p)})=0,$$
if $D_2^{(p)}\neq \mu_{p^\infty}.$ Hence, $Y\simeq X\oplus\Lambda$
by the second description of 4.5 b) in \cite{ja-is}\footnote{For
$\Gamma=\z^2_p$ this statement was proved by Jannsen (\cite{ja-is}
5.2 c): Though there the claimed isomorphism
$R^{ab}(p)\cong\Lambda^{d-1}$ is only correct if $\rho$ is
trivial, the  arguments (which we restated above)  still prove
$X\oplus\Lambda\simeq Y.$}.
\end{proof}

\noindent {\it Proof } \/(of the theorem).
  Let $U_m=p^m\Gamma\trianglelefteq G$ and
denote the fixed field of $U_m$ by $k_m.$ From the exact sequence
 $$\xymatrix@1{
    {\ 1 \ar[r] } &  {\ G_{k_\infty} \ar[r] } &  {\ G_{k_m} \ar[r] } &  {\ U_m \ar[r] } &  {\
    1 } \
 }$$
we obtain the associated homological Hochschild-Serre sequence
 $$\xymatrix@1@C-0.5pc{
    {\ 0=\H_2(k_m,\zp) \ar[r] } &  {\H_2(U_m,\zp) \ar[r] } &  {\ X_{U_m} \ar[r] } &  {\ G_{k_m}^{ab}(p) \ar[r] } &  {\
    \H_1(U_m,\zp)\ar[r] }&0. \
 }$$
After tensoring with $\qp,$ it follows that
 $$X_{U_m}\otimes\qp\oplus\H_1(U_m,\zp)\otimes\qp\cong\qp[\bar{G}]^n\oplus\qp\oplus\H_2(U_m,\zp)\otimes\qp,$$ where we used Maschke's
 theorem and $\bar{G}_{k_m}^{ab}(p)\otimes\qp\cong\qp[\bar{G}]^n\oplus\qp$ (cf.\
\cite{nsw} 7.4.3). On the other hand, the Euler characteristic of
the projective module $R^{ab}(p)$ can be calculated by means of
the Lyndon sequence:
 \begin{eqnarray*}
 [R^{ab}(p)_{U_m}\otimes\qp]&=&h_{U_m}(R^{ab}(p))\\
                            &=&h_{U_m}(\zp)+h_{U_m}(\Lambda^{n+1})\\
                            &=&[\qp]-[\H_1(U_m,\zp)\otimes\qp]+[\H_2(U_m,\zp)\otimes\qp]\\
                            & & +[\qp[\bar{G}]^{n+1}]
\end{eqnarray*}
and hence $X_{U_m}\otimes\qp\oplus\qp[\bar{G}]\cong
R^{ab}(p)_{U_m}\otimes\qp.$\\

Assume that $\mu(k_0)(p)=1,$ i.e.\ $\tor_\zp \la (k_0)=1$ and
$X_{U_0}$ is $\zp$-free. Therefore, since $t$ is prime to $p$, it
follows that $X_{U_0}$ is $\zp[\Delta]$-projective. If $\rho$ is
trivial, we conclude, by the calculation above under consideration
of $h_{U_m}(\zp)=0$ (by lemma \ref{euler}) and using the
Krull-Schmidt theorem, that $X_{U_0}\cong\zp[\Delta]^n$. Applying
lemma \ref{free-asympt}, gives the desired result in this case.
Anyway, these arguments show that $X$ is projective also in the
case with non-trivial $\rho,$ i.e.\ we obtain $X\oplus\Lambda\cong
R^{ab}(p)$ in the general case.\\

In order to prove (ii), we apply theorem \ref{local-d2}: Since
 $X\oplus\Lambda\simeq Y$ in this case (see the proof of lemma \ref{projective}), we
obtain
\begin{eqnarray*} \E^2\D(X)&\cong&\E^2\D(Y)\\
           &\cong&\E^2(\z_p(-1))\\
           &\cong&\z_p(1)(\chi),
\end{eqnarray*}
where we applied the lemma \ref{localtor} and \cite[2.6]{ja-is}.
Note that $\chi^{-1}(x)=\det(Ad x)=\det
\rho(x):G\to\Delta\stackrel{\det \rho}{\to}\z_p^\ast$ (cf.\
\cite{la} V 2.5.8.1). We still have to determine the module
$P=\E^0\E^0(X),$  which is projective according to lemma
\ref{projective}: it is easily seen that $P_{U_m}\otimes\qp\cong
X_{U_m}\otimes \qp,$ i.e.\ $P\oplus\Lambda\cong R^{ab}(p),$ by the
above calculations. If $\rho$ is trivial, lemma \ref{Rab} gives
the desired result.\\

The first statement of (iii) is just theorem \ref{local-d2} and
lemma  \ref{localtor}. By proposition \ref{reflexivepd3}, we
obtain an exact sequence
 $$\xymatrix@1{
    {\ 0 \ar[r] } &  {\ X \ar[r] } &  {\ P \ar[r] } &  {\ \Lambda^s \ar[r] } &  {\
    \mu(k_\infty)(p) } \
 }$$
 for some $s.$ Splitting up the sequence, taking the long exact
 $\H_i(U_m,-)$-sequences and using the above calculations, one immediately
 sees that
  $P_{U_m}\otimes\qp\cong R^{ab}(p)_{U_m}\otimes\qp\oplus\qp[\bar{G}]^{s-1},$
  i.e.\ $P\cong R^{ab}(p)\oplus\Lambda^{s-1}.$  After possibly changing the basis of $\La^d$  and using the Krull-Schmidt theorem,
  one easily sees that we can get rid off  the summand $\Lambda^{s-1}.$ \\ In order to prove
  the last statement, we assume that $\chi^{-1}\neq\chi_{cycl}$ and consider the exact sequence
 $$
 \xymatrix{
   0\ar[r] & {\E^1(X)^\vee\ar[r]}  &{\E^1(Y)^\vee\ar[r]\ar@2{-}[d]}  & {\E^1(I)^\vee\ar@2{-}[d]} \\
     &   & {\mu(k_\infty)(p)} & {\qp/\z_p(\chi^{-1}).} \
 }$$
The decomposition of the sequence with respect to the irreducible
$\qp$-\linebreak characters of $\Delta$ gives
$(\E^1(X)^\vee)^{\chi_{cycl}}=\mu(k_\infty)(p)^{\chi_{cycl}}=
\mu(k_\infty)(p).$ \nopagebreak \hspace*{\fill}
$\square$\smallskip\\

\subsubsection{Principal units}

When $l=p,$ we are also interested in the $\Lambda$-structure of
the inverse limit of the principal units
 $$\lu^1(k_\infty):=\projlim{k'} \lu^1(k'),$$
where $k'$ runs through all finite subextensions of $k_\infty|k$
and the limit is taken with respect to the norm maps.

\begin{prop}\label{localunits}
Let $k$ be a finite extension of $\qp$ and $k_\infty$ a Galois
extension of $k$.
\begin{enumerate}
\item If $k_\infty$ contains the maximal unramified $p$-extension
of $k,$ i.e.\ if $p^\infty$ divides the degree of the residue
field extension associated with $k_\infty|k,$ then
 $$\lu^1(k_\infty)\cong\la (k_\infty).$$
\item In the other case there is the following exact sequence
 $$\xymatrix@1{
    {\ 0 \ar[r] } &  {\ {\lu}^1(k_\infty) \ar[r] } &  {\la (k_\infty)\ar[r] } &  {\zp \ar[r] } &  {\
    0. } \
 }$$
\end{enumerate}
\end{prop}

\begin{proof}
For finite extensions $K'|K|k$ of $k$ with associated residue
field extensions $\lambda'|\lambda|\kappa$ consider the following
commutative diagram with exact rows

  $$  \xymatrix{
    0\ar[r] & {\lu^1(K')/p^m}\ar[r]\ar[d]^{N_{K'/K}} & {\la(K')\ar[r]\ar[d]^{N_{K'/K}}} & {\z_p/p^m \ar[r]\ar[d]^{[\lambda':\lambda]}}& 0 \\
   0\ar[r] & {\lu^1(K)/p^m}\ar[r] & {\la(K)\ar[r]} & {\z_p/p^m \ar[r]}&   0.
   \
  }$$
While in case (i) the inverse limit $\projlim{K,m} \zp/p^m$
vanishes, because for any $m$ and any $K$ there is an extension
$K'$ such that $p^m|[\lambda':\lambda],$ in the second case it is
isomorphic to $\zp.$
\end{proof}

\begin{thm}
Assume in the situation of theorem \ref{localAd2} that $k_\infty$
contains $\mu_{p^\infty}$ but not the  maximal unramified
$p$-extension of $k.$ Then there exists an exact sequence
$$\xymatrix@1{
    {\ 0 \ar[r] } &  {\lu^1(k_\infty)\oplus\Lambda \ar[r] } &  {\ R^{ab}(p) \ar[r] } &  {\ M \ar[r] } &  {\
    0, } \
 }$$
 where $M$ fits into the exact sequence
 $$\xymatrix@1{
    {\ 0 \ar[r] } &  {\ \zp \ar[r] } &  {\ M \ar[r] } &  {\ \zp(1)(\chi) \ar[r] } &  {\
    0. } \
 }$$
In particular, if $\rho$ is trivial, there exists an exact
sequence
 $$\xymatrix@1{
    {\ 0 \ar[r] } &  {\lu^1(k_\infty) \ar[r] } &  {\ \Lambda^n \ar[r] } &  {\ M \ar[r] } &  {\
    0 .} \
 }$$

\end{thm}

\begin{proof}
Evaluating the long exact $\E^i$-sequence associated with the
exact sequence from the proposition above and noting that
$\pd_\Lambda \lu^1(k_\infty)\leq 1$ due to $\pd_\Lambda \la
(k_\infty)\leq 1$ and $\pd_\Lambda \zp= 2,$ one obtains
 that
\begin{enumerate}
\item $\E^0(\lu^1(k_\infty))\cong\E^0(X),$
\item $\E^1\mathrm{D}(\lu^1(k_\infty))=0$ and  an exact sequence
\item $\xymatrix@1{
    {\ 0 \ar[r] } &  {\ \zp \ar[r] } &  {\E^2\mathrm{D}(\lu^1(k_\infty)) \ar[r] } &  {\ \zp(1)(\chi) \ar[r] } &  {\
    0. } \
 }$
\end{enumerate}
Here we used that $\E^2\E^2(\zp)\cong\zp,$ because $\zp$ is  a
Cohen-Macaulay module of dimension 2. The result follows from the
canonical sequence.
\end{proof}

\begin{rem}
In the situation of theorem \ref{localAd1} with trivial action of
$\Delta$ the structure of the principal units is described in
\cite{nsw} as follows:
\begin{enumerate}
\item If $\mu_{p^\infty}\subseteq k_\infty,$ then
$$\lu^1(k_\infty)\cong\Lambda^n\oplus\zp(1).$$
\item If $\mu(k_\infty)(p)$ is finite, then there is an exact
sequence $$\xymatrix@1{
   {\ 0 \ar[r] } &  {\ \lu^1(k_\infty) \ar[r] } &  {\ \Lambda^n \ar[r] } &  {\ \mu(k_\infty)(p).  }
}$$
\item If $k_\infty|k$ is unramified, then $$\lu^1(k_\infty)\cong\la
(k_\infty).$$
\end{enumerate}
\noindent But the  proof of \cite{nsw} works also if $\omega$ is
not trivial.
\end{rem}

\subsection{The local CM-case}

As a consequence of theorem \ref{localAd2} we can also determine -
up to homotopy -  the structure of $X_A=\H^1(k_\infty,A)^\vee$ in
the trivializing case, i.e.\ $k(A)\subseteq k_\infty:$

\begin{thm}\label{triviald2}
Let $n=[k:\mathbb{Q}_p]$  be the finite degree of $k$ over
$\mathbb{Q}_p$ and $k_\infty$  a Galois extension of $k$ with
Galois group $G\cong\Gamma\rtimes_\rho\Delta,$ where $\Gamma$
 is a pro-$p$ Lie group  of dimension $2$  and $\Delta$ is a finite group of
order $t$ prime to $p,$ which acts on $\Gamma$ via $\rho:\Delta\to
Aut(\Gamma).$ Let $k_0$ be the fixed field of $\Gamma$ and let
$\chi=det\rho^{-1}$ denote the inverse of the character which
determines the action on the $p$-dualizing module of $G$. For any
$A$ with $\rk_\zp A^\vee=r$ such that $k(A)\subseteq k_\infty$ the
following is true.
\begin{enumerate}
\item If $\mu(k_0)(p)=1$, then $X_A\oplus\Lambda^r\cong R^{ab}(p)[A],$ in particular, if
$\rho$ is trivial: $X_A\cong\Lambda^{nr}.$
\item If $\mu_{p^\infty}\subseteq k_\infty$ and $G$ is  $p$-torsion-free and
its dualizing module is not isomorphic to $\mu_{p^\infty},$ then
there is an exact sequence of $\Lambda$-modules

 $$ \xymatrix@1{
    {\ 0 \ar[r] } &  {\ X_A\oplus\Lambda^r \ar[r] } &  {\ R^{ab}(p)[A] \ar[r] } &  {\ A^\vee (1)(\chi) \ar[r] } &  {\
    0 .} \ } $$

In particular, if $\rho$ is trivial, then
 $$\xymatrix@1{
    {\ 0 \ar[r] } &  {\ X_A \ar[r] } &  {\ \Lambda^{nr} \ar[r] } &  {\ A^\vee(1) \ar[r] } &  {\
    0 } \
 } $$
is exact.
\item If $\mu(k_\infty)(p)$ is finite, then $X_A\cong\E^0\E^0(X_A)$ is
reflexive, i.e.\ there is an exact sequence
 $$\xymatrix@1{
    {\ 0 \ar[r] } &  {\ X_A \ar[r] } &  {\ R^{ab}(p)[A] \ar[r] } &  {\ \Lambda^{r} \ar[r] } &  {\
    \mu(k_\infty)(p)[A]. } \
 }$$
If, in addition, $\mu(k)(p)=1,$ but $\mu(k_\infty)(p)\neq1$
 and $\chi^{-1}\neq\chi_{cycl},$ then the right map
is also surjective (in particular, $X_A$ is not free in this
case).
\end{enumerate}

\end{thm}
\begin{proof}
In this case the subgroups $\h, R$ and $N$ act trivially on
$A=A(k_\infty),$ i.e.\ $X_A\cong X[A].$
\end{proof}

This result applies  to the following situation: Let $K$ be a
imaginary quadratic   number field, $F$ a finite, abelian
extension of $K$ and $E$ an elliptic curve defined over $F$ with
complex multiplication (CM) by the ring of integers $\O_K$ of $K$
such that $F(E_{tor})$ is an abelian extension of $K.$ Assume that
the rational prime $p$ splits in $K,$ i.e.\ $p\O_K=\p\bar{\p},$
$\p\neq\bar{\p},$ and that $E$ has good reduction at all places
lying over $p.$ Set $G=G(F(E(p))/F)_{\P}$ the decomposition group at
some $\P|\p.$ According to \cite[1.9]{deSh}, the prime $\P$
ramifies totally in $F(E(\p))|F$ and decomposes only finitely (and
is unramified) in $F(E(\bar{\p}))|F.$ Therefore the decomposition
group $G$ is an open subgroup of $G(F(E(p))/F),$ i.e.\ of type
$\z_p^2\times\Delta$ where $\Delta$ is a finite abelian group.
Thus we obtain an exact sequence
 $$\xymatrix@1{
    {\ 0 \ar[r] } &  {\ \H^1(F(E(p))_\P,E(p))^\vee \ar[r] } &  {\ \La(G)^{2n} \ar[r] } &  {\ T_pE \ar[r] } &  {\
    0, }
 }$$
where $n=[F_\p:\qp].$ By the same argument, but now using theorem
\ref{localAd1}(ii), there  exists an exact sequence
 $$\xymatrix@1{
    {\ 0 \ar[r] } &  {\ \H^1(F(E(\p)_\P,E(\p))^\vee \ar[r] } &  {\ \La(G')^{n} \ar[r] } &  {\ \mu(F(E(\p)_\P)[E(\p)] \ar[r] } &  {\
    0, } \
 }$$
 where $G'=G(F(E(\p))/F)_{\P},$ and a similar one for $\bar{\p}.$
\POP
\PUSH{global.tex}%
\section{Global Iwasawa modules}\label{global}

 Let  $k_\infty$ be a $p$-adic Lie
extension of the number field $k$ contained in $k_S$ with Galois
group $G$ and let $A$ be a $p$-divisible $p$-torsion abelian group
with $\zp$-corank $r$ and on which $G_S(k)=G(k_S/k)$ acts
continuously where $S$ is a finite set of places of $k$ containing
all places $S_p$ over $p$ and all infinite places $S_\infty$ (and
by definition all places at which $A$ is ramified). Here  $k_S$
denotes the maximal $S$-ramified extension of $k,$ i.e.\ the
maximal extension of $k$ which is unramified outside $S.$ In order
to derive information about the $\Lambda=\Lambda(G)$-modules
$\H^i(G(k_S/k_\infty),A)$ we would like to apply the diagram
(\ref{defR}) to the group $\G=\G_S:=G(k_S/k).$
  On the other hand we have to guarantee that $\G$ is
finitely generated as a profinite group which, unfortunately, is
not known for the group $\G_S.$ But using a theorem of Neumann,
i.e.\ the inflation maps are  isomorphisms
 $$\H^i(G(\Omega/k_\infty),A)\cong\H^i(G_S(k_\infty),A),\hspace{1cm}i\geq 0,$$ for any $(p,S)$-closed extension $\Omega$ of $k$
  (i.e.\ $\Omega$ is a
 $S$-ramified extension of $k$ which does not possess any
 non-trivial $S$-ramified $p$-extension) and for any
 $p$-torsion $G(\Omega/k_\infty)$-module $A,$
we are free to replace $G_S(k)$ for example by the Galois group
$\G:=G(\Omega/k)$ where $\Omega$ is the maximal $S$-ramified
$p$-extension of $k'(A)$ and $k'$ is a Galois subextension of
$k_\infty/k$ such that $G(k_\infty/k')$ is an open (normal)
pro-$p$-group.  Regarding this technical detail, we assume in what
follows that $k_\infty$ is contained in such a $(p,S)$-closed
field $\Omega.$ Then, since $\G$ has an open pro-$p$ Sylow group,
it is finitely generated and has $\cd_p(\G)\leq 2$ for odd $p.$
Note that $Y_{S,A}:=Y_{G(\Omega/k_\infty),A}$ (\ref{Y}) and
$X_{S,A}:=X_{G(\Omega/k_\infty),A}$ (\ref{X}) do not depend on the
choice of $\Omega.$ The next lemma shows among other things that
the corresponding module $Z$ (\ref{Z}) only depends on $k_\infty,$
$A$ and $S.$ Recall that $T_pA=\Hom(\qp/\zp,A)$ denotes the ``Tate
module" of $A.$ We shall write
$\H^\ast_{cts}(G_S(k),T_pA)\cong\projlimssc{n}
\H^\ast(G_S(k),{_{p^n}A})$ for the continuous cochain cohomology
groups (see \cite[II.\S3.]{nsw}).

\begin{lem}\label{globalZ}
 Let $k,$  $k_\infty$ and $A$ be as above. Then
 $$Z_{S,A}:=Z_{G(\Omega/k_\infty),A}\cong\projlim{k\subseteq k^\prime \subseteq k_\infty} H^2_{cts}(G_S(k^\prime), T_pA).$$
\end{lem}


A basic structure result is the following

\begin{thm}\label{nospeudonullglobalXA} Let\ $G$ a $p$-adic Lie group  without
$p$-torsion. If the ``weak Leopoldt conjecture holds for $A$ and
$k_\infty$", i.e.\ $\H^2(G_S(k_\infty),A)=0$, then neither
$Y_{S,A}$ nor $X_{S,A}\cong\H^1(G_S(k_\infty), A)^\vee$ have
non-zero pseudo-null submodules.
\end{thm}

\begin{proof}
The conditions of theorem  \ref{nospeudonullXA} are fulfilled.
\end{proof}

Furthermore, the \La-rank  of $X_{S,A}$  can be determined, using
the diagram \ref{diagram}:

\begin{thm}(Ochi \cite{ochi-ven2})\label{gl-rank} Let $k_\infty|k$ be a $p$-adic
pro-$p$ extension. Assume  that $k(A)|k$ is  a pro-$p$-extension
and that $\La$ is an integral domain. Then

$$\rk_\Lambda \H^1(G_S(k_\infty), A)^\vee-\rk_\Lambda \H^2(G_S(k_\infty), A)^\vee=r_2(k)r$$
Here $r_2(k)$ denotes as usual the number of complex places of
$k.$
\end{thm}

Thus, if the weak Leopoldt conjecture holds for $A$ and
$k_\infty,$ one obtains a simple formula for the \La-rank of
$\H^1(G_S(k_\infty), A)^\vee.$ So, we conclude with a brief
discussion and motivation  concerning this conjecture:

In \cite{ja-ladic}  Jannsen extended the {\em strong} Leopoldt
conjecture, which is equivalent to the vanishing of
$\H^2(G_S(k),\qp/\zp),$ to the following setting: Let $X$ be a
smooth projective variety of pure dimension over $k$ and assume
that $S$ contains $S_p,$ $S_\infty$ and all places $S_{bad}$ where
$X$ has bad reduction. Then the \'{e}tale cohomology
$T^i(n):=\H^i_{\acute{e}t}(X\times_k \bar{k},\zp(n))$ is a compact
$G_S(k)$-module which is finitely generated over $\zp;$ here
$\bar{k}$ denotes as usual an algebraic closure of $k.$ Hence
$A^i(n):=T^i(n)\otimes_\zp \qp/\zp$ is a $p$-divisible discrete
$G_S(k)$-module, for $X=Spec(k)$ and $i=0$ isomorphic to
$(\qp/\zp)(n).$ Assuming $p\neq2$ or that $k$ is totally imaginary
his conjecture (cf.\ \cite[Conjecture 1, Lem.\ 1]{ja-ladic})
predicts

\[\H^2(G_S(k),A^i(n))=0\;\; \mbox{ if } \left\{\begin{array}{cl}
  \mbox{(i)} & i+1<n, \mbox{ or } \\
  \mbox{(ii)} & i+1>2n.
\end{array}\right.\] Thus, if this conjecture true for fixed $X,$
$i$ as well as $n$ for {\em all} number fields contained in
$k_\infty,$ it implies obviously the weak Leopoldt conjecture for
$A^i(n)$ over $k_\infty.$  While in the ``unstable" range $n\leq
i+1 \leq 2n$ the cohomology group $\H^2(G_S(k),A^i(n))$ is
non-trivial in general, it is supposed to vanish after going up a
``nice" $p$-adic Lie-extension (cf.\ corollary \ref{qp1} for an
example of this phenomena).

It is a result of Iwasawa that over the cyclotomic $\zp$-extension
of {\em any} number field the {\em original} weak Leopoldt (i.e.\
for $A=\qp/\zp$) holds and consequently for  $(\qp/\zp)(n)$ for
all $n\in\mathbb{Z}$ (see \cite[10.3.25]{nsw} for a cohomological
proof). This leads to

\begin{rem}\label{weakLeopold}
The weak Leopoldt conjecture for $A$ and $k_\infty$ holds for
example if $k(A)$ and the cyclotomic $\zp$-extension of $k$  are
contained in $k_\infty.$   The claim follows by expressing
$\H^2(G_S(k_\infty),A)$ (considered as abelian group) as direct
limit $\dirlimssc{k'} \H^2(G_S(k'_{cyc}),\qp/\zp)^r,$ where $k'$
runs through the finite extensions of $k$ in $k_\infty.$
\end{rem}

For a discussion about the weak Leopoldt conjecture over the
cyclotomic $\zp$-extension of a number field for other $p$-adic
representations than above we refer the reader to section 1.3 and
Appendix B of \cite{perrin-Lfunc}.
\POP
\PUSH{gm.tex}%
\subsection{The multiplicative group $\mathbb{G}_m$}
\subsubsection{The maximal
abelian $p$-extension of $k_\infty$ unramified outside $S$}
\label{maxabp-ext}
 We still consider $p$-adic Lie extensions
$k_\infty|k$ with Galois group $G=G(k_\infty/k)$ such that
$k_\infty$ is contained in the maximal $S$-ramified  extension
$k_S$ of $k.$ Here, as before, $S$ denotes a finite set of places
of $k$ containing all places $S_p$ over $p$ and all infinite
places $S_\infty.$ For $K|k$ finite let $S_f(K)$ be the set of
finite primes in $K$ lying above $S.$ In this paragraph we
specialize to the case $A=\qp/\zp$ and we shall write  $X_S$ for
the $\La=\La(G)$-module $X_{S,\qp/\zp}$ (\ref{X})
$$X_S:=X_{S,\qp/\zp} =\H^1(G_S(k_\infty),\qp/\zp)^{\vee}\cong
G(k_S/k_\infty)^{ab}(p),$$ and respectively for $Y_S$ (\ref{Y})
and $Z_S$ (\ref{Z}).

In this case  theorem \ref{nospeudonullglobalXA} is  a
generalization of the theorems of Greenberg \cite{green} and
Nguyen-Quang-Do \cite{Ng}, who considered the case
$G\cong\mathbb{Z}_p^d.$ Indeed, it confirms Greenberg's suggestion
 that an analogous statement
also should hold for $p$-adic Lie extensions.

\begin{thm}\label{nospeudonullglobalX} Let $G$ be a $p$-adic Lie group  without
$p$-torsion. If the ``weak Leopoldt conjecture holds for
$k_\infty$", i.e.\ $\H^2(G_S(k_\infty),\qp/\zp)=0$, then
$X_{S}\cong G_S(k_\infty)^{ab}(p)$ has no non-zero pseudo-null
\La-submodule.
\end{thm}

\begin{rem}
Recall that the weak Leopoldt conjecture for  $k_\infty$ holds  if
 the cyclotomic $\zp$-extension of $k$ is contained in
$k_\infty.$
\end{rem}

We will also consider the \La-modules
\begin{eqnarray*}
X_{nr}&=&G(L/k_\infty), \\X_{cs}^S&=&G(L'/k_\infty),
\end{eqnarray*} where $L$
is the maximal abelian unramified pro-$p$-extension of $k_\infty$
and $L'$ is the maximal subextension in which every prime above
$S$ is completely decomposed.

 For an arbitrary number field $K$, we denote the ring
of integers  (resp.\ $S$-integers) by ${\mathcal{O}}_{K}$ (resp.\
${\mathcal{O}}_{K,S}$) and its units by
${E(K):=\mathcal{O}}_{K}^\times$ (resp.\
${E_S(K):=\mathcal{O}}_{K,S}^\times$). Then we define
\begin{eqnarray*}
\gu:&=&\projlim{k'} ({\mathcal{O}}_{k'}^\times\otimes\zp),\\
\gus:&=&\projlim{k'} ({\mathcal{O}}_{k',S}^\times\otimes\zp),
\end{eqnarray*} where the
limit  is taken with respect to the norm maps. This should not be
confused with the discrete module of units (resp.\ $S$-units)
$E(k_\infty)=\dirlim{k'} E(k')$ (resp.\ $E_S(k_\infty)=\dirlim{k'}
E_S(k')).$

Finally, we write for the local-global modules

\begin{eqnarray*}
 \gla &=& \bigoplus_{S_f(k)} \Ind^{G_\nu}_G \la_\nu,\\
 \glu &=& \bigoplus_{S_f(k)} \Ind^{G_\nu}_G \lu_\nu,
\end{eqnarray*}
where $\la_\nu=\la(k_{\infty,\nu})$ (resp.\
$\lu_\nu=\lu^1(k_{\infty,\nu})$) are the local modules introduced
in chapter 2.\\ The above modules are connected via global class
field theory and the Poitou-Tate sequence as follows

\begin{prop}(Jannsen)\label{tate-poitou} There are the following exact and commutative
diagrams of \La-modules:
\begin{enumerate}
\item
$$\hspace{0cm}
                \xymatrix{
  0\ar[r] & {\H^2(G_S(k_\infty),\qp/\zp)^\vee\ar[r]\ar@{=}[d]} &{\gu\ar@{^{(}->}[d]\ar[r] }& { \glu \ar@{^{(}->}[d]\ar[r]}& {X_S\ar@{=}[d]\ar[r]} & {X_{nr}}\ar@{->>}[d]\ar[r] & 0 \\
  0\ar[r]  &{\H^2(G_S(k_\infty),\qp/\zp)^\vee\ar[r]} & {\gus} \ar[r] & {\gla \ar[r]}& {X_S\ar[r]} & {X_{cs}^S\ar[r]} & 0
  }$$

\item
$$\hspace{0cm}
                \xymatrix@1{
   {\ 0 \ar[r] } &  {\ \gu \ar[r] } &  { \gus \ar[r] } &  { \bigoplus_{S_{un}(k)} \Ind^{G_\nu}_G\zp \ar[r] } &  {\  X_{nr} \ar[r] }&  {\  X_{cs}^S\ar[r] }&  {\   0, }
}$$ where $S_{un}:=\{\nu\in S(k)\:|\: p^\infty\nmid f_\nu \}$ and
$f_\nu=f(k_{\infty,\nu}/k_\nu)$ denotes the degree of the
extension of the corresponding residue class fields.
\item
$$ \hspace{0cm}
                \xymatrix@1{
   {\ 0 \ar[r] } &  {X_{cs}^S \ar[r] } &  {Z_{S,\qp/\zp(1)} \ar[r] } &  {\bigoplus_{S_f(k)} \Ind^{G_\nu}_G \zp\ar[r] } &  {\
   \zp\ar[r] } & 0,
}$$

and, if $\mu_{p^\infty}\subseteq k_\infty,$
$$ \hspace{0cm}
                \xymatrix@1{
   {\ 0 \ar[r] } &  {X_{cs}(-1) \ar[r] } &  {Z_S \ar[r] } &  {\bigoplus_{S_f(k)} \Ind^{G_\nu}_G \zp(-1)\ar[r] } &  {\
   \zp(-1)\ar[r] } & 0.
}$$
 In particular,  $X_{cs}^S=X_{cs}:=X_{cs}^{S_p}$ is independent of
 $S$ in this case.
\item $N^{ab}(p)$ (see (\ref{defR})-(\ref{lyndon})) is a finitely generated, projective
$\La(G({k_\infty}_S(p)/k))$-module and, if the free presentation
of $\G=G({k_\infty}_S(p)/k)$ (cf.\ section \ref{diagram})  is
chosen such that $d\geq r_1'+r_2+1,$ then
 $$N^{ab}(p)_{G_S(k_\infty)(p)}\cong\bigoplus_{S'_\infty}\Ind^{G_\nu}_G\zp\oplus\La(G)^{d-r_2-r_1'-1},$$
where $S'_\infty$ is the set of real places of $k$ which ramify
(i.e.\ become complex) in $k_\infty,$ $r_1'$ is the cardinality of
$S'_\infty,$ and $r_2$ is the number of complex places of $k.$
\end{enumerate}
\end{prop}

\begin{proof}
The assertions (i) and (iii) are obtained by taking inverse limits
of the Tate-Poitou sequence (see \cite[Thm 5.4]{ja-is}) and
recalling lemma \ref{globalZ} while (ii) follows from (i) by the
snake lemma and prop.\ \ref{localunits}.
\end{proof}

>From these diagrams and the fact that \La\ is Noetherian it
follows that the modules $X_{nr}, X_{cs}^S$ are finitely
generated. Furthermore, S. Howson \cite[7.14-7.16]{howson} and
independently Y. Ochi \cite[4.10]{ochi} proved that $X_{nr}$ and
$X_{cs}^S$ are \La-torsion. Actually, this result was first proved
by M. Harris \cite[thm 3.3]{harris} but, as S. Howson remarked,
his proof is incomplete because it relies on the false ``strong
Nakayama" lemma (loc.\ cit.\ lem 1.9), see the discussion in
\cite{bal-how}. However, in a recent correction M. Harris
\cite{harris-cor} has given a new proof of the result. In the case
$G\cong\mathbb{Z}_p^d,$ this result is originally proved by
Greenberg \cite{green-inv}.

\begin{cor}\label{qp1}
\begin{enumerate}
\item If $\H^2(G_S(k_\infty),\Qp/\zp(1))=0$ (e.g.\ if  $\dim(G_\nu)\geq 1$ for all $\nu\in S_f$), then $X_{cs}^S$ is  a
\La-torsion module.
\item If $\dim(G_\nu)\geq 1$ for all $\nu\in S_f,$ then $X_{nr}$ is a
\La-torsion module.
\end{enumerate}
\end{cor}

For example, the conditions of the corollary are satisfied if
$k_\infty$ contains the cyclotomic $\zp$-extension.

\begin{proof} (cf.\ \cite{ochi-ven2}) 
The first statement follows from \ref{DZ} while the second one is
a consequence of the first one and the above proposition (To
calculate the (co)dimension of $\Ind^{G_\nu}_G \zp$ use
\cite[4.8,4.9]{ven1}. Note that the condition ``$\dim(G_\nu)\geq
1$ for all $\nu\in S_f$" implies, using Tate-Poitou duality,
\begin{eqnarray*}
\H^2(G_S(k_\infty),\Qp/\zp(1))&=&\sha^2(G_S(k_\infty),\mu_{p^\infty})\\
                             &=& \dirlim{k',n}
                              \sha^1(G_S(k'),\mathbb{Z}/p^n)^\vee\\
                              &=&\dirlim{k'}
                              Cl_S(k')\otimes_\mathbb{Z}\qp/\zp\\
                              &=&0
\end{eqnarray*}
because $Cl_S(k')$ is finite.
\end{proof}

\begin{thm}\label{XcspseudotorX}
If  $\mu_{p^\infty}\subseteq k_\infty,$ and $\dim(G_\nu)\geq 2$
for all $\nu\in S_f,$ then
 $$X_{nr}(-1)\sim
 X_{cs}^S(-1)\sim\E^1(Y_S)\sim\E^1(\tor_\La Y_S)\cong\E^1(\tor_\La X_S).$$
If, in addition, $G\cong \mathbb{Z}^r_p,$ $r\geq2,$ then even the
following holds:
 $$X_{nr}(-1)\sim X_{cs}^S(-1)\sim (\tor_\La X_S)^\circ,$$
 where $^\circ$ means that  $G$ operates via the involution
 $g\mapsto g^{-1}.$
 \end{thm}

\begin{rem}
In case $\tor_\La X_S$ is isomorphic in $\mod/\mathcal{ PN}$ to a
direct sum of cyclic modules of the form $\La$ modulo a (left)
principal ideal the proposition \ref{circ} implies that
\[E^1(\tor_\La X_S)\equiv (\tor_\La X_S)^\circ \ \mbox{\rm  mod
 }{\mathcal{ PN}}\] holds under the conditions of the theorem.
\end{rem}

\begin{proof} Note that $\H^2(G_S(k_\infty),\qp/\zp)=0,$ since
remark \ref{weakLeopold} applies. The first two
pseudo-isomorphisms follow again from proposition
\ref{tate-poitou} using \cite[4.8,4.9]{ven1} and \ref{DZ}. The
third one is just \cite[Prop.\ 3.13]{ven1}. Note that there is
even an isomorphism $\tor_\La Y_S\cong\tor_\La X_S$ because the
augmentation ideal $I_G$ is torsion-free.
\end{proof}

The following consequence generalizes a result of McCallum
\cite[thm 8]{mcCal} who considered the $\mathbb{Z}_p^r$-case:

\begin{cor}
With the assumptions of the theorem the following holds.
\begin{enumerate}
\item There is a pseudo-isomorphism $$\tor_\La X_S\sim\E^1(X_{cs}^S(-1)).$$

\item If $\dim(G)\geq 3,$ then there is an isomorphism $$\tor_\La X_S\cong\E^1(X_{cs}^S(-1)).$$
\end{enumerate}
\end{cor}

\begin{proof}
The cokernel $K:=coker(X_{cs}^S(-1)\hookrightarrow
Z_S\cong\E^1(Y_S)$ is pseudo-null, i.e.\ $\E^1(K)=0.$ If
$\dim(G)\geq 3,$ then $\E^2(K)=0,$ too, as can be calculated using
\cite[Prop.\ 2.7]{ven1}. Now, the long exact $\E$-sequence gives
the result observing $\E^1\E^1(Y_S)\cong\E^1\du Y_S\cong\tor_\La
Y_S\cong\tor_\La X_S.$
\end{proof}

\begin{rem}
The condition ``$\dim(G_\nu)\geq 2$ for all $\nu\in S_f,$" is
known to hold in ``most" extensions arising from geometry, e.g.
for the set $S_f=S_{bad}\cup S_p,$  if $k_\infty=k(\A(p))$ arises
by adjoining the $p$-division points of
 an abelian variety $\A$ over
$k$ with good reduction at all places dividing $p$ and such that
$G(k_\infty/k)$ is a pro-$p$-group without $p$-torsion,  see (the
proof of) corollary \ref{supersingular-ab-var} below. The latter
condition is satisfied if for instance $k$ contains $k({_p}\A)$
for $p\neq2$ or
$k({_{p^2}}\A)$  for $p=2,$  see at the beginning of section \ref{not}.\\
 Other important cases are the following ones:\\
(a) Let $k_\infty$ be the maximal multiple $\zp$-extension $
\tilde{k}$ of $k,$ i.e.\ the composite of all
$\zp$-extensions of $k,$ and assume that $\mu_{2p}\subseteq k$ or\\
(b) let $k_\infty$ be a multiple $\zp$-extension with
$G\cong\mathbb{Z}^r_p,$ $r\geq2,$ and assume that there is only
one prime of $k$ lying over $p.$\\ Then, as has been observed
independently by T. Nguyen-Quang-Do \cite[thm 3.2]{Ng-preprint}
and McCallum \cite[proof of thm 7]{mcCal}, the condition holds for
$S=S_p\cup S_\infty.$ Indeed, since $\Q(\mu_{2p})$ has only one
prime dividing $p,$ it suffices to consider the second case. But
then all inertia groups $T_\nu,$ $\nu\in S_p,$ are conjugate, thus
they are all equal and hence an open subgroup of $G$ due to the
finiteness of the ideal class group.
\end{rem}

With respect to the composite $\tilde{k}$ of all $\zp$-extensions
of $k$ there is the following outstanding

\begin{bconj}[R. Greenberg]
For any number field $k,$ the $\La(G(\tilde{k}/k))$-module
$X_{nr}$ is pseudo-null.
\end{bconj}

Recently, W. McCallum \cite{mcCal} proved this conjecture for the
base field $k=\Q(\mu_p)$ under some mild assumptions. For a list
of other cases in which this conjecture is known to hold, see
\cite[rem 4.6]{Ng-preprint}. Assume that $\mu_p\subseteq k$ and
that the condition ``$\dim(G_\nu)\geq 2,$ for all $\nu\in S_f,$"
holds. Then, by the above theorem and theorem
\ref{nospeudonullglobalX}, Greenberg's conjecture is equivalent to
the statement that $X_S$ is \La-torsion-free, compare with
\cite[4.4]{Ng-preprint} and \cite[Cor 13]{mcCal}.\\

The observation of the previous proof leads also to:

\begin{prop}
If $\dim(G_\nu)\geq 2$ for all $\nu\in S_f,$ then
 $$\gu\sim\gus.$$
\end{prop}

We are also interested in the  (Pontryagin duals of the) direct
limits
 \begin{eqnarray*}
 Cl_S(k_\infty)(p)&=&\dirlim{k'} Cl_S(k')(p),\\
 {\mathcal{ E}}_S(k_\infty)&:=&(E_S(k_ \infty)\otimes_\z\qp/\zp)^\vee,
 \end{eqnarray*}
of the $p$-part of the ideal class group, resp.\ of the global
($S$-)units of finite extensions $k'$ of $k$ inside $k_\infty.$
\begin{prop}\label{ClSClT}
Let $T$ be a set of places of $k$ such that $S_\infty\subseteq
T\subseteq S.$ Assume that $\dim (T_\nu)\geq 1$ for all $\nu\in
S\setminus T,$ where $T_\nu\subseteq G_\nu$ denotes the inertia
group of $\nu.$
\begin{enumerate}
\item There is an exact sequence of \La-modules
 $$\xymatrix@1@C=12pt{
    {\ 0 \ar[r] } &  {\ Cl_S(k_\infty)(p)^\vee \ar[r]^\psi } &  {\ Cl_T(k_\infty)(p)^\vee \ar[r] } &  {\ \e_S(k_\infty) \ar[r]^{\varphi} } &  {\ \e_T(k_\infty) \ar[r] } & {\
    0. } \
 }$$
\item Assume that $S\setminus T=\{\nu\}.$ Then, if $\dim(G_\nu)\geq1$ (resp.\ $\dim(G_\nu)\geq2$), then $coker(\psi)\cong
ker(\varphi)$ is \La-torsion (resp.\ pseudo-null).
\item If $\dim(G_\nu)\geq2$ for every $\nu\in S\setminus T,$ then
there are canonical pseudo-isomorphisms
 $$Cl_S(k_\infty)(p)^\vee\sim Cl_T(k_\infty)(p)^\vee,\hspace{1cm}
 \e_S(k_\infty)\sim\e_T(k_\infty).$$

\end{enumerate}
\end{prop}

\begin{proof}
Consider the canonical exact diagram for a finite extension $k'$
of $k$ in $k_\infty$
 $$\xymatrix@1@C=12pt{
 {\ E_T(k')\otimes_\mathbb{Z}\zp \ar@{^{(}->}[r]^{i_{k'}} } &  {\ E_S(k')\otimes_\mathbb{Z}\zp \ar[r] } &  {\ \bigoplus_{(S\setminus T)(k')}\zp \ar[r] } &  {\
    Cl_T(k')(p)\ar@{->>}[r]^{\pi_{k'}} } & {\ Cl_S(k')(p)}. \
 }$$
Setting $C(k'):=coker(i_{k'})$ (resp.\ $D(k'):=ker(\pi_{k'})$),
$C_\infty=\dirlim{} C(k')$ (resp.\ $D_\infty=\dirlim{} D(k')$) and
tensoring with $\qp/\zp$ we get the following exact sequences
 \begin{eqnarray*}
&& \xymatrix@1@C=12pt{
    {\ 0 \ar[r] } & {\ E_T(k')\otimes_\mathbb{Z}\qp/\zp \ar[r] } &  {\ E_S(k')\otimes_\mathbb{Z}\qp/\zp \ar[r] } &  {\ C(k')\otimes_\mathbb{Z}\qp/\zp \ar[r] } &  {\
    0, } \
 }\\
&&\xymatrix@1{
  {\  0 \ar[r]}  & {\ D(k') \ar[r] } &  {\ C(k')\otimes_\mathbb{Z}\qp/\zp \ar[r] } &  {\ \bigoplus_{(S\setminus T)(k')}\qp/\zp \ar[r] } &  {\ 0,  }
 }\\
&&\xymatrix@1{
   {\ 0 \ar[r] } &  {\ D(k') \ar[r] } &  {\ Cl_T(k')(p) \ar[r] } &  {\ Cl_S(k')(p) \ar[r] } &  {\
   0. }
}
 \end{eqnarray*}
Taking the direct limit over all finite subextensions $k',$ we get
an isomorphism $D_\infty\cong C_\infty\otimes\qp/\zp$ because the
transition maps for the sum of the $\qp/\zp$'s is just the
multiplication with the ramification index. The first result
follows after taking the Pontryagin dual.\\
 Now assume that $S\setminus T$ consists of a single prime and set $\bar{G}:=G(k'/k).$ Since
 then \linebreak $\bar{G}_\nu=G_\nu G(k_\infty/k')/G(k_\infty/k')$ acts
 trivial on $\bigoplus_{(S\setminus
 T)(k')}\zp\cong\Ind^{\bar{G}_\nu}_{\bar{G}}\zp$ and therefore also on $C(k')\otimes\qp/\zp,$ it
 follows that $G_\nu$ acts trivial on
 $(C_\infty\otimes\qp/\zp)^\vee.$ But then any surjection
 $\La^n\twoheadrightarrow (C_\infty\otimes\qp/\zp)^\vee$ factors
 through $(\Ind^{G_\nu}_G \zp)^n$ which is torsion (resp.\ pseudo-null) if $\dim
 (G_\nu)\geq 1 $ (resp.\ $\dim
 (G_\nu)\geq 2).$
The last statement is a consequence of the second one.
\end{proof}

The \La-modules $Cl_S(k_\infty)(p)^\vee$ and $\e_S(k_\infty)$ are
related to each other and to $X_S$ via Kummer theory:

\begin{prop}\label{ClS-eS}
Assume that $\mu_{p^\infty}\subseteq k_\infty.$ Then the following
holds:
\begin{enumerate}
\item There are  exact sequences of \La-modules

$$\xymatrix{
    {\ 0 \ar[r] } &  {\ Cl_S(k_\infty)(p)^\vee \ar[r] } &  {\ X_S(-1) \ar[r] } &  {\ \e_S(k_\infty) \ar[r] } &  {\
    0 } }$$
and, if $k_\infty$ is contained in $k_\Sigma,$ where
$\Sigma=S_p\cup S_\infty,$
 $$\xymatrix{{\ 0 \ar[r] } &  {\ Cl(k_\infty)(p)^\vee \ar[r] } &  {\ X_\Sigma(-1) \ar[r] } &  {\ \e(k_\infty) \ar[r] } &  {\
    0. } \
 }$$
In particular, $Cl_S(k_\infty)(p)^\vee$ and $Cl(k_\infty)(p)^\vee$
do not contain any pseudo-null submodule in these cases.
\item $Cl_S(k_\infty)(p)^\vee$ is
\La-torsion. If $\dim(G_\nu)\geq 1$ for every $\nu\in S_p,$ then
\linebreak $Cl(k_\infty)(p)^\vee$ is \La-torsion, too. In
particular, there are exact sequences
\begin{eqnarray*}
\xymatrix{
    {\ 0 \ar[r] } &  {\ Cl_S(k_\infty)(p)^\vee \ar[r] } &  {\ \tor_\La X_S(-1) \ar[r] } &  {\ \tor_\La\e_S(k_\infty) \ar[r] } &  {\
    0, } \\
    {\ 0 \ar[r] } &  {\ Cl(k_\infty)(p)^\vee \ar[r] } &  {\ \tor_\La X_\Sigma(-1) \ar[r] } &  {\ \tor_\La\e(k_\infty) \ar[r] } &  {\
    0. }
 }&&
\end{eqnarray*}
\end{enumerate}
\end{prop}

\begin{proof}
The long exact $\H^i(G(k_S/k_\infty),-)$ -sequence of
 $$\xymatrix@1@C=12pt{
    {\ 0 \ar[r] } &  {\ \mu_{p^n} \ar[r] } &  {\ E_S(k_S) \ar[r]^{p^n} } &  {\ E_S(k_S) \ar[r] } &  {\
    0 } \
 }$$
induces the short exact sequence
 $$\xymatrix@1@C=12pt{
   {\ 0 \ar[r] } & {\ E_S(k_\infty)/p^n \ar[r] } &    {\ \H^1(G(k_S/k_\infty),\mu_{p^\infty}) \ar[r] } &  {\ _{p^n}\H^1(G(k_S/k_\infty),E_S(k_s)) \ar[r] } &  {\
    0, } \
 }$$
i.e.\ after taking the direct limit with respect to $n$
 $$\xymatrix@1@C=12pt{
   {\ 0 \ar[r] } & {\ E_S(k_\infty)\otimes_\mathbb{Z}\qp/\zp \ar[r] } &    {\ \H^1(G(k_S/k_\infty),\qp/\zp)(1) \ar[r] } &  {\ Cl_S(k_\infty)(p) \ar[r] } &  {\
    0.} \
 }$$
Taking the dual, we obtain  the first statement. A canonical
map\linebreak $Cl(k_\infty(p)^\vee)\to X_S(-1)$ which is
compatible with the inclusion $Cl_S(k_\infty)(p)^\vee\to X_S(-1)$
from the first sequence can be defined exactly as in the
$\zp$-case, see \cite[11.4.2 and errata]{nsw}. Then the exactness
of the second sequence at the first term is obtained from the
first one and prop.\ \ref{ClS-eS}:
 $$ Cl(k_\infty)(p)^\vee/Cl_\Sigma(k_\infty)(p)^\vee\subseteq \e_\Sigma\cong
 X_\Sigma(-1)/Cl_\Sigma(k_\infty)(p)^\vee,$$
i.e.\ $Cl(k_\infty)(p)^\vee$ can be considered as submodule of
$X_\Sigma(1)$ and then its quotient is $\e.$\\ Comparing the ranks
of $X_S$ and $\e_S$ (see \ref{E0}) (with respect to an arbitrary
open subgroup $H\subseteq G$ such that $\La(H)$ is integral),  we
conclude that $Cl_S(k_\infty)(p)^\vee$ is \La-torsion while the
analogous result for $Cl(k_\infty)(p)^\vee$ follows from prop.\
\ref{ClS-eS}. Now, the last sequences can be derived from the
prior ones by rank considerations or by applying  the snake lemma
to the canonical sequence of homotopy theory
(\ref{canonicalsequ}).
\end{proof}

\begin{ques}
Is it true for any  $p$-adic Lie extension $k_\infty$ (of
dimension at least one) that $Cl(k_\infty)(p)^\vee$ and
$Cl_S(k_\infty)(p)^\vee$ don't have no non-zero pseudo-null
$\La$-submodules?
\end{ques}

In the $\zp$-case there exists a remarkable duality between the
inverse and direct limit of the ($S$-) ideal class groups in the
$\zp$-tower, viz the  pseudo-isomorphisms
\begin{eqnarray*}
X_{nr}&\sim&\E^1(Cl(k_\infty)(p)^\vee)\sim
(Cl(k_\infty)(p)^\vee)^\circ,\\
X_{cs}^S&\sim&\E^1(Cl_S(k_\infty)(p)^\vee)\sim
(Cl_S(k_\infty)(p)^\vee)^\circ.
\end{eqnarray*}
Therefore it seems natural (though maybe very optimistic) to pose
the following

\begin{ques}\label{quescl-xnr}
Is it true that for any $p$-adic Lie extension (at least under the
assumption ``$\dim(G_\nu)\geq 2,$ for all $\nu\in S_f,$") there
exist pseudo-isomorphisms \begin{eqnarray*}
X_{nr}&\sim&\E^1(Cl(k_\infty)(p)^\vee)\equiv
(Cl(k_\infty)(p)^\vee)^\circ \ \mbox{\rm  mod
 }{\mathcal{ PN}},\\
X_{cs}^S&\sim&\E^1(Cl_S(k_\infty)(p)^\vee)\equiv
(Cl_S(k_\infty)(p)^\vee)^\circ \ \mbox{\rm  mod
 }{\mathcal{ PN}}
\end{eqnarray*}
\end{ques}

 Observe, that $X_{nr}\sim X_{cs}^S$ and
$Cl_S(k_\infty)(p)^\vee\sim Cl(k_\infty)(p)^\vee$ by
\ref{tate-poitou}, \ref{ClSClT}. Hence, it would  suffice to prove
the existence of one of the pseudo-isomorphisms. By  prop.\
\ref{ClS-eS}(ii) and theorem \ref{XcspseudotorX} the question
would be true if one could show that the \La-torsion  of
$\e_S(k_\infty)$ is pseudo-null. But it seems difficult to prove
the latter statement directly. In fact, in the  case of  a
multiple $\zp$-extension  $k_\infty|k$ where
$\mu_{p^\infty}\subseteq k_\infty$ and $k$ has only one prime
above $p$ ,  W. McCallum \cite[thm 7]{mcCal} answers the above
question positively and then derives $\tor_\La\e_S(k_\infty)=0$
 just from the desired pseudo-isomorphism. This is the only case to the knowledge of the author where a
 positive answer to this question is known. Also J. Nekovar \cite[0.14.2]{nek} proved partial results in the
direction of the question. In a forthcoming paper \cite{ven3} we
will present the first non-abelian example (for
$G\cong\zp\rtimes\zp$ the semidirect product of two copies of \zp), in which such a duality holds.\\




The next result generalizes theorem 11.3.7 of \cite{nsw}.

\begin{thm}\label{free-prop}
Let $k_\infty|k$ be a  $p$-adic pro-$p$ Lie extension such that
$G$ is
without $p$-torsion.
\linebreak Then $\G=G(k_S(p)/k_\infty)$ is a free pro-$p$-group if
and only if $\mu(X_S)=0$ and the weak Leopoldt conjecture holds:
$\H^2(G_S(k_\infty),\qp/\zp)=0.$
\end{thm}

\begin{proof}
 Since $\G$ is pro-$p$ it is free if and
only if $\H^2(\G,\mathbb{Z}/p)=0,$ i.e.\ if and only if $_p(X_S)$
and $\H^2(G_S(k_\infty),\qp/\zp)$ vanish. But, by remark 3.33 of
\cite{ven1}  and since $X_S$ does not contain any pseudo-null
submodule, these two conditions are equivalent to the vanishing of
$\mu(X_S)$ and the validity of the weak Leopoldt conjecture.
\end{proof}

The next theorem, which generalizes theorem 11.3.8 in \cite{nsw},
shows that the validity of the weak Leopoldt conjecture and the
vanishing of the $\mu$-invariant  are properties which should be
considered simultaneously if one studies the behavior of $X_S$
under change of the base field.

\begin{thm}\label{base-change}
Let $K|k$ be a finite Galois $p$-extension inside $k_S,$
$k_\infty|k$ a $p$-adic pro-$p$ Lie extension such that
 \begin{center}
$(\ast)$ \ $G=G(k_\infty/k)$ is without $p$-torsion.
\end{center}

Set $K_\infty=Kk_\infty$ and $G'=G(K_\infty/K).$ Then $G'$
satisfies the condition $(\ast),$ too, and the following holds
 $$\left\{\begin{array}{c}
   \mu(X_S(k_\infty/k))=0\mbox{ and } \\
   \H^2(G_S(k_\infty),\qp/\zp)=0 \
 \end{array}\right\}\Leftrightarrow\left\{\begin{array}{c}
   \mu(X_S(K_\infty/K))=0\mbox{ and } \\
   \H^2(G_S(K_\infty),\qp/\zp)=0 \
 \end{array}\right\}.$$
In particular, if $ k_\infty$ contains the cyclotomic
$\zp$-extension, then
 $$\mu(X_S(k_\infty/k))=0\Leftrightarrow\mu(X_S(K_\infty/K))=0.$$
\end{thm}

\begin{proof}
Let $\h':=\h\cap G(k_S(p)/K).$ Then, by theorem \ref{free-prop},
the statements to be compared are equivalent to the freeness of
$\h,$ resp.\ $\h',$ thus equivalent to $\cd_p(\h)=1,$ resp.\
$\cd_p(\h')=1.$ But, since $\h'$ is open in $\h$ and
$\cd_p(\h)<\infty,$ we have $\cd_p(\h')=\cd_p(\h)$ by \cite{nsw}
3.3.5,(ii).
\end{proof}

The same arguments prove the following

\begin{thm}\label{relativmu}
Let $K_\infty|k_\infty|k$ be $p$-adic pro-$p$ Lie extensions
(inside $k_S$) such that for both $G(K_\infty/K)$ and
$G(k_\infty/k)$ the condition $(\ast)$ of the previous theorem
holds.  Then
 $$\left\{\begin{array}{c}
   \mu(X_S(k_\infty/k))=0\mbox{ and } \\
   \H^2(G_S(k_\infty),\qp/\zp)=0 \
 \end{array}\right\}\Rightarrow\left\{\begin{array}{c}
   \mu(X_S(K_\infty/k))=0\mbox{ and } \\
   \H^2(G_S(K_\infty),\qp/\zp)=0 \
 \end{array}\right\}.$$
\end{thm}

The next theorem, which generalizes theorem 11.3.5 in \cite{nsw},
describes the ``difference" if we vary the finite set of places
$S$ defining the module $X_S.$ By $T(K/k)\subseteq G(K/k)$ we
shall denote the inertia subgroup for a Galois extension $K|k$ of
local fields and, for an arbitrary set of places $S$ of $k$ and a
$p$-adic analytic extension $k_\infty|k,$ we write $S^{cd}(k)$ for
the subset of finite places which decompose completely in
$k_\infty|k$.

\begin{thm}\label{S-variation}
Let $S\supseteq T\supseteq S_p\cup S_\infty$ be finite
sets of places of $k$ and let $k_\infty|k$  be a  $p$-adic pro-$p$
Lie extension  inside $k_T$ with Galois group $G.$ Assume that $G$
does not contain any $p$-torsion element and that the weak
Leopoldt conjecture holds for $k_\infty|k.$ Then there exists a
canonical exact sequence of \La-modules
  $$\xymatrix{
    {0 \ar[r] } &  { \bigoplus_{(S\setminus
T)(k)}\Ind^{G_\nu}_G T(k_\nu(p)/k_\nu)_{G_{k_{\infty,\nu}}} \ar[r]
} &  {\ X_S \ar[r] } &  {\ X_T \ar[r] } &  {\
    0 } \
 }$$
and the direct sum on the  left is isomorphic to
 $$\bigoplus_{\stackrel{(S\setminus T)(k)}{p^\infty|f_\nu,\
 \mu_p\subseteq
 k_\nu}}\Ind^{G_\nu}_G\zp(1)\oplus\bigoplus_{(S\setminus
 T)^{cd}(k)}\La/p^{t_\nu},$$
 where $p^{t_\nu}=\#\mu(k_\nu)(p)$ and, as before,
$f_\nu=f(k_{\infty,\nu}/k_\nu)$ denotes the degree of the
extension of the corresponding residue class fields. In
particular, there is an exact sequence of \La-torsion modules
$$\xymatrix{
    {0 \ar[r] } &  { \bigoplus_{(S\setminus
T)(k)}\Ind^{G_\nu}_G T(k_\nu(p)/k_\nu)_{G_{k_{\infty,\nu}}} \ar[r]
} &  {\ \tor_\La X_S \ar[r] } &  {\ \tor_\La X_T \ar[r] } &  {\
    0 .} \
 }$$
\end{thm}

\begin{proof}
Since $\H^2(G_T(k_\infty)(p),\qp/\zp)=0,$ we have an exact
sequence
 $$\xymatrix@1{
    {\ 0 \ar[r] } &  {\ G(k_S(p)/k_T(p))^{ab}_{G_T(k_\infty)} \ar[r] } &  {\ X_S \ar[r] } &  {\ X_T \ar[r] } &  {\
    0 .} \
 }$$
Setting $\G=G_T(k)(p)$ and using \cite[10.5.4,10.6.1]{nsw} as well
as lemma \ref{kuz}, we obtain
 \begin{eqnarray*}
G(k_S(p)/k_T(p))^{ab}_{G_T(k_\infty)}&\cong&(\bigoplus_{(S\setminus
T)(k)} \Ind^{\G_\nu}_\G T(k_T(p)_\nu(p)/k_T(p)_\nu)\
)_{G_T(k_\infty)}\\
 &\cong&\bigoplus_{(S\setminus
T)(k)}\Ind^{G_\nu}_G T(k_\nu(p)/k_\nu)_{G_{k_{\infty,\nu}}}.
 \end{eqnarray*}
Observe that, for $\nu\in S\setminus T,$
 $$T(k_\nu(p)/k_\nu)\cong\left\{\begin{array}{cl}
  \mathbb{Z}_p(1) & \mbox{ if } \mu_p\subseteq k_\nu ,\\
  0 & \mbox{ otherwise. }
\end{array}\right.$$
Since $G$ is without $p$-torsion and $\nu\in S\setminus T$ is
unramified in $k_\infty|k,$ there are only two possibilities for
$G_\nu:$
 $$G_\nu=\left\{\begin{array}{cl}
   0 & \mbox{if }\nu\mbox{ is completely decomposed in } k_\infty|k, \\
   \mathbb{Z}_p & \mbox{if }p^\infty |f_\nu,
 \end{array}\right.$$ respectively,
 $$G_{k_\infty,\nu}(p)\cong\left\{\begin{array}{cl}
   \mathbb{Z}_p(1)\rtimes\zp, & \mbox{if }\nu\mbox{ is completely decomposed in } k_\infty|k, \\
   \mathbb{Z}_p(1), & \mbox{if }p^\infty |f_\nu.
 \end{array}\right.$$
It follows that
 $$G(k_S(p)/k_T(p))^{ab}_{G_T(k_\infty)}\cong\bigoplus_{\stackrel{(S\setminus T)(k)}{p^\infty|f_\nu,\
 \mu_p\subseteq
 k_\nu}}\Ind^{G_\nu}_G\zp(1)\oplus\bigoplus_{(S\setminus
 T)^{cd}(k)}\La/p^{t_\nu}.$$
In particular, this module is \La-torsion and therefore the second
statement follows from the first.
\end{proof}

Recalling that $\mu$ is additive on short exact sequences of
\La-torsion modules  we obtain the following

\begin{cor}\label{S-variationmu}
Under the assumptions of the theorem,
 $$\mu(X_S)=\mu(X_T) +\sum_{(S\setminus T)^{cd}(k)} t_\nu,$$
where  $ p^{t_\nu}=\#\mu(k_\nu)(p).$
\end{cor}
\POP
\PUSH{units.tex}%

\subsubsection{Global units} \label{globalunits}

 We still consider $p$-adic Lie extensions
$k_\infty|k$ with Galois group $G=G(k_\infty/k)$.

Recall that we denote the norm compatible $S$-units of $k_\infty$
by $\gus:=\projlim{k'} ({\mathcal{O}}_{k',S}^\times\otimes\zp).$
Noting that $\gus\cong\projlim{k'} \H^1(G_S(k'),\zp(1))$ by Kummer
theory and the finiteness of the $S$-ideal class group, we are
going to derive some relations between $\gus$ and
$\H^1(G_S(k_\infty),\mu_{p^\infty})^\vee$ by interpreting
Jannsen's spectral sequence (\cite{ja-sps}, see also \cite[Thm.\
4.5]{ven1}) or for Iwasawa adjoints with respect to
$A=\Qp/\zp(1)=\mu_{p^\infty}(k_S)$. We assume that $G$ does not
have any $p$-torsion, i.e.\ $G$ is a Poincar\'{e} group at $p,$ and we
denote the character which gives the operation of $G$ on the
dualizing module by $\chi^{-1}.$

\begin{prop}

\begin{enumerate}
\item If $\mu_{p^\infty}\subseteq k_\infty,$ then
\begin{enumerate}
\item if $\cd_p(G)=1:$
\begin{eqnarray*}
{\phantom{\projlim{k'} }}
\gus&\cong&\mathbb{Z}_p(1)(\chi)\oplus\E^0(X_S(-1))\\
\projlim{k'} H^2(G_S(k'),\zp(1))&\cong&\E^1(X_S(-1)),\\
\E^i(X_S(-1))&=&0\ \; \mbox{for }\ i\geq 2.
\end{eqnarray*}
\item if $\cd_p(G)=2:$ there is an exact sequence
\begin{eqnarray*}&\hspace{-3cm}\xymatrix@1{
   {\ 0 \ar[r] } &  {\ \gus \ar[r] } &  {\ \E^0(X_S(-1)) \ar[r] } &  {\zp(1)(\chi) \ar[r] } &
   }\\
   &\hspace{2cm}\xymatrix@1{ {\
   \projlim{k'}^{\phantom{\stackrel{M}{N}}}\hspace{-6pt} H^2(G_S(k'),\zp(1)) }\ar[r] &
   {\E^1(X_S(-1))\ar[r]}
   & 0,}
\end{eqnarray*}
 and $$\E^i(X_S(-1))=0\ \mbox{for }\ i\geq 2.$$
\item if $\cd_p(G)=3:$ there is an exact sequence
 $$\xymatrix@1@C=12pt{
    {\ 0 \ar[r] } &  {\ \projlim{k'}^{\phantom{\stackrel{M}{N}}}\hspace{-6pt} H^2(G_S(k'),\zp(1)) \ar[r] } &  {\E^1(X_S(-1)) \ar[r] } &  {\ \zp(1)(\chi) \ar[r] } &  {\
   0, }
 }$$
 and  \begin{eqnarray*}
\gus&\cong&\E^0(X_S(-1)),\\ \E^i(X_S(-1))&=&0\ \; \mbox{for }\
i\geq 2.
 \end{eqnarray*}
\item if $\cd_p(G)\geq 4:$
 \begin{eqnarray*}& \phantom{\projlim{k'}}\gus\cong\E^0(X_S(-1)),\\
&\projlim{k'} H^2(G_S(k'),\zp(1))\cong\E^1(X_S(-1)),\\
&\E^i(X_S(-1))=\left\{ \begin{array}{cl}
  \mathbb{Z}_p(1)(\chi) &    \mbox{if }\ i=\cd_p(G)-2, \\
  0 &  \mbox{otherwise,}
\end{array}\right.\ \mbox{for }\ i\geq 2.\end{eqnarray*}
\end{enumerate}

\noindent
 Similar results hold for arbitrary $A$ with
$k(A)\subseteq k_\infty$ if $\gus$ is replaced by \linebreak
$\projlim{k'} \H^1(G_S(k'),T_p A)$, $X_S(-1)$ by $X_S[A]$, ...

\item If $\mu(k_\infty)(p)$ is finite, then
\begin{enumerate}
\item if $\cd_p(G)=1:$ then there is an exact sequence
 {
 $$ \hspace{-0cm}\xymatrix@1@C=12pt{
   { 0 \ar[r] } &  { \gus \ar[r] } &  { \E^0(\H^1(G_S(k_\infty),\mu_{p^\infty})^\vee) \ar[r] } &  { \mu(k_\infty)(p)^\vee(\chi) \ar[r] } &  {
  \projlim{k'}^{\phantom{\stackrel{M}{N}}}\hspace{-6pt} H^2(G_S(k'),\zp(1)). }
 }$$ }
{\rm ($\mbox{a}_1$)} If in addition
$\H^2(G_S(k_\infty),\mu_{p^\infty})=0,$ then the cokernel of the
sequence is $\E^1(\H^1(G_S(k_\infty),\mu_{p^\infty})^\vee)$ and\\
 $\E^i(\H^1(G_S(k_\infty),\mu_{p^\infty})^\vee)=0\ \mbox{for }\ i\geq
 2.$\\

{\rm ($\mbox{a}_2$)} If  in addition
$\H^2(G_S(k_\infty),\Qp/\zp)=0,$ then there is a short exact
sequence
 $$ \hspace{-1.5cm}\xymatrix@1@C=16pt{
   { 0 \ar[r] } &  {\gus \ar[r] } &  { \E^0(\H^1(G_S(k_\infty),\mu_{p^\infty})^\vee) \ar[r] } &  { \mu(k_\infty)(p)^\vee(\chi) \ar[r] } &  {
  0. }
}$$

\item if $\cd_p(G)=2,$ then
$\gus\cong\E^0(\H^1(G_S(k_\infty),\mu_{p^\infty})^\vee).$\\

 If in
addition $\H^2(G_S(k_\infty),\mu_{p^\infty})=0,$ then there is an
exact sequence  {
$$\hspace{0cm}{ \xymatrix@1@C=12pt{
   { 0 \ar[r] } &  {\projlim{k'}^{\phantom{\stackrel{M}{N}}}\hspace{-6pt} H^2(G_S(k'),\zp(1)) \ar[r] } &  {\E^1(\H^1(G_S(k_\infty),\mu_{p^\infty})^\vee) \ar[r] } &  {\mu(k_\infty)(p)^\vee(\chi) \ar[r] } &  {
  0 }
}}$$} and $$\E^i(\H^1(G_S(k_\infty),\mu_{p^\infty})^\vee)=0\
\mbox{for }\ i\geq 2.$$
\item if $\cd_p(G)\geq 3,$ then
$\gus\cong\E^0(\H^1(G_S(k_\infty),\mu_{p^\infty})^\vee).$\\

 If in
addition $\H^2(G_S(k_\infty),\mu_{p^\infty})=0,$ then
\begin{eqnarray*}
\gus&\cong&\E^0(\H^1(G_S(k_\infty),\mu_{p^\infty})^\vee),\\
\projlim{k'}
H^2(G_S(k'),\zp(1))&\cong&\E^1(\H^1(G_S(k_\infty),\mu_{p^\infty})^\vee),\\
\end{eqnarray*}
{
$\E^i(\H^1(G_S(k_\infty),\mu_{p^\infty})^\vee)=\left\{
\begin{array}{cl}
 \mu(k_\infty)(p)(\chi) &    \mbox{if }\ i=\cd_p(G)-1, \\
  0 &  \mbox{otherwise},
\end{array}\right.\ \mbox{for }\ i\geq 2.$}
\end{enumerate}

\item If $\mu(k_\infty)(p)=0$, then there is in addition to the
results for finite $\mu(k_\infty)(p)$ the following exact
sequence: \begin{eqnarray*} &\hspace{-2.5cm}\xymatrix@1@C=12pt{
    {\ 0 \ar[r] } &  { \E^1(\H^1(G_S(k_\infty),\mu_{p^\infty})^\vee) \ar[r] } &  {\projlim{k'}^{\phantom{\stackrel{M}{N}}}\hspace{-6pt}
H^2(G_S(k'),\zp(1)) \ar[r] } & } \\
&\hspace{2cm}\xymatrix@1@C=12pt{
{\E^0(\H^2(G_S(k_\infty),\mu_{p^\infty})^\vee) \ar[r] } & {\
    \E^2(\H^1(G_S(k_\infty),\mu_{p^\infty})^\vee)\ar[r] }& 0, }
\end{eqnarray*}
and
$$\E^i(\H^1(G_S(k_\infty),\mu_{p^\infty})^\vee)\cong\E^{i-2}(\H^2(G_S(k_\infty),\mu_{p^\infty})^\vee)).$$
 \end{enumerate}
\end{prop}

For the proof apply Jannsen's  theorem (see \cite[Thm.\ 4.5,Cor.\
4.6]{ven1}) and note the following facts:
$\H^1(G_S(k_\infty),A)^\vee\cong X_S[A]$ if $k(A)\subseteq
k_\infty,$ $\H^2(G_S(k_\infty),A)=0$ if $\mu_{p^\infty}\subseteq
k_\infty$ because the weak Leopoldt conjecture is true for the
cyclotomic extension of any number field. Furthermore, we applied
several times \cite[2.6]{ja-is}. Also observe, that the reflexive
module $\E^0(X_S(-1))$ is projective in the case $\cd_p(G)=1$
regarding the defining sequence of the transpose functor $\du$ and
using that $\pd(\La)=cd_p(G)+1=2.$ The last statement of (ii)(a)
is proved in \cite{nsw} 11.3.9.

These results bear a lot of information about the structure of
 $\H^1(G_S(k_\infty),\mu_{p^\infty})^\vee$ and $\gus,$  e.g. one can
derive the projective dimension of the latter module (using
corollary \cite[Cor.\ 6.3]{ven1}) and some information about the
dimensions of the modules occurring above, in particular whether a
module is torsion(free). Furthermore, we see that $\gus$ is
reflexive for almost all cases with $\cd_p(G)\geq 2$ by
proposition 3.11 of \cite{ven1}.\\ In order to relate $\gus$ to
the finitely generated $\Lambda$-module
$$\e_S(k_\infty)=(E_S(k_ \infty)\otimes_\z\qp/\zp)^\vee$$ we need
some technical lemmas.

\begin{lem}
\begin{enumerate}
\item Let $G=G(k_\infty/k)\cong\mathbb{Z}_p^d,\ d\geq 1,$ and
$G_n:=p^nG.$
\begin{enumerate}
\item If $\mu_{p^\infty}\subseteq k_\infty,$ then with
$\Gamma=G(k(\mu_{p^\infty})$ and $\Gamma_n=p^n\Gamma$ the
following holds
 $$\H^i(G_n,\mu_{p^\infty})=\mu(k_n)(p)^{
   d-1 \choose
   i },$$
 where $k_n=k(\mu_{p^\infty})^{\Gamma_n}.$
 \item If $\mu(k_\infty)(p)$ is finite, then for any $n$ such that
 $\mu(k_\infty)(p)^{G_n}=\mu(k_\infty)(p)$ it holds
   $$\H^i(G_n,\mu_{p^\infty})=\mu(k_\infty)(p)^{
   d \choose
   i }.$$
\end{enumerate}
\item Let $G$ be a finitely generated pro-$p$ Lie group without $p$-torsion which fits into
a exact sequence $$\xymatrix@1{
   {\ 1 \ar[r] } &  {\ U \ar[r] } &  {\ G \ar[r]^\pi } &  {\ \Gamma \ar[r] } &  {\ 1
   },
}$$ with  $\Gamma\cong\zp$ and let  $G_n$ be an open subgroup.
Assume that $\Gamma_n:=\pi(G_n)$ acts via a splitting  trivially
on $U_n=G_n\cap U.$ Then $\H^2(G_n,\mu(k_\infty)(p))$ is finite
and the following holds
\begin{enumerate}
\item If $\mu_{p^\infty}\subseteq k_\infty$  and $\Gamma=G(k(\mu_{p^\infty}),$ then
 $$
 \H^1(G_n,\mu(k_\infty)(p))\cong\mu(k_n)(p)^s\oplus\bigoplus_i\mu_{p^{\nu_i}}(k_n),$$
where $U_n^{ab}\cong\zp^s\oplus\bigoplus_i \zp/p^{\nu_i}$ with
$U_n=U\cap G_n.$
\item If $\mu(k_\infty)(p)$ is finite, then for any $n$ such that
 $\mu(k_\infty)(p)^{G_n}=\mu(k_\infty)(p)$ there is an exact
 sequence
  $$\hspace{1cm}\xymatrix@1@C=12pt{
     { 0 \ar[r] } &  { \mu(k_\infty)(p) \ar[r] } &  {\H^1(G_n,\mu(k_\infty)(p)) \ar[r] } &  { \mu(k_\infty)(p)^s\oplus\bigoplus_i\mu_{p^{\nu_i}}(k_\infty) \ar[r] } &  {
     0. }
  }$$
\item If\ $\cd_p(G)=2,$ then $$\H^2(G_n,
  \mu(k_\infty)(p))\cong\left\{\begin{array}{cl}
    0 & \mbox{if } \mu_{p^\infty}\subseteq k_\infty, \\
    \mu(k_\infty)(p) & \mbox{otherwise.} \\
  \end{array}\right.$$
\end{enumerate}
\end{enumerate}
\end{lem}

\begin{proof}
Consider the exact sequence
 $$\xymatrix@1{
   {\ 1 \ar[r] } &  {\ U \ar[r] } &  {\ G \ar[r]^{\pi} } &  {\ \Gamma \ar[r] } &  {\ 1
   },
}$$ and let $U_n=G_n\cap U$ and $\Gamma_n=\pi (G_n)$. The
Hochschild-Serre spectral sequence gives  $$\xymatrix@1{
    { \H^1(\Gamma_n,\H^i(U_n,\mu(k_\infty)(p))) \ar@{^{(}->}[r] } &  { \H^{i+1}(G_n,\mu(k_\infty)(p)) \ar@{->>}[r] } &  { \H^{i+1}(U_n,\mu(k_\infty)(p))^{\Gamma_n}  }   \
 }$$ for $i\geq 0.$\\ Let us first assume that $\mu_{p^\infty}\subseteq
 k_\infty:$ Since $U_n$ acts trivially on $\mu_{p^\infty},$ we get
  $$\H^i(U_n,\Qp/\zp(1))=\H^i(U_n,\Qp/\zp)(1)=(\Qp/\zp)^{
   d-1 \choose
   i }.$$ in the abelian case by the K\"{u}nneth formula.
 As $\Qp/\zp(1)_{\Gamma_n}=0$ it follows that
 $\H^i(G_n,\mu_{p^\infty})=\H^i(U_n,\mu_{p^\infty})^{\Gamma_n}=\mu(k_n)^{
   d-1 \choose
   i }.$ In the non-abelian case we calculate
 \begin{eqnarray*}
 \H^1(G_n,\mu_{p^\infty})&=&\H^1(U_n,\Qp/\zp)(1)^{\Gamma_n}\\
                 &=&(U_n^{ab})^\vee(1)^{\Gamma_n}\\
                 &=&\mu(k_n)(p)^s\oplus\bigoplus\mu_{p^{\nu_i}}(k_n).
\end{eqnarray*}
Hence $ \H^1(\Gamma_n,\H^1(U_n,\mu(k_\infty)(p)))$ is finite and
the finiteness of $\H^2(G_n,\mu_{p^\infty})$ follows because
$\H^2(U_n,\mu_{p^\infty})^{\Gamma_n}\cong\H^2(U_n,\Qp/\zp)(1)^{\Gamma_n}$
is also finite  ($H^2(U_n,\Qp/\zp)$ is a cofinitely generated
abelian group).\\ Now we consider the case of finite
$\mu(k_\infty)(p):$ Here
$\H^1(\Gamma_n,\mu(k_\infty)(p))=\mu(k_\infty)(p)$ and the abelian
case follows again using the K\"{u}nneth formula. In the non-abelian
case the finiteness of $\H^2(G_n,\mu_{p^\infty})$ is trivial while
 $\H^1(U_n,\mu(k_\infty)(p))^{\Gamma_n}$ can be
calculated similarly as above. For the last assertion just note
that $U_n\cong\zp.$
\end{proof}

\noindent
\begin{lem}\label{finite submodules}\begin{enumerate}
\item In the situation of the previous lemma \mbox{(ii)} it holds
\hspace{-1.5cm}\begin{enumerate}
\item $\projlim{m,n}
{_{p^m}\H^1(G_n,E_S(k_\infty)/\mu(k_\infty))}\cong\projlim{m,n}
{_{p^m}\H^1(G_n, E_S(k_\infty))}=0,$

\item $\projlim{n}\H^1(G_n, E_S(k_\infty))\subseteq X_{cs}^S,$
\item $\begin{array}[t]{rcl}
\E^0(\e_S(k_\infty))&\cong&\projlim{m,,n}
{_{p^m}(E_S(k_\infty)\otimes \Qp/\zp)^{G_n}}\\       &\cong&
\projlim{m,n} (E_S(k_\infty)/\mu(k_\infty))^{G_n}/p^m,\end{array}$
\item $\t0 (\projlim{n}\H^1(G_n, E_S(k_\infty)/\mu(k_\infty)))=\t0
(\E^1((E_S(k_\infty)\otimes_\mathbb{Z}\Q_p/\zp)^\vee)),$
\item that the following sequence is exact:{\footnotesize
$$\hspace{0.5cm}\xymatrix@1@C=12pt{
   {\ 0 \ar[r] } &  {\projlim{n}^{\phantom{\stackrel{M}{N}}}\hspace{-6pt}
\H^1(G_n, E_S(k_\infty)/\mu(k_\infty))} \ar[r]  &
{\E^1(\e_S(k_\infty)) \ar[r] } &
{\projlim{m,n}^{\phantom{\stackrel{M}{N}}}\hspace{-6pt}
{_{p^m}\H^2(G_n, E_S(k_\infty)/\mu(k_\infty))} \ar[r] } &  {\
   0. }
}$$}
\end{enumerate}
\item If, in addition, $\cd_p(G)\leq 2,$ then with  $\kappa=\left\{\begin{array}{cl}
  1 & \mbox{if } \mu(k_\infty)(p)\mbox{ is finite,} \\
  0 & \mbox{otherwise,e}
\end{array}\right.$
there are the following exact sequences
\hspace{0cm}\begin{enumerate}
\item if $\cd_p(G)=2:$
{\footnotesize
\begin{eqnarray*}
&&\hspace{-0cm}\xymatrix@1@C=12pt{
    { 0 \ar[r] } &  {\projlim{n}^{\phantom{\stackrel{M}{N}}}\hspace{-6pt}
\H^1(G_n, E_S(k_\infty)) \ar[r] } &
{\projlim{m,n}^{\phantom{\stackrel{M}{N}}}\hspace{-6pt}
\H^1(G_n,E_S(k_\infty)/\mu(k_\infty))/p^m \ar[r] } &
{\mu(k_\infty)(p)^\kappa \ar[r] } & D\ar[r]&0,   \
 }\\
 &&\hspace{-0cm}\xymatrix@1@C=12pt{
    { 0 \ar[r] } &  {\projlim{m,n}^{\phantom{\stackrel{M}{N}}}\hspace{-6pt} {_{p^m}\H^2(G_n, E_S(k_\infty))} \ar[r] } &  {\projlim{m,n}^{\phantom{\stackrel{M}{N}}}\hspace{-6pt} {_{p^m}\H^2(G_n, E_S(k_\infty)/\mu(k_\infty)) }\ar[r] } &  {D\ar[r]} &  } \\
  &&\hspace{2cm}\xymatrix@1@C=12pt{ {\projlim{m,n}^{\phantom{\stackrel{M}{N}}}\hspace{-6pt} {\H^2(G_n, E_S(k_\infty))/{p^m}} \ar[r]
    }& {\projlim{m,n}^{\phantom{\stackrel{M}{N}}}\hspace{-6pt} {\H^2(G_n, E_S(k_\infty)/\mu(k_\infty))/{p^m} }\ar[r] }
    &0, }
\end{eqnarray*}}
where $D$ is some finite module.
\item if $\cd_p(G)=1:$
{\footnotesize
$$\hspace{-0.2cm}\xymatrix@1@C=12pt{
    { 0 \ar[r] } & {\E^2\E^1(\gus)\ar[r]}& {\mu(k_\infty)(p)^\kappa  }\ar[r] &{\projlim{n}^{\phantom{\stackrel{M}{N}}}\hspace{-6pt}
\H^1(G_n, E_S(k_\infty)) \ar[r] } &
{\projlim{m,n}^{\phantom{\stackrel{M}{N}}}\hspace{-6pt}
\H^1(G_n,E_S(k_\infty)/\mu(k_\infty))/p^m \ar[r] } & 0   \
 }$$}
 and $$\projlim{m,n} {_{p^m}\H^2(G_n,
E_S(k_\infty))}\cong\projlim{m,n} {_{p^m}\H^2(G_n,
E_S(k_\infty)/\mu(k_\infty)) }.$$
\end{enumerate}
\end{enumerate}
\end{lem}

\begin{proof}
If we split the long exact cohomology sequence  induced by
 $$\xymatrix@1{
    {\ 0 \ar[r] } &  {\ \mu(k_\infty)  \ar[r] } &  {\ E_S(k_\infty) \ar[r] } &  {\ E_S(k_\infty)/\mu(k_\infty) \ar[r] } &  {\
    0, } \
 }$$
we get the following short exact sequences
 \begin{eqnarray*}
 && \xymatrix@1{
    {\ 0 \ar[r] } &  {\ F_n \ar[r] } &  {\ \H^1(G_n, \mu(k_\infty)) \ar[r] } &  {\ A_n \ar[r] } &  {\
    0, } \
 }\\
&& \xymatrix@1{
    {\ 0 \ar[r] } &  {\ A_n \ar[r] } &  {\ \H^1(G_n, E_S(k_\infty)) \ar[r] } &  {\ B_n \ar[r] } &  {\
    0, } \
 }\\
&& \xymatrix@1{
    {\ 0 \ar[r] } &  {\ B_n \ar[r] } &  {\ \H^1(G_n,E_S(k_\infty)/\mu(k_\infty))  \ar[r] } &  {\ C_n \ar[r] } &  {\
    0 } \
 }
 \end{eqnarray*}
and furthermore a map
$\xymatrix{C_n\ar@{^{(}->}[r]&{\H^2(G_n,\mu(k_\infty)(p)).}}$
Evaluating the associated long exact sequences  of
 $p^m$-torsion (snake lemma) and noting the finiteness of $A_n$
 and  $C_n$ according to the previous lemma, we get
  \begin{eqnarray*}
 && \projlim{m} {_{p^m}B_n}\cong \projlim{m}
  {_{p^m}\H^1(G_n,E_S(k_\infty)/\mu(k_\infty)),} \\
 && \xymatrix@1{
     {\ 0 \ar[r] } &  {\ \projlim{m}^{\phantom{\stackrel{M}{N}}}\hspace{-6pt} {_{p^m}\H^1(G_n,E_S(k_\infty))}  \ar[r] } &  {\ \projlim{m}^{\phantom{\stackrel{M}{N}}}\hspace{-6pt} {_{p^m}B_n} \ar[r] } &  {\ A_n,  }  \\
  }
  \end{eqnarray*}
  and therefore
  $$\xymatrix@1@C=12pt{
   {\ 0 \ar[r] } &  {\ \projlim{m,n}^{\phantom{\stackrel{M}{N}}}\hspace{-6pt}{_{p^m}\H^1(G_n,E_S(k_\infty))} \ar[r] } &  {\ \projlim{m,n}^{\phantom{\stackrel{M}{N}}}\hspace{-6pt}
 { _{p^m}\H^1(G_n,E_S(k_\infty)/\mu(k_\infty))} \ar[r] } &  {\ \projlim{n} A_n
  }
}$$ is exact.

But $\projlim{n} A_n$ is a quotient of $$\projlim{n}
\H^1(G_n,\mu(k_\infty)(p))=\left\{\begin{array}{cl}
  \mu(k_\infty)(p) & \mbox{if } d=1 \mbox{ and } \mu(k_\infty)(p)\mbox{ is finite, }  \\
 0 & \mbox{otherwise.}
\end{array} \right.$$ (See the previous lemma and note that the
transition maps are  partially norm maps besides the non-trivial
case where they are  the natural projections, i.e.\ identities for
$n$ sufficiently big.). Since the middle term is $\zp$-torsion
free, we get the desired isomorphism, because, by the
Hochschild-Serre spectral sequence, it can be seen in any case
that the first group is contained in $\projlim{m,n}
{_{p^m}Cl_S(k_n)=0.}$ This proves (i)(a) while (b) is again the
cited spectral sequence.

The first equality of (i)(c) is just theorem 4.7 (iii) of
\cite{ven1} because $\e_S(k_\infty)$  has no $\zp$-torsion while
the second one follows by the exact sequence
{
$$\hspace{0cm}\xymatrix@1{
      {\ (E_S(k_\infty)/\mu(k_\infty))^{G_n}/p^m  \ar@{^{(}->}[r] } &  {\ _{p^m}(E_S(k_\infty)\otimes_\mathbb{Z}\Q_p/\zp)^{G_n} \ar@{->>}[r] } &  {\ _{p^m}\H^1(G_n,E_S(k_\infty)/\mu(k_\infty))  }  \
 }$$}
and (a). Similar arguments apply for (i)(e), i.e.\
$$\E^1(\e_S(k_\infty))\cong\projlim{m,n}\H^1(G_n,{_{p^m}(E_S(k_\infty)\otimes_\mathbb{Z}\Q_p/\zp)}).$$
The assertion (d) is a direct consequence of (e), because
\linebreak $\projlim{m,n} {_{p^m}\H^2(G_n,
E_S(k_\infty)/\mu(k_\infty))}$ is $\zp$-torsion-free.

Now let us assume that $\cd_p(G)\leq 2.$ With the  notation as
above and recalling that $A_n,B_n$ and $C_n$ are finite, we get
exact sequences
 $$\xymatrix@1{
    {\ 0 \ar[r] } &  {\ A_n \ar[r] } &  {\ \H^1(G_n, E_S(k_\infty)) \ar[r] } &  {\ B_n \ar[r] } &  {\
    0, } \
 }$$

 $$\xymatrix@1{
    {\ 0 \ar[r] } &  { B_n \ar[r] } &  {\projlim{m}\H^1(G_n,E_S(k_\infty)/\mu(k_\infty))/p^m \ar[r] } &  {\ C_n \ar[r] } &  {\
    0 } \
 }$$
 and
 $$\xymatrix@1{
    {\ 0 \ar[r] } &  {\ C_n \ar[r] } &  {\ \H^2(G_n, \mu(k_\infty)) \ar[r] } &  {\ D_n \ar[r] } &  {\
    0. } \
 }$$

Passing to the limit gives the first exact sequence in (ii)(a)
(Note that the transition maps of the system $\{C_n\}$ are the
canonical projections, i.e.\ identities for $n$ sufficiently
large). The second one is proved similarly using
 $$\xymatrix@1@C=12pt{
      {\ D_n \ar@{^(->}[r] } &  {\H^2(G_n,E_S(k_\infty)) \ar[r] } &  { \H^2(G_n,E_S(k_\infty)/\mu(k_\infty)) \ar[r] } &  {
    \H^3(G_n,\mu(k_\infty)(p))=0 }
 }$$
and $\xymatrix{{\H^2(G_n,\mu(k_\infty)(p))\ar@{->>}[r]}&D_n}.$ The
proof of (ii)(b) is completely analogous, just note that $\projlim
n F_n\cong\E^2\E^1(\gus)$ because the latter module is the
cokernel of $\gus\to\E^0\E^0(\gus)\cong\E^0(\e_S(k_\infty)).$
\end{proof}

\begin{prop}\label{gus-structure}
There is an exact sequence { $$\xymatrix@1{
    {\  0 \ar[r] } &  {\ \zp(1)^\delta \ar[r] } &  {\ \gus \ar[r] } &  {\ \E^0(\e_S(k_\infty)) \ar[r] } &  {\
    C }
 }$$}
with  $$\ C=\left\{\begin{array}{ll}
 \mu(k_\infty)(p) & \mbox{if } d=1\mbox{ and }\mu(k_\infty)(p)\mbox{ finite} \\
 \mathbb{Z}_p(1) & \mbox{if } d=2\mbox{ and }\mu_{p^\infty}\subseteq k_\infty\\
 f.g.\ \zp\mbox{-module} & d\geq 3\mbox{ and  $G$ non-abelian}\\
  0 & otherwise
\end{array}\right.$$ and  $$\ \delta=\left\{\begin{array}{cl}
  1 & \mbox{if } d=1,\mu_{p^\infty}\subseteq k_\infty, \\
  0 & \mbox{otherwise.}
\end{array}\right.$$
If in addition the weak Leopoldt conjecture holds, the right map
 is onto  in the case  $d=1$ and  $\mu(k_\infty)(p)$ finite.
\end{prop}

\begin{proof}
Taking $G_n$-invariants of the exact sequence
 $$\xymatrix@1{
    {\ 0 \ar[r] } &  {\ \mu(k_\infty)(p) \ar[r] } &  {\ E_S(k_\infty)\otimes_\mathbb{Z}\zp \ar[r] } &  {\ (E_S(k_\infty)/\mu(k_\infty))\otimes_\mathbb{Z}\zp \ar[r] } &  {\
    0 } \
 }$$
 and passing to the inverse limit, we get
{\footnotesize $$\xymatrix@1@C=12pt{
    {\ 0 \ar[r] } &  {\ \projlim{n}^{\phantom{\stackrel{M}{N}}}\hspace{-6pt}\mu(k_n)(p) \ar[r] } &  {\ \gus \ar[r] } &  {\ \projlim{m,n}^{\phantom{\stackrel{M}{N}}}\hspace{-6pt} (E_S(k_\infty)/\mu(k_\infty))^{G_n}/p^m  \ar[r] } &  {\
    \projlim{n}^{\phantom{\stackrel{M}{N}}}\hspace{-6pt}
\H^1(G_n,\mu(k_\infty)(p))} \
 }$$} The result follows except  the fact that $\E^0$ maps onto the finite group of roots of unity in the case when
 $d=1.$ But this is proved in \cite{nsw} 11.3.9 under the
 assumption that
 the weak Leopoldt conjecture holds.
\end{proof}

\begin{cor}\label{E0}Let $k_\infty|k$ be a $p$-adic Lie extension
such that $G$ does not have any $p$-torsion. Then
$$\E^0(\gus)\cong\E^0\E^0(\e_S(k_\infty))\cong\E^0(\H^1(G_S(k_\infty),\mu_{p^\infty})^\vee).$$
 In particular, if $G$ is in addition pro-$p$ and
 $\H^2(G_S(k_\infty),\mu_{p^\infty})=0$ (e.g. if $\mu_{p^\infty}\subseteq
 k_\infty$), then $$\rk_\La\gus=\rk_\La\e_S=r_2(k).$$
\end{cor}

Now the question arises whether the  module $\E^0(\gus)$ is not
only reflexive but also projective. While in the case $\cd_p(G)=1$
this is always true, in  higher dimensions one needs additional
conditions. We will only get a satisfying answer in the two
dimensional case:

\begin{prop}\label{E0projectiv}
Let $k_\infty|k$ be a $p$-adic Lie extension such that
$\cd_p(G)=2$ and assume that the weak Leopoldt conjecture holds
for $k_\infty.$ Then the following is equivalent:
\begin{enumerate}
\item $\E^0(\gus)$ is projective,
\item  $\t0
 \E^1(\e_S(k_\infty))=0.$
\end{enumerate}

\end{prop}

\begin{rem}
These equivalent statements hold for example, if either
$\mu_{p^\infty}\subseteq k_\infty$ or $\mu(k_\infty)(p)=0$, and
$\t0(X_{cs}^S)=0,$ i.e.\ if $X_{cs}^S$ does not have any non-zero
finite submodule, because then $\t0
 \E^1(\e_S(k_\infty))=0$ by  lemma \ref{finite
submodules}.\end{rem}

\begin{proof}
Since we already know that $\pd(\E^0(\gus))\leq 1,$ because
$\E^0(\gus)$ is the  second syzygy of $\du\gus$, the projectivity
is equivalent to the vanishing of $\E^1\E^0(\gus)$. Now the
equivalence stated above follows from the next lemma.
\end{proof}

\begin{lem} In the situation of the proposition it holds
 $$\t0
 \E^1(\e_S(k_\infty))\cong\E^1\E^0(\gus)\cong\E^3\E^1(\gus)$$

\end{lem}

\begin{proof}
Set $M:=\e_S(k_\infty)$ and consider the exact sequence
 $$\xymatrix@1{
    {\ 0 \ar[r] } &  {\ M/\t1 (M) \ar[r] } &  {\ \E^0(\gus) \ar[r] } &  {\ \E^2D(M) \ar[r] } &  {\
    0. } \
 }$$
The long exact sequence for $\E^i$ gives
 $$\xymatrix@1{
    {\ 0=\E^1\E^2D(M) \ar[r] } &  {\ \E^1\E^0(M)  \ar[r] } &  {\ \E^1(M/\t 1(M)) \ar[r] } &  {\ \E^2\E^2D(M) . }    \
 }$$
On the other hand there is the exact sequence
 $$\xymatrix@1{
    {\ 0=\E^0(\t1 (M) ) \ar[r] } &  {\ \E^1(M/\t1 (M) ) \ar[r] } &  {\ \E^1(M) \ar[r] } &  {\ \E^1\E^1D(M).  }    \
 }$$
Since $\E^i\E^iD(M)$ is pure of codimension $i$, the isomorphism
follows. But \linebreak $\E^1\E^0(\gus)\cong\E^3\E^1(\gus)$ by the
spectral sequence due to Bj\"{o}rk, see \ref{reflexivepd3}.\end{proof}

The proposition above should be compared with  the following
result which has already been observed by Kay Wingberg
(unpublished):

\begin{prop}
If $\cd_p(G)=1,$ then for sufficiently large $n$ there is a
canonical exact sequence
 $$\xymatrix@1{
    {\ 0 \ar[r] } &  {{\e_S(k_\infty)}^{G_n} \ar[r] } &  {\e_S(k_\infty) \ar[r] } &  {\E^0(\gus) \ar[r] } &  {\
    C\ar[r] } &0 \
 }$$
where $C=\E^2\mathrm{D}(\e_S(k_\infty))$ is connected with
$\E^2\mathrm{D}(\gus)$ by the exact sequence
 $$\xymatrix@1{
    {\ 0  \ar[r] } &  {\ \E^2\mathrm{D}(\gus) \ar[r] } &  {\mu(k_\infty)(p)^\kappa \ar[r] } &  {\ \t0 X_{cs}^S \ar[r] } &  {\
    C^\vee \ar[r]} & 0. \
 }$$
\end{prop}

\begin{proof}The first sequence is just the canonical sequence
\ref{canonicalsequ} for the module $\e_S(k_\infty)$ while the
second one already occurred in lemma \ref{finite submodules}
(ii)(b) as we show now: The fact that $\t0
(X_{cs}^S)\cong\projlimssc{n}\H^1(G_n, E_S(k_\infty))$ is well
known (see for example \cite[XI.\S3.]{nsw}). Recall that
$\E^2\E^1(\e_S(k_\infty))\cong\t0 (\E^1(\e_S(k_\infty)))^\vee$ and
apply lemma \ref{finite submodules} (i)(d) to recover $C$. Using
\ref{finite submodules}, (i)(e) and (ii)(b) we see that
$\E^1\E^1(\e_S(k_\infty))\cong\E^1(\projlim{m,n} {_{p^m}\H^2(G_n,
E_S(k_\infty))}),$ which we will determine by means of \cite[Thm.\
4.7 (iii)]{ven1}:
$$M:=\projlim{m,n} {_{p^m}\H^2(G_n, E_S(k_\infty))}^\vee\\  \cong
 \dirlim{m,n} {\e_S(k_\infty)}^{G_n}/p^m\\
  = \dirlim{m}
{\e_S(k_\infty)}^{G_n}/p^m$$
for $n$ sufficiently large, because $\e_S(k_\infty)$ is a finitely
generated $\Lambda$-module. Hence $$ \E^1(M)\cong  \projlim{m,n}
({_{p^m}M})_{G_n}={\e_S(k_\infty)}^{G_n}.$$
 for $n$ large enough.
\end{proof}

\begin{prop}\label{E0gus-structure}
Let $k_\infty/k$ be a $p$-adic Lie extension such that
$G\cong\Gamma\times\Delta$ where $\Gamma$ is a pro-$p$-Lie group
of $\cd_p(\Gamma)=2$, $\Delta$ is a finite group of order prime to
$p.$ Assume that the weak Leopoldt conjecture holds for
$k_\infty.$ Then the following is true:
\begin{enumerate}\item There is an exact sequence
\hspace{-1cm}\begin{eqnarray*}\hspace{-1cm}
&\hspace{-2cm}\xymatrix@1{
   {\ 0  \ar[r] } &  {\E^0\E^0(\gus)   \ar[r] } & {\Lambda^{r_2+r_1-r_1'-s}\oplus\bigoplus_{S^{cd}\cup S'_\infty}\Ind^{G_\nu}_G(\zp)\ar[r]}& }\\
   & \hspace{6.5cm}\xymatrix@1{ {\ \Lambda^s \ar[r] } &  {\ \t0 \E^1(\e_S(k_\infty)) \ar[r] } &  {\
   0. }}
\end{eqnarray*}
\item If $\E^0\E^0(\gus)$ is projective, then
 $$\E^0\E^0(\gus)\cong\Lambda^{r_2+r_1-r_1'}\oplus\bigoplus_{S^{cd}\cup S'_\infty}\Ind^{G_\nu}_G(\zp).$$
\end{enumerate}
\end{prop}

\begin{proof}
We calculate the Euler characteristic with respect to an arbitrary
open normal subgroup $U\unlhd\Gamma$ using lemma \ref{euler},
proposition \ref{diagram-prop},\cite{ja-is} 5.4 b),
 \begin{eqnarray*}
 h_U(\E^0\E^0(\gus))&=&h_U(\gus)\\
                    &=&h_U(\gas)-h_U(X_S)+h_U(X_{cs}^S)\\
                    &=&h_U(\gas)-h_U(Y_S)+h_U(I_G)\\
                    &=&h_U(\gas)-h_U(\Lambda^d)+h_U(N^{ab}_{\h}(p))+h_U(\Lambda)-h_U(\zp)\\
                    &=&h_U(\gas)-h_U(\Lambda^{r_2+r_1'})+h_U(\bigoplus_{S'_\infty}\Ind^{G_\nu}_G(\zp))\\
                    &=&\sum_S\Ind^{G_\nu}_G h_{U\cap G_\nu}(\la_\nu)-h_U(\Lambda^{r_2+r_1'})+h_U(\bigoplus_{S'_\infty}\Ind^{G_\nu}_G(\zp))\\
                    &=&\sum_{S^{cd}}\Ind^{G_\nu}_G
                    h_U(\zp)+h_U(\Lambda^{r_2+r_1-r_1'})+h_U(\bigoplus_{S'_\infty}\Ind^{G_\nu}_G(\zp)).
 \end{eqnarray*}
Therefore, if $\E^0\E^0(\gus)$ is projective, it follows that
  $$\E^0\E^0(\gus)\cong\Lambda^{r_2+r_1-r_1'}\oplus\bigoplus_{S^{cd}\cup S'_\infty}\Ind^{G_\nu}_G(\zp).$$
  This proves (ii) while (i) follows easily applying proposition
  \ref{reflexivepd3}.
\end{proof}
\POP
\PUSH{selmer.tex}%

\subsection{Selmer groups of abelian varieties}
\label{selmergroups}

 In this section let $k$ be a number field,
\A\ a $g$-dimensional abelian variety defined over $k$ and  $p$ a
fixed rational {\em odd} prime number. For a non-empty, finite set
$S$ of places of $k$ containing the places $S_{bad}$  of bad
reduction of \A, the places $S_p$ lying over $p$ and the places
$S_\infty$ at infinity we write $\H^i(G_S(k),\A)$, respectively
$\H^i(k_\nu,\A)$, for the cohomology groups
$\H^i(G_S(k),\A(k_S))$, respectively $\H^i(G_\nu,
\A(\bar{k}_\nu))$,  where $G_S(k)$ denotes the Galois group of the
maximal outside $S$ unramified extension of $k$, $\bar{k}_\nu$ the
algebraic closure of the completion of $k$ at $\nu$ and $G_\nu$
the corresponding decomposition group. The ($p^m)$-)Selmer group
$\Sel(\A,k,p^m)$ and the Tate-Shafarevich group $\ts(\A,k,p^m)$ %
fit by definition into the following commutative exact diagram

$$\xymatrix{
          &        & 0\ar[d] & 0 \ar[d] &   \\
  0 \ar[r]& {\A(k)/p^m}\ar@{=}[d]\ar[r] &{\Sel(\A,k,p^m)}\ar[d]\ar[r] & {\ts(\A,k,p^m)}\ar[d]\ar[r] & 0 \\
  0 \ar[r] & {\A(k)/p^m}\ar[r] & {\H^1(G_S(k),{_{p^m}\A})} \ar[d]\ar[r]& {_{p^m}\H^1(G_S(k),\A)}\ar[d]\ar[r] & 0 \\
          &  & {\bigoplus_{S(k)}^{\phantom{S}}\H^1(k_\nu,\A)(p)}\ar@{=}[r] & {\bigoplus_{S(k)}^{\phantom{S}}\H^1(k_\nu,\A)(p)}. &
}$$

If $k_\infty$ is an infinite Galois extension of $k$  with Galois
group  $G=G(k_\infty/k,)$ we get the following commutative exact
diagram by passing to the direct limit with respect to $m$ and
finite subextensions $k'$ of $k_\infty/k$:

{\footnotesize $$\hspace{-0cm}\xymatrix@C=8pt{
          &        & 0\ar[d] & 0 \ar[d] &   \\
  0 \ar[r]& {\A(k_\infty)\otimes\Qp/\zp}\ar@{=}[d]\ar[r] &{\Sel(\A,k_\infty,p^\infty)}\ar[d]\ar[r] & {\ts(\A,k_\infty,p^\infty)}\ar[d]\ar[r] & 0 \\
  0 \ar[r] & {\A(k_\infty)\otimes \Qp/\zp}\ar[r] & {\H^1(G_S(k_\infty),\A(p))} \ar[d]\ar[r]& {\H^1(G_S(k_\infty),\A)}(p)\ar[d]\ar[r] & 0 \\
          &  &{\bigoplus_{S(k)}
\Coind^{G_\nu}_G \H^1(k_{\infty,\nu},\A)(p)}\ar@{=}@<0.4ex>[r] &
{\bigoplus_{S(k)} \Coind^{G_\nu}_G \H^1(k_{\infty,\nu},\A)(p)}. &
 } $$}

Note that $\dirlim{k'}
\bigoplus_{S(k')}\H^1(k'_\nu,\A)(p)\cong\bigoplus_{S(k)}
\Coind^{G_\nu}_G \H^1(k_{\infty,\nu},\A)(p)$. Alternatively, we
can pass to the inverse limits and we will get the following
commutative exact diagram $$\xymatrix{
          &         & 0\ar[d] & 0 \ar[d] &   \\
  0 \ar[r]& {\widehat{\A}_{k_\infty}}\ar@{=}[d]\ar[r] &{\csel (k_\infty,\A)}\ar[d]\ar[r] & {\projlim{k',m}^{\phantom{\stackrel{M}{N}}}\hspace{-6pt}\ts(\A,k',p^m)}\ar[d]\ar[r] & 0 \\
  0 \ar[r] &  {\widehat{\A}_{k_\infty}}\ar[r] & {\projlim{k'}^{\phantom{\stackrel{M}{N}}}\hspace{-6pt}\H^1(G_S(k'),T_p \A)} \ar[d]\ar[r]& {\projlim{k'}^{\phantom{\stackrel{M}{N}}}\hspace{-6pt}T_p\H^1(G_S(k'),\A)}\ar[d]\ar[r] & 0 \\
          &  &{\projlim{k'} \bigoplus_{S(k')}T_p\H^1(k'_\nu,\A)}\ar@{=}@<1ex>[r] & {\projlim{k'} \bigoplus_{S(k')}T_p\H^1(k'_\nu,\A).} &
} $$ where $ \widehat{\A}_{k_\infty}:=\projlim{k',m} \A(k')/p^m$
and $\csel (k_\infty,\A):=\projlim{k',m} \Sel (k',\A,p^m)$ (The
limits are taken with respect to corestriction  maps and
multiplication by $p$).

Henceforth we will drop   the $p$ from the notation of the Selmer
group: $$\Sel(\A,k_\infty):=\Sel(\A,k_\infty,p^\infty).$$
Furthermore, we shall use the following notation for the
local-global modules
 \begin{eqnarray*}
 \lu_{S,\A}&:=&\bigoplus_{S_f(k)}
 \Ind^{G_\nu}_G\H^1(k_{\infty,\nu},\A)(p)^\vee,
 \\
 \la_{S,\A}&:=&\bigoplus_{S_f(k)} \Ind^{G_\nu}_G
 \H^1(k_{\infty,\nu},\A(p))^\vee,\\
 \tab_{S,\A}&:=&\bigoplus_{S_f(k)} \Ind^{G_\nu}_G
 (\A(k_{\infty,\nu})\otimes\Qp/\zp)^\vee.
 \end{eqnarray*}

As a consequence of the long exact  sequence of the Tate-Poitou
duality theorem we have the following (compact) analogue of
proposition \ref{tate-poitou}, where we shall write $\A^d$ for the
dual abelian variety of $\A$ and $\Sha^1_S(k_\infty,\A(p))$ for
the kernel of the localization map
$$\H^1(G_S(k_\infty),\A(p))\to\bigoplus_{S(k)} \Coind^{G_\nu}_G
\H^1(k_{\infty,\nu},\A(p)).$$

\begin{prop}\label{selmerTate-Poitou}Let $k_\infty|k$ be a $p$-adic Lie extension with Galois group $G.$ Then,
there are the following exact commutative diagrams of
$\La=\La(G)$-modules
\begin{enumerate}
\item

 $$ \xymatrix{
 0& 0\\
{\Sel(\A,k_\infty)^\vee}\ar@{->>}[r]\ar[u] &
{\Sha^1_S(k_\infty,\A(p))^\vee\ar[u]}\\
{\H^1(G_S(k_\infty),\A(p))^\vee\ar@{=}[r]\ar[u]}&
{\H^1(G_S(k_\infty),\A(p))^\vee\ar[u]}\\
{\lu_{S,\A}\ar@{^{(}->}[r]\ar[u]}&
{\la_{S,\A}\ar@{->>}[r]\ar[u]}&{\tab_{S,\A} } \\ {\csel
(k_\infty,\A^d)\ar@{^{(}->}[r]\ar[u]
}&{\projlimsc{k'}^{\phantom{\stackrel{M}{N}}}\hspace{-6pt}\H^1(G_S(k'),T_p
(\A^d))} \ar[r]\ar[u]&{\tab_{S,\A} ,\ar@{=}[u]}\\
{\H^2(G_S(k_\infty),\A(p))^\vee\ar[u]\ar@{=}[r]}&{\H^2(G_S(k_\infty),\A(p))^\vee\ar[u]}\\
 0\ar[u]& 0\ar[u]}$$

\item
\begin{eqnarray*}
&\hspace{-2.3cm}\xymatrix@1@C-0pt{
   {\ 0 \ar[r] } &  {\ \csel (k_\infty,\A^d) \ar[r] } &  { \projlimsc{k'}^{\phantom{\stackrel{M}{N}}}\hspace{-6pt}\H^1(G_S(k'),T_p \A^d) \ar[r] } &  { \tab_{S,\A} \ar[r] } & {\phantom{M}}}\\
   &\hspace{6.5cm}\xymatrix@1{  {\  \Sel(\A,k_\infty)^\vee \ar[r] }&  {\  \Sha^1_S(k_\infty,\A(p))^\vee\ar[r] }&  {0, }}
\end{eqnarray*}
\item
{\small $$ \hspace{0.6cm}\xymatrix@C=10pt{
   {0 \ar[r] } &  {\Sha_S^1(k_\infty,\A(p))^\vee \ar[r] } &  { Z_{S,\A^d(p)} \ar[r] } &  {\bigoplus_{S_f(k)} \Ind^{G_\nu}_G  (\A(k_{\infty,\nu})(p))^\vee\ar[r] } &  {
   \A(k_\infty)(p)^\vee\ar[r] } & 0.
}$$}
\end{enumerate}
\end{prop}

For the proof, just note that by virtue of local Tate duality
(\cite[Cor.3.4]{milne}),
 the Weil-pairing and \ref{globalZ},
\begin{enumerate}
\item $\H^1(k_{\infty,\nu},\A)(p)^\vee\cong\widehat{(\A^d)}_{\infty,\nu}:=\projlim{k',m}\A^d(k'_\nu)/p^m,$
\item  $Z_{S,\A^d(p)}\cong\projlimsc{k'}\H^2(G_S(k'),T_p(\A^d)),$
\item $(\A(k_{\infty,\nu})\otimes\Qp/\zp)^\vee\cong\projlim{k'}T_p\H^1(k'_\nu,\A^d) $ and
\item $\H^1(k_{\infty,\nu},\A(p))^\vee\cong\projlimsc{k'}\H^1(k'_\nu,T_p(\A^d))$
\end{enumerate}
hold.\\

By a well-known theorem  of Mattuck, we have an isomorphism
$$\A(k'_\nu)\cong\mathbb{Z}_l^{g[k'_\nu :\mathbb{Q}_l]}\times
\mbox{(a finite group)},$$ for any finite extension $k'_\nu$ of
$\mathbb{Q}_l$. Recall that $g$ denotes the dimension of the
abelian variety $\A.$ Clearly
 $$\A(k'_\nu)\otimes_\mathbb{Z} \Qp/\zp=0$$
for all $l\neq p$ and $\nu\mid l,$ i.e.\
$$\H^1(k'_\nu,\A)(p)\cong\H^1(k'_\nu,\A(p)),$$ respectively
$$\H^1(k'_{\infty,\nu},\A)(p)\cong\H^1(k'_{\infty,\nu},\A(p)),$$
in this case. On the other hand, Coates and Greenberg proved that
for primes $\nu\mid p$ with good reduction
 $$\H^1(k_{\infty,\nu},\A)(p)\cong\H^1(k_{\infty,\nu},\widetilde{\A}(p))$$
holds, if $k_\infty$ is a deeply ramified, where $\widetilde{\A}$
denotes the reduction of $\A$ (see \cite[Prop.\
 4.8]{co-gr}). We recall that an algebraic extension  $k$ of
 $\qp$ is called {\em deeply ramified}
if $\H^1(k,\overline{\m})$ vanishes,   where $\overline{\m}$ is
the maximal ideal of the ring of integers of an algebraic closure
$\overline{\mathbb{Q}_p}$ of $\qp;$ see \cite[p.\ 143]{co-gr} for
equivalent conditions and for the following statement (loc.\ cit.\
thm.\ 2.13): A field $k_\infty$ which is a $p$-adic Lie extension
of a finite extension $k$ of $\qp$ is deeply ramified if the
inertial subgroup of $G(k_\infty/k)$ is infinite.\\ For arbitrary
reduction at $\nu\mid p$, the same
 result as above holds, if one replaces $\widetilde{\A}_{p^\infty}$ by the
 quotient $\A(p)/{\mathcal{F}}_\A(\overline{\m})(p)$, where
  ${\mathcal{F}}_\A$ denotes
 the formal group associated with the
 Neron model of $\A$ over a possibly  finite extension of
 $k_\nu$, such that the Neron model has semi-stable reduction.
Taking these facts into account, we get the following description
for $\lu_{S,\A},$ where $T(k_{\infty,\nu}/k_\nu)$ denotes the
inertia subgroup of $G_\nu.$

\begin{prop}\label{reducedUSA} (cf.\ {\rm \cite[lemma 5.4]{ochi-ven}})
Assume that $\dim(T(k_{\infty,\nu}/k_\nu)\geq 1$ for all $\nu\in
S_p.$ Then there is an isomorphism of \La-modules
 $$\lu_{S,\A}\cong\bigoplus_{S_p(k)}
 \Ind^{G_\nu}_G\H^1(k_{\infty,\nu},\widetilde{\A}(p))^\vee\oplus\bigoplus_{S_f\setminus S_p(k)}
 \Ind^{G_\nu}_G\H^1(k_{\infty,\nu},\A(p))^\vee.$$
In particular, if $\dim(G_\nu)\geq2$ for all $\nu\in S_f,$ then
 $$\lu_{S,\A}\cong\bigoplus_{S_p(k)}
 \Ind^{G_\nu}_G\H^1(k_{\infty,\nu},\widetilde{\A}(p))^\vee$$
and $\lu_{S,\A}$ is \La-torsion-free.
\end{prop}

\begin{proof}
The first assertion has been explained above while the second
statement follows from the local calculations \ref{local-d2},
\ref{localcd>2} with respect to the $p$-adic representations
$A=\widetilde{\A}(p)$ respectively $A=\A(p)$ and the comment
before \ref{localcd>2}.
\end{proof}

Before going on we would like to recall some well-known facts
about abelian varieties:

\begin{rem}
\label{schneider} (i) $\rk_\zp (\A(p)^\vee)= 2g,$ where $g$
denotes the dimension of $\A.$\\ (ii) There exists always an
isogeny from $\A$ to its dual $\A^d,$ by which the Weil-pairing
induces a non-degenerate skew-symmetric pairing on the Tate-module
$T_p\A$ of $\A,$ (combine \cite[cor.\ 7.2,lem.\ 16.2(e),prop.\
16.6]{mil-ab}). If $\A=E$ is an elliptic curve this isogeny can be
chosen as a canonical isomorphism between $E$ and $E^d.$ Again for
an arbitrary abelian variety it follows that
$k(\mu_{p^\infty})\subseteq k(\A(p))=k(\A^d(p))$ (see \cite[\S 0
lem. 7]{schneider}).
\end{rem}

\begin{thm}\label{billot}
Assume that $\H^2(G_S(k_\infty),(\A^d)(p))=0.$ If
 $\dim(G_\nu)\geq 2$ for all $\nu\in S_f,$ then
 $$\Sha^1_S(k_\infty,\A(p))^\vee\sim\E^1(Y_{S,\A^d(p)})\sim\E^1(\tor_\La Y_{S,\A^d(p)})\cong\E^1(\tor_\La X_{S,\A^d(p)}).$$
If, in addition, $G\cong \mathbb{Z}^r_p,$ $r\geq2,$ then  the
following holds:
 $$\Sha^1_S(k_\infty,\A(p))^\vee\sim (\tor_\La X_{S,\A^d(p)})^\circ,$$
 where $^\circ$ means that the $G$ acts via the involution
 $g\mapsto g^{-1}.$
 \end{thm}

\begin{rem}\label{rem}
In case $\tor_\La X_{S,\A^d(p)}$ is isomorphic in $\mod/\mathcal{
PN}$ to a direct sum of cyclic modules of the form $\La$ modulo a
(left) principal ideal  proposition \ref{circ} implies that
\[\Sha^1_S(k_\infty,\A(p))^\vee\equiv (\tor_\La X_{S,\A^d(p)})^\circ \ \mbox{\rm  mod
 }{\mathcal{ PN}}\] holds under the conditions of the theorem.
\end{rem}

\begin{proof}
The first condition implies
$Z_{S,\A^d(p)}\cong\E^1(Y_{S,\A^d(p)})$ while the other condition
grants that $\bigoplus_{S_f(k)} \Ind^{G_\nu}_G
(\A(k_{\infty,\nu})(p))^\vee$ is pseudo-null because
$\A(k_{\infty,\nu})(p)^\vee$ is a finitely generated (free)
$\zp$-module. Now everything follows as in \ref{XcspseudotorX}
using here prop.\ \ref{selmerTate-Poitou}.
\end{proof}

\begin{cor} \label{supersingular-ab-var} Let $\A$ be an abelian variety over $k$ with good  {\em
supersingular} reduction, i.e.\ $\widetilde{\A_{k_\nu}}(p)=0,$ at
all places $\nu$ dividing $p.$  Set $k_\infty=k(\A(p))$ and assume
that $G(k_\infty/k)$ is a pro-$p$-group without any $p$-torsion.
Then, for $\Sigma_{bad}:=S_{bad}\cup S_p\cup S_\infty$ the
following holds:
 $$X_{cs}[\A^d(p)]\cong\Sha^1_{\Sigma_{bad}}(k_\infty,\A^d(p))^\vee\sim\E^1(\tor_\La(\Sel(\A,k_\infty)^\vee)).$$
In particular, if  $\A$ has CM, then there is even a
pseudo-isomorphism
$$X_{cs}[\A^d(p)]\sim(\tor_\La(\Sel(\A,k_\infty)^\vee))^\circ.$$
\end{cor}

Therewith, in the case of an elliptic curve with CM,  we reobtain
a theorem of P. Billot \cite[3.23]{billot}. Over a $\zp$-extension
an analogous statement was proved by K. Wingberg \cite[cor.
2.5]{wing}. Of course, remark \ref{rem} applies literally to
$\tor_\La\Sel(\A,k_\infty)^\vee,$ i.e.\ under the conditions
mentioned there it holds
 $$X_{cs}[\A^d(p)]\equiv(\tor_\La\Sel(\A,k_\infty)^\vee)^\circ \ \mbox{\rm  mod
 }{\mathcal{ PN}}.$$

\begin{proof} First note that by the N\'{e}ron-Ogg-Shafarevich
criterion the sets of bad reduction of $\A$ and its dual $\A^d$
coincide. Therefore, it suffices to prove that $\dim(G_\nu)\geq 2$
for all $\nu\in S_{bad}\cup S_p$ because then the theorem applies
to $\A^d$ and proposition \ref{reducedUSA} shows that
$\lu_{S,\A}=0,$ i.e.\ $X_{S,\A(p)}\cong\Sel(\A,k_\infty)^\vee.$ \\
 So, let $\nu$ be either in $S_p$ or in $S_{bad}.$ Since
$k_\nu(\A(p))$ contains $k_\nu(\mu_{p^\infty}),$ we only have to
show that $G(k_\nu(\A(p))/k_\nu(\mu_{p^\infty}))$ is not trivial
because then it automatically has to be infinite as
$G_\nu\subseteq G$ has no finite subgroup by assumption.\\ If
$\nu|p,$ by a theorem of Imai\footnote{ I owe to John Coates the
idea to use Imai's theorem here.} \cite{imai}
$\A(k_\nu(\mu_{p^\infty}))(p)$ is finite and thus
$k_\nu(\A(p))\neq k_\nu(\mu_{p^\infty}).$\\ If $\nu\in S_{bad},$
then the N\'{e}ron-Ogg-Shafarevich criterion implies that
$G(k_\nu(\A(p))/k_\nu(\mu_{p^\infty}))=T(k_\nu(\A(p))/k_\nu)$ is
non-trivial.
\end{proof}


By  remark \ref{weakLeopold} and \ref{schneider} the conditions of
theorem \ref{nospeudonullglobalXA} are fulfilled for the
$p$-torsion points $\A(p)$ and its  trivializing extension of $k,$
i.e.\ the extension which is obtained by adjoining the $p$-torsion
points of $\A:$

\begin{thm}\label{nospeudonullglobalXAbvariety} Let $k_\infty=k(\A(p))$ and assume that $G$ does
 not have any $p$-torsion. Then
$\H^1(G_S(k_\infty), \A(p))^\vee$ has no non-zero pseudo-null
submodule.
\end{thm}

Recall that $G$ does not have any $p$-torsion if $p\geq
2\dim(\A)+2.$ Otherwise one only has to replace $k$ by a finite
extension inside $k_\infty.$

We should mention that the rank of the global module
$\H^1(G_S(k_\infty),\A(p))^\vee$ is $g[k:\Q],$ which was
determined by Y. Ochi   who also calculated the ranks  and
torsion-submodules of the local, respectively local-global modules
(i.e. those global modules which are induced from local ones) that
occur in proposition \ref{selmerTate-Poitou} (cf.
\cite[5.7,5.11,5.12]{ochi}). See also the results in S. Howson's
PhD-thesis \cite[5.30,6.1,6.5-6.9,6.13-6.14,7.3]{howson}.

Furthermore, in the case of elliptic curves S. Howson proved the
following result.

\begin{prop}(Howson {\rm \cite[6.14-15]{howson}})
Let $E$ be an elliptic curve over $k$ without complex
multiplication and with good {\em ordinary} reduction at all
places over $p.$ Assume that\/ $G=G(k(E(p)/k)$ is  pro-$p$ without
any $p$-torsion. Then
 $$\tab_{S,E}\cong\la_{S,\widetilde{E}}\cong\bigoplus_{S_f(k)}
 \Ind^{G_\nu}_G\projlimsc{k'}\H^1(k'_\nu,T_p(\widetilde{E}))$$
and these modules are $\La(G)$-torsion-free. Furthermore, there is
an isomorphism
 $$\lu_{S,E}\cong\E^0(\tab_{S,E}).$$
\end{prop}
\POP

\INPUT{xbib.bib}   
\INPUT{iwasawa.bbl} 

\bibliographystyle{amsplain}
\bibliography{xbib}
\end{document}